\documentclass[twoside,10pt,english]{amsart}

\usepackage[matrix,arrow]{xy}
\usepackage{amssymb,amsmath}
\usepackage[mathcal]{euscript}
\usepackage{mathrsfs}
\usepackage{url}
\title{Mixed Weil cohomologies}

\author{Denis-Charles Cisinski}
\address{LAGA\\
CNRS~(UMR 7539)\\
Institut~Galil\'ee\\
Universit\'e Paris~13\\
Avenue Jean-Baptiste Cl\'ement\\
93430~Villetaneuse\\France}
\email{cisinski@math.univ-paris13.fr}
\urladdr{http://www.math.univ-paris13.fr/~cisinski/}

\author{Fr\'ed\'eric D\'eglise}
\address{LAGA\\
CNRS (UMR 7539)\\
Institut Galil\'ee\\
Universit\'e Paris~13\\
Avenue Jean-Baptiste Cl\'ement\\
93430~Villetaneuse\\France}
\email{deglise@math.univ-paris13.fr}
\urladdr{http://www.math.univ-paris13.fr/~deglise/}
\thanks{Partially supported by the ANR (grant No. ANR-07-BLAN-042)}

\newtheorem{thmm}{Theorem}
\newtheorem{thm}{Theorem}[subsection]
\newtheorem{prop}[thm]{Proposition}
\newtheorem{lm}[thm]{Lemma}
\newtheorem{cor}[thm]{Corollary}

\theoremstyle{remark} 
\newtheorem{rem}[thm]{Remark}

\newtheorem{ex}[thm]{Example}

\theoremstyle{definition} 
\newtheorem{df}[thm]{Definition}
\newtheorem*{dff}{Definition}

\newtheorem{paragr}[thm]{}
\newtheorem{sch}[thm]{Scholium}
\numberwithin{equation}{thm}

\DeclareFontFamily{OT1}{pzc}{}
\DeclareFontShape{OT1}{pzc}{m}{it}%
             {<-> s * [1,150] pzcmi7t}{}
\DeclareMathAlphabet{\mathpzc}{OT1}{pzc}%
                                 {m}{it}
\input cyracc.def
\DeclareFontFamily{U}{russian}{}
\DeclareFontShape{U}{russian}{m}{n}
        { <5><6> wncyr5
        <7><8><9> wncyr7
        <10><10.95><12><14.4><17.28><20.74><24.88> wncyr10 }{}
\DeclareSymbolFont{Russian}{U}{russian}{m}{n}
\DeclareSymbolFontAlphabet{\mathcyr}{Russian}
\makeatletter
\let\@math@cyr\mathcyr
\renewcommand{\mathcyr}[1]{\@math@cyr{\cyracc #1}}
\makeatother

\renewcommand{\leq}{\leqslant}
\renewcommand{\geq}{\geqslant}
\renewcommand{\mathscr}{\mathcal}

\newcommand{\odmte}{\otimes^{\derL}}
\newcommand{\odmt}{\otimes^\derL}

\newcommand{\adams}[1]{\mathbf{\Psi}^{#1}}
\newcommand{\specialisation}{\mathit{sp}}
\newcommand{\motcoh}{\mathbf{H}_{\mathcyr{B}}}
\newcommand{\hypercoh}{H}

\newcommand{\DM}{\mathit{DM}}
\newcommand{\DMt}{\smash{\Der_{\AA^1}}}

\newcommand{\DMte}{\smash{\Der^{\mathit{eff}}_{\AA^1}}}

\newcommand{\DMtp}{\smash{\Der^{\vee}_{\AA^1}}}
\newcommand{\DMB}{\mathit{DM}_{\mathcyr{B}}}

\newcommand{\BGLQ}{\mathit{KGL}_{\QQ}}
\newcommand{\HQ}{\mathbf{HQ}}

\newcommand{\Sm}{\mathit S\mspace{-2.mu}m}

\newcommand{\tot} {\mathrm{Tot}}

\newcommand{\sm}{\Sm}

\newcommand{\Sh}{\operatorname{\mathit{Sh}}}

\newcommand{\pic}{\mathrm{Pic}}
\newcommand{\chow}{\mathit{CH}}
\newcommand{\tR}{\underline{R}}
\newcommand{\tKK}{\underline{\KK}}
\newcommand{\tQQ}{\underline{\QQ}}

\newcommand{\TDer} {\mathsf T}
\newcommand{\ste}{\mathcal E}
\newcommand{\tate}{\mathit{Tate}}
\newcommand{\Spt}{\smash{\mathrm{Sp}^{}_{\tate}}}
\newcommand{\thom}{\mathit{Th}\,}
\newcommand{\bundle}[1]{\mathscr #1}

\newcommand{\KK}{\mathbf K}
\newcommand{\dlog}{d\log} 
\newcommand{\dR}{\mathit{dR}}
\newcommand{\rigid} {\mathit{rig}}
\newcommand{\rig} {\rigid}
\newcommand{\an} {\mathit{an}}
\newcommand{\mw} {\mathit{MW}}
\newcommand{\idealmax}{\mathfrak{m}}

\renewcommand{\varinjlim}{\limind}%

\def\limproj{\mathop{\oalign{\rm lim\cr
\hidewidth$\longleftarrow$\hidewidth\cr}}}%
\def\limind{\mathop{\oalign{\rm lim\cr
\hidewidth$\longrightarrow$\hidewidth\cr}}}%
\def\TO#1{\mathrel{\hbox to #1pt{\rightarrowfill}}}
\def\OT#1{\mathrel{\hbox to #1pt{\leftarrowfill}}}
\def\hocolim{\mathop{\oalign{\rm holim\cr
\hidewidth$\TO{26.5}$\hidewidth\cr}}}%
\newcommand{\Mod}{{\mathrm{\text{-}\mathrm Mod}}}
\newcommand{\derL}{\mathbf{L}}
\newcommand{\derR}{\mathbf{R}}
\newcommand{\sHom}{{\mathbf{Hom}}}

\newcommand{\op}[1]{{ #1 }^{\mathit{op}}}

\newcommand{\Hom}{{\mathrm{Hom}}}

\newcommand{\Comp}{\mathit{Comp}}
\newcommand{\Htp} {\mathit K}
\newcommand{\Der} {\mathit D}

\newcommand{\spec}[1] { \mathrm{Spec}\left( #1 \right) }

\newcommand{\dual}[1]{#1^{\vee}}
\newcommand{\unit}{\mathbf 1}

\newcommand{\du}{\mathit{dual}}

\newcommand{\To}{\longrightarrow}

\newcommand{\e}{\varepsilon}
\newcommand{\s}{\sigma}
\newcommand{\A}{\mathscr A}

\newcommand{\C}{\mathscr C}
\newcommand{\D}{\mathscr D}

\newcommand{\V}{\mathfrak V}
\newcommand{\T}{\mathscr T}

\newcommand{\X}{\mathscr X}

\newcommand{\G}{\mathcal G}
\renewcommand{\H}{\mathcal H}

\newcommand{\ZZ} {\mathbf Z}
\newcommand{\QQ} {\mathbf Q}

\renewcommand{\AA} { \mathbf A }
\newcommand{\PP} {\mathbf P}
\newcommand{\GG} { \mathbf G }

\newcommand{\HH} {H}

\newcommand{\zar}{{\mathrm{Zar}}}
\newcommand{\nis}{{\mathrm{Nis}}}

\newcommand{\etale}{{\mathrm{\acute{e}t}}}

\begin{document}

\begin{abstract}
We define, for a regular scheme $S$ and a given field
of characteristic zero $\KK$, the notion of $\KK$-linear
mixed Weil cohomology on smooth
$S$-schemes by a simple set of properties, mainly:
Nisnevich descent, homotopy invariance, stability
(which means that the cohomology of $\GG_{m}$
behaves correctly), and K\"unneth formula. We prove
that any mixed Weil cohomology defined on smooth $S$-schemes induces
a symmetric monoidal realization of some suitable triangulated
category of motives over $S$ to the derived category of the field
$\KK$. This implies a finiteness theorem and a Poincar\'e duality
theorem for such a cohomology with respect to smooth and projective
$S$-schemes (which can be extended to smooth $S$-schemes
when $S$ is the spectrum of a perfect field). This formalism
also provides a convenient tool to understand the comparison
of such cohomology theories.
\end{abstract}

\maketitle
\setcounter{tocdepth}{3}
\tableofcontents
\section*{Introduction}

Weil cohomologies were introduced
by Grothendieck in the 1960's as the cohomologies defined
on smooth and projective varieties over a field with
enough good properties (mainly, the existence of cycle class maps,
K\"unneth formula and Poincar\'e duality)
to prove the Weil conjectures (i.e. to understand
the $L$-functions attached to smooth and projective varieties
over a finite field). According to the philosophy of Grothendieck,
they can be seen as the fiber functors of the (conjecturaly) tannakian
category of pure motives. From this point of view,
a mixed Weil cohomology should define an exact tensor
functor from the (conjectural) abelian category
of mixed motives to the category of (super) vector spaces
over  a field of characteristic zero, such that, among other things,
its restriction to pure motives would be a Weil cohomology.

The purpose of these notes is to provide a simple
set of axioms for a cohomology theory to induce a
symmetric monoidal realization functor of a suitable
version of the triangulated category of mixed motives
to the derived category of vector spaces over
a field of characteristic zero.
Such a compatibility with symmetric monoidal
structures involves obviously a K\"unneth
formula for our given cohomology. And the main
result we get here says that this property is
essentially sufficient to get a realization functor.
Moreover, apart from the K\"unneth formula,
our set of axioms is very close to that of Eilenberg and Steenrod 
in algebraic topology.

Let $k$ be a perfect field and $\KK$ a field
of characteristic zero. Let $\V$ be the
category of smooth affine $k$-schemes.
Consider a presheaf of commutative differential graded
$\KK$-algebras $E$ on $\V$.
Given any smooth affine scheme $X$, 
any closed subset $Z \subset X$ such that $U=X-Z$ is affine, 
and any integer $n$,
we put:
$$H^n_Z(X,E)=H^{n-1}\big(\mathrm{Cone}(E(X) \rightarrow E(U))\big)$$
\begin{dff}
A $\KK$-linear \emph{mixed Weil theory} is 
 a presheaf of differential graded
  $\KK$-algebras $E$ on $\V$ 
   satisfying the following properties:
\begin{enumerate}
\item[] \textit{Dimension}.--- 
$\mathrm{dim}_{\KK} H^{i}(\spec k,E)=
\begin{cases}
1&\text{if $i=0$,}\\
0&\text{otherwise;}
\end{cases}$
\item[] \textit{Homotopy}.---
$\mathrm{dim}_{\KK} H^i(\AA^1_{k},E)=
\begin{cases}
1&\text{if $i=0$,}\\
0&\text{otherwise;}
\end{cases}$
\item[] \textit{Stability}.---
$\mathrm{dim}_{\KK} H^i(\GG_{m},E)=
\begin{cases}
1&\text{if $i=0$ or $i=1$,}\\
0&\text{otherwise;}
\end{cases}$
\item[] \textit{Excision}.--- Consider a commutative diagram of $k$-schemes
$$
\xymatrix{
T\ar^j[r]\ar_g[d] & Y\ar^f[d] \\
Z\ar^i[r] & X
}
$$
such that $i$ and $j$ are closed immersions,
the schemes $X,Y,X-Z,Y-T$ are smooth and affine,
  $f$ is \'etale
 and $f^{-1}(X-Z)=Y-T$, $g$ is an isomorphism. Then the induced morphism
$$
H^*_T(Y,E) \To H^*_Z(X,E)
$$
is an isomorphism;
\item[] \textit{K\"unneth formula}.---
For any smooth affine $k$-schemes $X$ and $Y$, the exterior cup product
induces an isomorphism
$$
\bigoplus_{p+q=n}H^p(X,E) \otimes_\KK H^q(Y,E)
 \overset{\sim}{\To} H^{n}(X\times_{k}Y,E) \ .
$$
\end{enumerate}
\end{dff}
The easiest example of a mixed Weil theory the reader might
enjoy to have in mind
is algebraic de Rham cohomology over a field of characteristic zero.
The homotopy axiom in this setting is rather called the 
\emph{Poincar\'e lemma}.

We will prove that the excision axiom on a presheaf of differential
graded $\KK$-algebras $E$ is equivalent to the following property~: \\
\textit{Nisnevich Descent}.--- For any smooth affine scheme $X$,
the cohomology groups of the complex $E(X)$ are isomorphic 
to the Nisnevich hypercohomology groups of $X$ with coefficients
in $E_{\nis}$ under the canonical map.

Given a mixed Weil theory $E$
and any smooth scheme $X$, we denote by $H^n(X,E)$ the Nisnevich 
hypercohomology groups of $X$ with coefficients in $E_{\nis}$.
According to the previous assertion,
this extends the definition given above to the case where $X$ is affine.
We define for a $\KK$-vector space $V$
and an integer $n$
$$V(n)=\begin{cases}
V\otimes_{\KK}\Hom_{\KK}(H^1(\GG_{m},E)^{\otimes n},\KK)&\text{if $n\geq 0$,}\\
V\otimes_{\KK}H^1(\GG_{m},E)^{\otimes (-n)}&\text{if $n\leq 0$.}
\end{cases}$$

Note that any choice of a generator of $H^1(\GG_{m},E)$ defines an isomorphism
$V(n)\simeq V$.
The introduction of these Tate twists allows us to make canonical
constructions, avoiding the choice of a generator for $H^1(\GG_{m},E)$.
Our main results can now be summarized as follows.

\begin{thmm}\label{thm1}
The cohomology groups $H^q(X,E)$ have the following properties.
\begin{itemize}
\item[1.]\emph{Finiteness}.--- For any smooth $k$-scheme $X$,
the $\KK$-vector space
$\oplus_{n}H^n(X,E)$ is finite dimensional.
\item[2.]\emph{Cycle class map}.--- For any smooth $k$-scheme $X$,
there is a natural map which is compatible with cup product
$$H^q(X,\QQ(p))\To H^q(X,E)(p)$$
(where $H^q(X,\QQ(p))$ is motivic cohomology, as
defined by Voevodsky; in particular, we have $H^{2n}(X,\QQ(n))=\chow^n(X)_{\QQ}$).
\item[3.]\emph{Compact support}.---  For any smooth $k$-scheme $X$,
there are cohomology groups $H^q_{c}(X,E)$, which satisfies all
the usual functorialities of a cohomology with compact support,
and there are natural maps
$$H^q_{c}(X,E)\To H^q(X,E)$$
which are isomorphisms whenever $X$ is projective.
\item[4.]\emph{Poincar\'e duality}.--- For any smooth $k$-scheme $X$
of pure dimension $d$, there is a natural perfect pairing
of finite dimensional $\KK$-vector spaces
$$\HH^q_{c}(X,E)(p)\otimes^{}_{\KK}
\HH^{2d-q}_{}(X,E)(d-p)\To\KK\, .$$
\end{itemize}
\end{thmm}

\begin{thmm}[Comparison]\label{thm2}
Let $E'$ be a presheaf of commutative differential graded $\KK$-algebras
satisfying the dimension, homotopy, stability and excision axioms 
and such that for any smooth $k$-scheme $X$, the K\"unneth map
$$
\bigoplus_{p+q=n}H^p(X,E') \otimes_\KK H^q(Y,E')
 \To H^{n}(X\times_{k}Y,E')
$$
is an isomorphism for $Y=\AA^1_{k}$ or $Y=\GG_{m}$ (e.g. $E'$ might be a mixed Weil theory).
Then a morphism of presheaves of differential graded $\KK$-algebras
$E\To E'$ is a quasi-isomorphism (locally for the Nisnevich topology)
if and only if the map $H^1(\GG_{m},E)\To H^1(\GG_{m},E')$
is not trivial\footnote{The main point here is in fact that
the map $H^1(\GG_{m},E)\To H^1(\GG_{m},E')$ controls
the compatibility with cycle class maps, which in turns
ensures the compatibility with Poincar\'e duality.}.
\end{thmm}
Remark this comparison theorem is completely similar to the classical comparison
 theorem of Eilenberg-Steenrod.

\begin{thmm}[Realization]\label{thm3}
There is a symmetric monoidal triangulated functor
$$R:\DM_{\mathit{gm}}(k)_{\QQ}\To\Der^b(\KK)$$
(where $\DM_{\mathit{gm}}(k)_{\QQ}$
is Voevodsky's triangulated category of mixed motives over $k$,
and $\Der^b(\KK)$ denotes the bounded derived category of the
category of finite dimensional $\KK$-vector spaces)
such that for any smooth $k$-scheme $X$, one has the
following canonical identifications (where $\dual M$
denotes the dual of $M$):
$$R(\dual{M_{\mathit{gm}}(X)})\simeq
\dual{R(M_{\mathit{gm}}(X))}\simeq\derR\Gamma(X,E)\, .$$
Moreover, for any object $M$ of $\DM_{\mathit{gm}}(k)_{\QQ}$,
and any integer $p$, one has
$$R(M(p))=R(M)(p)\, .$$
\end{thmm}

These statements are proved using the homotopy theory
of schemes of Morel and Voevodsky.
We work in the stable homotopy
category of motivic symmetric spectra with rational coefficients,
denoted by $\DMt(\spec k,\QQ)$.
We associate canonically to a mixed Weil theory $E$
a commutative ring spectrum $\ste$ such that for any
smooth $k$-scheme $X$ and integers $p$ and $q$,
we get a natural identification
$$H^q(X,E)(p)=H^q(X,\ste(p))\, .$$
We also consider the triangulated category $\DMt(\spec k,\ste)$,
which might be thought of as the category of `motives
with coefficients in $\ste$' (this is simply the
localization of the category of $\ste$-modules by
stable $\AA^1$-equivalences).
We obviously have a symmetric monoidal triangulated
functor
$$\DMt(\spec k,\QQ)\To\DMt(\spec k,\ste)\quad , \qquad
M\longmapsto\ste\otimes^\derL_{\QQ}M\, .$$
If $\Der(\KK)$ is the unbounded derived category
of the category of $\KK$-vector spaces, the result
hidden behind Theorems \ref{thm1} and \ref{thm2} is

\begin{thmm}[Tilting]\label{thm4}
The homological realization functor
$$\DMt(\spec k,\ste)\To\Der(\KK)\quad , \qquad
M\longmapsto\derR\Hom_{\ste}(\ste,M)$$
is an equivalence of symmetric monoidal triangulated categories.
\end{thmm}

To obtain Theorem \ref{thm3}, we
interpret the cycle class map as a map of ring spectra $\HQ\To\ste$,
and use a result of R{\"o}ndigs and {\O}stv{\ae}r
which identifies $\DM(k)_{\QQ}$ with the homotopy category
of modules over the motivic cohomology spectrum $\HQ$.
Note that, by Theorem \ref{thm4},
the homological realization functor
of Theorem \ref{thm3} is essentially
the derived base change functor
$M\longmapsto\ste\otimes^\derL_{\QQ}M$.
This means that the theory of motivic realization functors
is part of (a kind of) tilting theory.

Most of our paper is written over a general regular base $S$ rather than just a perfect field. The first reason for this is that a big part of this
machinery works \emph{mutatis mutandis} over a regular base\footnote{In fact, one could drop the regularity assumption, but then, the formulation of
some of our results about the existence of a cycle class map are a little more involved: $K$-theory is homotopy invariant only for regular schemes.
This is not a serious problem, but we decided to avoid the extra complications due to the fact algebraic $K$-theory is not representable in the
$\AA^1$-homotopy theory of singular schemes.} once we are ready to pay the price of slightly weaker or modified results (we essentially can
say interesting things only for smooth and projective $S$-schemes). When the base is a perfect field, the results announced above are obtained from
the general ones using de~Jong resolution of singularities by alterations. The second reason is that mixed Weil theories defined on smooth schemes
over a complete discrete valuation ring $V$ are of interest: the analog of Theorem \ref{thm2} gives a general way to compare the cohomology of the
generic fiber and of the special fiber of a smooth and projective $V$-scheme (or, more generally, of a smooth $V$-scheme with good properties near
infinity).

Here is a more detailed account on the contents of this paper.

\bigskip

These notes are split into three parts.
The first one sets the basic constructions we need.
That is we construct the `effective' $\AA^1$-derived category
$\DMte(S,R)$ of a scheme $S$ with coefficients in a ring $R$ and recall its main geometrical
and formal properties. We then introduce the Tate object
$R(1)$ and define the category of Tate spectra as the category of symmetric $R(1)$-spectra.
The `non effective' $\AA^1$-derived category $\DMt(S,R)$
of a scheme $S$ with coefficients in a ring $R$ is then the localization of the category
of Tate spectra by stable $\AA^1$-equivalences. We finish
the first part by introducing the $\AA^1$-derived category $\DMt(S,\ste)$ of a scheme $S$
with coefficients in a (commutative) ring spectrum $\ste$
(that is a (commutative) monoid object in the category of symmetric Tate spectra).
The category $\DMt(S,\ste)$ is just defined as the localization of
the category of $\ste$-modules by the class of stable $\AA^1$-equivalences.

The second part is properly about mixed Weil cohomologies. We also define a slightly weaker notion which we call a stable cohomology theory over a
given regular scheme $S$. We associate canonically to any stable cohomology $E$ a commutative ring spectra $\ste$ and a canonical isomorphism
$$E\To\derR\Omega^\infty(\ste)\, .$$
In other words, $E$ can be seen as a kind of `Tate infinite loop space' in the category $\DMte(S,\KK)$. This means in particular that $\ste$
represents in $\DMt(S,\KK)$ the cohomology theory defined by $E$. We get essentially by definition a $1$-periodicity property for $\ste$, that is the
existence of an isomorphism $\ste(1)\simeq\ste$. We then study the main properties of the triangulated category $\DMt(S,\ste)$. In particular, we
prove that Thom spaces are trivial in $\DMt(S,\ste)$, which imply that there is a simple theory of Chern classes and of Gysin maps in $\DMt(S,\ste)$.
Using results of J. Riou, this allows to produce a canonical cycle class map
$$K_{2p-q}(X)^{(p)}_{\QQ}\To H^q(X,\ste(p))=H^q(X,E)(p)$$
(where $K_{q}(X)^{(p)}_{\QQ}$ denotes the part of $K_q(X)$ where the $k^{\text{th}}$ Adams operation acts by multiplication by $k^p$). The good
functoriality properties of Gysin maps implies a Poincar\'e duality theorem in $\DMt(S,\ste)$ for smooth and projective $S$-schemes. In particular,
for any smooth and projective $S$-scheme $X$, the object $\Sigma^\infty(\QQ(X))\otimes^\derL_\QQ\ste$ has a strong dual in $\DMt(S,\ste)$ (i.e. it is
a rigid object). We then prove a weak version of Theorem \ref{thm4}: if $\DMtp(S,\ste)$ denotes the localizing subcategory of $\DMt(S,\ste)$
generated by the objects which have a strong dual, and if $E$ is a mixed Weil theory, then the homological realization functor induces an equivalence
of symmetric monoidal triangulated categories
$$\DMtp(S,\ste)\simeq\Der(\KK)\, . $$
Given a stable theory $E'$, if $\ste'$ denotes the associated ring spectrum, we associate to any morphism of presheaves of differential graded
algebras $E\To E'$ a base change functor
$$\DMt(S,\ste)\To\DMt(S,\ste')$$
whose restriction to $\DMtp(S,\ste)$ happens to be fully faithful whenever $E$ is a mixed Weil theory. In particular, the cohomologies defined by $E$
and $E'$ have then to agree on smooth and projective $S$-schemes.

In the case where $S$ is the spectrum of a perfect field, we prove that
the cycle class map $\HQ\To\ste$, from Voevodsky's rational motivic cohomology spectrum
to our given mixed Weil cohomology $\ste$, is a morphism of commutative
ring spectra. This is achieved by interpreting the cycle class map as an isomorphism
$\ste\otimes^\derL_\QQ\HQ\simeq\ste$ in the homotopy category of $\ste$-modules.
We then observe that the theory of de Jong alterations implies the equality:
$$\DMtp(S,\ste)=\DMt(S,\ste)\, .$$
Using the equivalence of categories $\DMt(S,\HQ)\simeq\DM(k,\QQ)$,
we deduce the expected realization functor from the triangulated category
of mixed motives to the derived category of the category of $\KK$-vector
spaces. In a sequel,
 we also provide a shorter argument which relies on an unpublished result
  of F.Morel stated in \cite{ratmot}.

The last part is an elementary study of some classical mixed Weil theories
We prove that, over a field of characteristic zero, algebraic de Rham
cohomology is a mixed Weil theory, and we explain how Grothendieck's
(resp. Kiehl's) comparison theorem between algebraic and complex analytic (resp.
rigid analytic) de Rham cohomology fits in this picture.

We proceed after this to the study of Monsky-Washnnitzer cohomology as a mixed Weil theory,
and revisit the Berthelot-Ogus Comparison Theorem, which relates de~Rham cohomology
and crystalline cohomology: given a complete discrete valuation ring $V$, with field of
fractions of characteristic zero, and
perfect residue field, for any smooth and proper $V$-scheme $X$ the de Rham cohomology
of the generic fiber of $X$ and the crystalline cohomology of the
special fiber of $X$ are canonically isomorphic. Our proof of this fact also provides a simple
argument to see that the triangulated category of
geometrical mixed motives over $V$ cannot be rigid: we see that, otherwise,
for \emph{any} smooth $V$-scheme $X$, the de Rham cohomology of the
generic fiber of $X$ and the Monsky-Washnitzer cohomology of the special fiber of $X$
would agree, and this is very obviously false in general (for instance, the
special fiber of $X$ might be empty). We also explain how to
define rigid cohomology from
the Monsky-Washnitzer complex using
the natural functorialities of
$\AA^1$-homotopy theory of schemes.

We finally explain an elementary construction of \'etale cohomology
as a mixed Weil cohomology.

As a conclusion, let us mention that this paper deals only with the elementary part of the story:
in a sequel of this paper~\cite{TCMM}, we shall improve this
constructions. In particular, we shall prove that
any mixed Weil cohomology extends naturally to $k$-schemes
of finite type, satisfies h-descent (in particular \'etale descent and
proper descent), and defines a system of triangulated categories on which
the six operations of Grothendieck act. By appyling this construction
to rigid cohomology, this will define a convenient fundation
for a good notion of $p$-adic coefficients.

\bigskip

This paper takes its origins from a seminar on $p$-adic regulators organized by J~.Wildeshaus,
D.~Blotti\`ere and the second named author at university Paris 13. This
is where the authors went to the problem of representing rigid cohomology as a realization
functor of the triangulated category of mixed motives, and
the present paper can be seen as a kind of answer (among others, see for
example \cite{levine,huber}). We would like to thank deeply Y.~Henrio for
all the time he spent to explain the arcanes of rigid analytic geometry and
of $p$-adic cohomology to us. We benefited of valuable discussions with
J.~Ayoub, L.~Breen, W.~Messing and J.~Riou. We feel very grateful to J. Wildeshaus
for his constant warm support and enthusiasm.
We also thank J.~Wildeshaus and J.~I.~Burgos, for their joint careful reading
and valuable comments.


\section{Motivic homological algebra}

All schemes are assumed to be noetherian
and of finite Krull dimension.
We will say `$S$-scheme' for `separated scheme of finite type over $S$'.

If $\A$ is an abelian category, we let $\Comp(\A)$, $\Htp(\A)$ and
$\Der(\A)$ be respectively the category of unbounded cochain
complexes of $\A$, the same category modulo cochain homotopy equivalence,
and the unbounded derived category of $\A$.

\subsection{$\AA^1$-invariant cohomology}

\begin{paragr}\label{axiomesVV}
We suppose given a scheme $S$.
We consider a full subcategory $\V$ of the category $\sm/S$
of smooth $S$-schemes satisfying the following
properties\footnote{In practice, the category $\V$
will be $\sm/S$ itself or the full subcategory of
smooth affine $S$-schemes.}.
\begin{itemize}
\item[(a)] $\AA^n_S$ belongs to $\V$ for $n\geq 0$.
\item[(b)] If $X'\To X$ is an \'etale morphism
and if $X$ is in $\V$, then there exists a Zariski
covering $Y\To X'$ with $Y$ in $\V$.
\item[(c)] For any pullback square of $S$-schemes
$$\xymatrix{
X'\ar[r]\ar[d]_{u'}&X\ar[d]^u\\
Y'\ar[r]&Y}$$
in which $u$ is smooth, if $X$,$Y$ and $Y'$ are in $\V$,
so is $X'$.
\item[(d)] If $X$ and $Y$ are in $\V$, then their disjoint
union $X\amalg Y$ is in $\V$.
\item[(e)] For any smooth $S$-scheme $X$, there exists
a Nisnevich covering $Y\To X$ of $X$
with $Y$ in $\V$.
\end{itemize}
We recall that a Nisnevich covering is a surjective and completely
decomposed \'etale morphism. This defines the
Nisnevich topology on $\V$; see e.g. \cite{nis,KS,TT,MV}.

The last property (e) ensures that
the category of sheaves on $\V$ is equivalent to the category
of sheaves on the category of smooth $S$-schemes as far as we consider
sheaves for the Nisnevich topology (or any stronger one).
\end{paragr}

\begin{paragr}\label{defidsisqu}
A \emph{distinguished square} is a pullback square of schemes
\begin{equation}\label{distsquare}\begin{split}\xymatrix{
W\ar[r]^i\ar[d]_g&V\ar[d]^f\\
U\ar[r]_j&X}
\end{split}\end{equation}
where $j$ is an open immersion and $f$ is an \'etale morphism
such that the induced map from $f^{-1}((X-U)_{\mathit{red}})$ to
$(X-U)_{\mathit{red}}$ is an isomorphism. For such a distinguished
square, the map $(j,f) : U\amalg V\To X$ is a Nisnevich
covering. A very useful property of the Nisnevich topology is that
any Nisnevich covering can be refined by a covering coming from
a distinguished square (as far as we work with noetherian schemes).
This leads to the following characterization of the Nisnevich sheaves.\\
\indent A presheaf $F$ on $\V$ is a sheaf for the Nisnevich topology
if and only if for any distinguished square of shape \eqref{distsquare},
we obtain a pullback square
\begin{equation}\label{distsquare2}\begin{split}\xymatrix{
F(X)\ar[r]^{f^*}\ar[d]_{j^*}&F(V)\ar[d]^{i^*}\\
F(U)\ar[r]_{g^*}&F(W)}
\end{split}\end{equation}
This implies that Nisnevich sheaves are stable by filtering colimits
in the category of presheaves on $\V$. In other words,
if $I$ is a small filtering category,
and if $F$ is a functor from $I$ to the category
of presheaves on $\V$ such that $F_{i}$ is a Nisnevich sheaf
for all $i\in I$, then the presheaf $\limind_{i\in I}F_{i}$
is a Nisnevich sheaf.
\end{paragr}

\begin{paragr}\label{nissheavespremodcat}
We fix a commutative ring $R$.
Let $\mathit{Sh}(\V,R)$ be the category of Nisnevich sheaves
of $R$-modules on $\V$. For a presheaf (of $R$-modules) $F$,
we denote by $F_\nis$ the Nisnevich sheaf associated to $F$.
We can form its derived category
$$\Der(\V,R)=\Der(\mathit{Sh}(\V,R)) \ .$$
More precisely, the category $\Der(\V,R)$ is obtained
as the localization of the category $\Comp(\V,R)$ of
(unbounded) complexes of the Grothendieck abelian category
$\mathit{Sh}(\V,R)$ by the class of
quasi-isomorphisms. As we have an equivalence of categories
$$\mathit{Sh}(\V,R)\simeq\mathit{Sh}(\sm/S,R) \ , $$
we also have a canonical equivalence of categories
\begin{equation}\label{Nispremodcat1}
\Der(\V,R)\simeq\Der(\sm/S,R) \ .
\end{equation}
We have a canonical functor
\begin{equation}
R:\V\To\mathit{Sh}(\V,R)\quad , \qquad X\longmapsto R(X)
\end{equation}\label{Nispremodcat2}
where $R(X)$ denotes the Nisnevich sheaf associated to the
presheaf
$$Y\longmapsto\text{free $R$-module generated by $\Hom^{}_\V(Y,X)$.}$$

Note that according to \ref{defidsisqu},
for any $X$ in $\V$, and any small filtering system $(F_i)_{i \in I}$
of $\mathit{Sh}(\V,R)$, the canonical map
\begin{equation}\label{filtglobsection}
\limind_{i\in I}\Hom_{\mathit{Sh}(\V,R)}(R(X),F_{i})\To
\Hom_{\mathit{Sh}(\V,R)}(R(X),\limind_{i\in I}F_{i})
\end{equation}
is an isomorphism (we can even take
$X$ to be any smooth $S$-scheme according to the previous equivalence).

For a complex $K$ of presheaves of $R$-modules on $\V$, we have
a canonical isomorphism
\begin{equation}\label{Nispremodcat3}
\hypercoh^n_\nis(X,K_\nis)=\Hom^{}_{\Der(\V,R)}(R(X),K_\nis[n])
\end{equation}
where $X$ is an object of $\V$, $n$ is an integer, and
$\hypercoh^n_\nis(X,K_\nis)$ is the
Nisnevich hypercohomology with coefficients in $K$.
\end{paragr}

\begin{paragr}\label{predefcmfnis}
For a sheaf of $R$-modules $F$, and an integer $n$,
we denote by $D^nF$ the complex concentrated in degrees $n$
and $n+1$ whose only non trivial differential is the identity of $F$.
We write $S^{n}F$ for the sheaf $F$ seen as a complex concentrated
in degree $n$. We have a canonical inclusion of $S^{n+1}F$ in $D^nF$.
We say that a morphism of complexes of sheaves of $R$-modules
is a \emph{$\V$-cofibration} if it is contained in the smallest
class of maps stable by pushout, transfinite composition and retract
that contains the maps of the form $S^{n+1}R(X)\To D^nR(X)$
for any integer $n$ and any $X$ in $\V$. For example, for any $X$ in $\V$,
the map $0\To R(X)$ is a $\V$-cofibration (where $R(X)$ is seen
as a complex concentrated in degree $0$). A complex of presheaves $K$ is
\emph{$\V$-cofibrant} if $0\To K$ is a $\V$-cofibration.

A complex of presheaves of $R$-modules $K$ on $\V$
is \emph{$\V_\nis$-local} if for any $X$ in $\V$, the canonical map
$$H^n(K(X))\To \hypercoh^n_\nis(X,K_\nis)$$
is an isomorphism of $R$-modules.

A morphism $p:K\To L$ of complexes of presheaves of $R$-modules on $\V$
is $\V$-surjective if for any $X$ in $\V$, the map
$K(X)\To L(X)$ is surjective.
\end{paragr}

\begin{prop}\label{existmodcatlocalnis}
The category of complexes of Nisnevich sheaves of $R$-modules on $\V$
is a proper Quillen closed model category structure whose weak equivalences are the
quasi-isomorphisms, whose cofibrations are the
$\V$-cofibrations and whose fibrations are the $\V$-surjective
morphisms with $\V_{\nis}$-local kernel.
In particular, for any $X$ in $\V$, $R(X)$ is $\V$-cofibrant.
\end{prop}

\begin{proof}
If $\X$ is a simplicial object of $\V$, we denote by
$R(\X)$ the associated complex. If $\X$ is a Nisnevich hypercovering
of an object $X$ of $\V$, we have a canonical morphism
from $R(\X)$ to $R(X)$ and we define
$$\widetilde{R}(\X)={\mathrm{Cone}}\big(R(\X)\To R(X)\big) \ . $$
Let $\G$ be the collection of the $R(X)$'s for $X$ in $\V$,
and $\H$ the class of the $\widetilde{R}(\X)$'s
for all the Nisnevich hypercoverings of any object $X$ of $\V$.
Then $(\G,\H)$ is a descent structure
on $\mathit{Sh}(\V,R)$ as defined in \cite[Definition 1.4]{HCD},
so that we can apply \cite[Theorem 1.7 and Corollary 4.9]{HCD}.
\end{proof}

\begin{paragr}
The model structure above will be called the \emph{$\V$-local
model structure}.
\end{paragr}

\begin{cor}\label{flasqueresolution}
For any complex of Nisnevich sheaves of $R$-modules $K$,
there exists a quasi-isomorphism $K\To L$
where $L$ is $\V_\nis$-local.
\end{cor}

\begin{proof}
We just have to choose a factorization of
$K\To 0$ into a quasi-isomorphism $K\To L$
followed by a fibration $L\To 0$ for the above
model structure.
\end{proof}

\begin{prop}\label{exactdistsquare}
For any distinguished square
$$\xymatrix{
W\ar[r]^i\ar[d]_g&V\ar[d]^f\\
U\ar[r]_j&X}$$
the induced commutative square of sheaves of $R$-modules
$$\xymatrix{
R(W)\ar^{i*}[r]\ar_{g_*}[d]&R(V)\ar^{f_*}[d]\\
R(U)\ar_{j_*}[r]&R(X)}$$
is exact; this means that it is cartesian and cocartesian,
or equivalently that it gives rise to a short exact
sequence in the category of sheaves of $R$-modules
$$0\To R(W)
 \xrightarrow{g_*-i_*} R(U)\oplus R(V)
 \xrightarrow{(j_*,f_*)} R(X)\To 0 \ .$$
\end{prop}

\begin{proof}
The characterization of Nisnevich sheaves given in \ref{defidsisqu}
implies that the sequence
$$0\To R(W)\To R(U)\oplus R(V)\To R(X)\To 0$$
is right exact. So the result comes from the injectivity of the
map from $R(W)$ to $R(V)$ induced by $i$.
\end{proof}

\begin{paragr}
Let $K$ be a complex of presheaves of $R$-modules on $\V$.

A \emph{closed pair} will be a couple $(X,Z)$ such that
$X$ is a scheme in $\V$, $Z \subset X$ is a closed subset
and $X-Z$ belongs to $\V$.  Let $j$ be the immersion
of $X-Z$ in $X$. We put
$$
K_Z(X)=\mathrm{Cone}\big(K(X) \xrightarrow{j^*} K(X-Z)\big)[-1].
$$
A morphism of closed pairs
$f:(Y,T) \rightarrow (X,Z)$ is a morphism of schemes
$f:Y \rightarrow X$ such that $f^{-1}(Z) \subset T$.
The morphism of closed pairs $f$ will be called \emph{excisive}
when the induced square
$$\xymatrix{
Y-T\ar[r]^/4pt/i\ar[d]_g&Y\ar[d]^f\\
X-Z\ar[r]_/4pt/j&X}$$
is distinguished.

The complex $K_Z(X)$ is obviously functorial with respect
to morphisms of closed pairs. We will say that $K$
has the \emph{excision property on $\V$}
if for any excisive morphism $f:(Y,T) \rightarrow (X,Z)$
the map $K_T(Y) \rightarrow K_Z(X)$ is a quasi-isomorphism.

We will say that $K$ has the \emph{Brown-Gersten property on $\V$
with respect to the Nisnevich topology}, or the \emph{B.-G.-property} for short,
if for any distinguished square
$$\xymatrix{
W\ar[r]^i\ar[d]_g&V\ar[d]^f\\
U\ar[r]_j&X}$$
in $\V$, the square
$$\xymatrix{
K(X)\ar[r]^{f^*}\ar[d]_{j^*}&K(V)\ar[d]^{i^*}\\
K(U)\ar[r]_{g^*}&K(W)}$$
is a homotopy pullback (or equivalently a
homotopy pushout) in the category of complexes of $R$-modules.
The latter condition means that the commutative square
of complexes of $R$-modules
obtained from the distinguished square above by applying $K$
leads canonically to a long exact sequence ``\`a la Mayer-Vietoris''
$$H^n(K(X))
 \xrightarrow{j^*+f^*} H^n(K(U))\oplus H^n(K(V))
 \xrightarrow{g^*-i^*} H^n(K(W)) \To H^{n+1}(K(X))$$
 The complexes satisfying the B.-G.-property are in fact
the fibrant objects of the model structure of Proposition \ref{existmodcatlocalnis}.
This is shown by the following result which is essentially due to
Morel and Voevodsky.
\end{paragr}

\begin{prop}\label{bgtilde}
Let $K$ be a complex of presheaves of $R$-modules on $\V$.
Then the following conditions are equivalent.
\begin{itemize}
\item[(i)] The complex $K$ has the B.-G.-property.
\item[(i$^{\prime}$)] The complex $K$ has the excision property.
\item[(ii)] For any $X$ in $\V$, the canonical map
$$H^n(K(X))\To \hypercoh^n_\nis(X,K_\nis)$$
is an isomorphism of $R$-modules (\emph{i.e.} $K$ is $\V_\nis$-local).
\end{itemize}
\end{prop}
\begin{proof}
The equivalence of (i) and (i$^{\prime}$) follows from the definition
of a homotopy pullback.

As any short exact sequence of sheaves of $R$-modules
defines canonically a distinguished triangle in $\Der(\V,R)$,
the fact that (ii) implies (i) follows easily from
proposition \ref{exactdistsquare}.
To prove that (i) implies (ii), we need a little more
machinery. First, we can choose a monomorphism of complexes $K\To L$
which induces a quasi-isomorphism between $K_\nis$ and $L_\nis$,
and such that $L$ is $\V_\nis$-local. For this, we first choose
a quasi-isomorphism $K_\nis\To M$ where $M$ is $\V_\nis$-local
(which is possible by Corollary \ref{flasqueresolution}).
We have a natural embbeding of $K$ into
the mapping cone of its identity $\mathrm{Cone}(1_K)$.
But $\mathrm{Cone}(1_K)$ is obviously $\V_\nis$-local
as it is already acyclic as a complex of presheaves.
This implies that the direct sum $L=\mathrm{Cone}(1_K)\oplus M$
is also $\V_\nis$-local. Moreover, as
$K$ and $L$ both satisfy the B.-G.-property,
one can check easily that the quotient presheaf $L/K$ also
has the B.-G.-property. Hence it is sufficient to prove that
$H^n(L(X)/K(X))=0$ for any $X$ in $\V$ and any integer $n$.\\
\indent Let us fix an object $X$ of $\V$ and an integer $n$.
One has to consider for any $q \geq 0$ the presheaf $T_q$ on the small
Nisnevich site $X_\nis$ of $X$ defined by
$$T_q(Y)=H^{n-q}(L(Y)/K(Y)) \ . $$
These are B.-G.-functors as defined by Morel and
Voevodsky \cite[proof of Prop. 1.16, page 101]{MV},
and for any integer $q\geq 0$, the Nisnevich sheaf associated to $T_q$
is trivial (this is because $K_\nis\To L$ is a quasi-isomorphism
of complexes of Nisnevich sheaves by construction of $L$).
This implies by virtue of \cite[Lemma 1.17, page 101]{MV}
that $T_q=0$ for any $q\geq 0$. In particular,
we have $H^n(L(X)/K(X))=0$. Therefore $L/K$ is an acyclic
complex of presheaves over $\V$, and $K\To L$ is a quasi-isomorphism
of complexes of presheaves. This proves that $K$ is $\V_\nis$-local
if and only if $L$ is, hence that $K$ is $\V_\nis$-local.
\end{proof}

\begin{cor}\label{filtcolimNiscoh}
Let $I$ be a small filtering category, and
$K$ a functor from $I$ to the category of complexes of
Nisnevich sheaves of $R$-modules. Then for any smooth
$S$-scheme $X$, the canonical maps
$$\limind_{i\in I}\hypercoh^n_\nis(X,K_{i})\To
\hypercoh^n_\nis(X,\limind_{i\in I}K)$$
are isomorphisms for all $n$.
\end{cor}

\begin{proof}
We can suppose that $K_{i}$ is $\V_{\nis}$-local for all $i\in I$
(we can take a termwise fibrant replacement of $K$ with respect to the
model structrure of Proposition \ref{existmodcatlocalnis}).
It then follows from Proposition \ref{bgtilde}
that $\limind_{i\in I}K_{i}$ is still $\V_{\nis}$-local:
it follows from the fact that the filtering colimits are exact
that the presheaves with the B.-G.-property are
stable by filtering colimits. The map
$$\limind_{i\in I}H^n(K_{i}(X))\To
H^n(\limind_{i\in I}K_{i}(X))$$
is obviously an isomorphism for any $X$ in $\V$.
As we are free to take $\V=\sm/S$, this proves the assertion.
\end{proof}

\begin{paragr}\label{defNiscompact}
Remember that if $\T$ is a triangulated category with
small sums, an object $X$ of $\T$ is \emph{compact}
if for any small family $(K_{\lambda})_{\lambda\in \Lambda}$
of objects of $\T$, the canonical map
$$\bigoplus_{\lambda\in\Lambda}\Hom_{\T}\big(X, K_{\lambda}\big)
\To\Hom_{\T}\big(X,\bigoplus_{\lambda\in\Lambda} K_{\lambda}\big)$$
is bijective (as this map is always injective, this is
equivalent to say it is surjective). One denotes by $\T_{c}$
the full subcategory of $\T$ that consists of compact objects.
It is easy to see that $\T_{c}$ is a thick subcategory of $\T$
(which means that $\T_{c}$ is a triangulated subcategory of $\T$ stable by direct factors).
\end{paragr}

\begin{cor}\label{filtcolimNiscoh2}
For any smooth $S$-scheme $X$, $R(X)$
is a compact object of the derived category
of Nisnevich sheaves of $R$-modules.
\end{cor}

\begin{proof}
As any direct sum is a filtering colimit of
finite direct sums, this follows from \eqref{Nispremodcat3}
and Corollary \ref{filtcolimNiscoh}.
\end{proof}

\begin{paragr}\label{defA1equiv}
Let $\Der$ be a triangulated category.
Remember that a \emph{localizing subcategory} of $\Der$
is a full subcategory $\TDer$ of $\Der$ with the following properties.
\begin{itemize}
\item[(i)] $A$ is in $\TDer$ if and only if $A[1]$ is in $\TDer$.
\item[(ii)] For any distinguished triangle
$$A'\To A\To A''\To A'[1] \ , $$
if $A'$ and $A''$ are in $\TDer$, then $A$ is in $\TDer$.
\item[(iii)] For any (small) family $(A_i)_{i\in I}$
of objects of $\TDer$, $\bigoplus_{i\in I}A_i$ is in $\TDer$.
\end{itemize}
If $\T$ is a class of objects of $\Der$, the \emph{localizing
subcategory of $\Der$ generated by $\T$} is the smallest
localizing subcategory of $\Der$ that contains $\T$ (\emph{i.e.}
the intersection of all the localizing subcategories of $\Der$
that contain $\T$).\\
\indent Let $\T$ be the class of complexes of shape
$$\dots\To 0\To R(X\times_S\AA^1_S)\To R(X)\To 0\To\dots$$
with $X$ in $\V$ (the non trivial differential is induced
by the canonical projection). Denote by $\TDer(\V,\AA^1,R)$ the localizing
subcategory of $\Der(\V,R)$ generated by $\T$. We define the
triangulated category \smash{$\Der^{\mathit{eff}}_{\AA^1}(\V,R)$}
as the Verdier quotient of $\Der(\V,R)$ by $\TDer(\V,\AA^1,R)$.
$$\Der^{\mathit{eff}}_{\AA^1}(\V,R)=\Der(\V,R)/\TDer(\V,\AA^1,R)$$
We know that $\Der(\V,R)\simeq\Der(\sm/S,R)$,
and an easy Mayer-Vietoris argument for the Zariski topology
shows that the essential image of $\TDer(\V,\AA^1,R)$ in $\Der(\sm/S,R)$
is precisely $\TDer(\sm/S,\AA^1,R)$. Hence we get a canonical equivalence
of categories
\begin{equation}\label{eqdmdmV}
\Der^{\mathit{eff}}_{\AA^1}(\V,R)\simeq\Der^{\mathit{eff}}_{\AA^1}(\sm/S,R) \ .
\end{equation}
We simply put:
$$\DMte(S,R)=\Der^{\mathit{eff}}_{\AA^1}(\sm/S,R) \ .$$
According to F. Morel insights, the category $\DMte(S,R)$
is called the \emph{triangulated category of effective real
motives}\footnote{This terminology comes from the fact
$\DMte(S,R)$ give quadratic informations on $S$, which implies
it is bigger than Voevodsky's triangulated category of
mixed motives; see \cite{pizero,ratmot}.
The word `real' is meant here as opposed to `complex'.}
(with coefficients in $R$). In the sequel of this paper, we will
consider the equivalence \eqref{eqdmdmV} as an
equality\footnote{The role of the category $\V$ is only to define
model category structures on the category of complexes of
$\mathit{Sh}(\V,R)\simeq\mathit{Sh}(\sm/S,R)$ which
depend only on the local behaviour of
the schemes in $\V$ (e.g. the smooth
affine shemes over $S$).}.
We thus have a canonical localization functor
\begin{equation}\label{deflocdmequat}
\gamma : \Der(\V,R)\To\DMte(S,R) \ .
\end{equation}
We will say that a morphism of complexes of $\mathit{Sh}(\V,R)$
is an \emph{$\AA^1$-equivalence} if its image
in \smash{$\DMte(S,R)$}
is an isomorphism.\\
\indent A complex of presheaves of $R$-modules $K$ over $\V$
is \emph{$\AA^1$-homotopy invariant} if for any $X$ in $\V$,
the projection of $X\times_S\AA^1_S$ on $X$ induces a quasi-isomorphism
$$K(X)\To K(X\times_S\AA^1_S) \ .$$
\end{paragr}

\begin{prop}\label{A1cmfeff}
The category of complexes of $\mathit{Sh}(\V,R)$ is endowed with
a proper Quillen model category structure whose weak equivalences are the $\AA^1$-equivalences,
whose cofibrations are the $\V$-cofibrations, and whose fibrations
are the $\V$-surjective morphisms with $\AA^1$-homotopy
invariant and $\V_{\nis}$-local kernel.
In particular, the fibrant objects of this model structure are exactly
the $\AA^1$-homotopy invariant and $\V_\nis$-local complexes.
The corresponding homotopy category is the triangulated
category of effective real motives
\smash{$\DMte(S,R)$}.
\end{prop}

\begin{proof}
This is a direct application of \cite[Proposition 3.5 and Corollary 4.10]{HCD}.
\end{proof}

\begin{paragr}
Say that a complex $K$ of presheaves of $R$-modules
on $\V$ is \emph{$\AA^1$-local}
if for any $X$ in $\V$, the projection of $X\times_S\AA^1_S$
on $X$ induces isomorphisms in Nisnevich hypercohomology
$$\hypercoh^n_\nis(X,K_\nis)\simeq
\hypercoh^n_\nis(X\times_S\AA^1_S,K_\nis) \ . $$
It is easy to see that a $\V_\nis$-local complex
is $\AA^1$-local if and only if it is $\AA^1$-homotopy
invariant. In general, a complex of sheaves $K$
is $\AA^1$-local if and only if, for any quasi-isomorphism
$K\To L$, if $L$ is $\V_\nis$-local, then $L$ is $\AA^1$-homotopy
invariant. We deduce from this the following result.
\end{paragr}

\begin{cor}
The localization functor $\Der(\V,R)\To\DMte(S,R)$
has a right adjoint that is fully faithful and whose essential
image consists of the $\AA^1$-local complexes. In other words,
\smash{$\DMte(S,R)$} is canonically equivalent
to the full subcategory of $\AA^1$-local complexes in $\Der(\V,R)$.
\end{cor}

\begin{proof}
For any complex of sheaves $K$, one can produce functorially
a map $K\To L^{}_{\AA^1}K$ which is an $\AA^1$-equivalence
with $L^{}_{\AA^1}K$ a $\V_\nis$-local and $\AA^1$-homotopy
invariant complex (just consider a functorial fibrant resolution
of the model category of \ref{A1cmfeff}). Then the functor
$L^{}_{\AA^1}$ takes $\AA^1$-equivalences to quasi-isomorphisms
of complexes of (pre)sheaves, and induces a functor
$$L^{}_{\AA^1} : \DMte(S,R)\To\Der(\V,R)$$
which is the expected right adjoint of the localization
functor.
\end{proof}

\begin{ex}\label{exampleR}
The constant sheaf $R$ is $\AA^1$-local (if it is considered as a complex
concentrated in degree $0$).
\end{ex}

\begin{paragr}\label{changeringcoeff}
Let $R'$ be another commutative ring, and $R\To R'$
a morphism of rings. The functor $K\longmapsto K\otimes^{}_{R}R'$
is a symmetric monoidal left Quillen functor from $\Comp(\V,R)$ to $\Comp(\V,R')$
for the model structures of Proposition \ref{A1cmfeff}. Hence it has a total
left derived functor
$$\DMte(S,R)\To\DMte(S,R')\qquad ,
\quad K\longmapsto K\otimes^\derL_{R}R'$$
whose right adjoint is the obvious forgetful functor.
I.e. for a complex of sheaves of $R$-modules $K$
and a complex of sheaves of $R'$-modules $L$, we have a canonical
isomorphism
$$\Hom_{\DMte(S,R)}(K,L)\simeq
\Hom_{\DMte(S,R')}(K\otimes^\derL_{R}R',L) \ .$$
\end{paragr}

\begin{ex}\label{exampleGm}
Let $\GG_{m}=\AA^1_{S}-\{0\}$ be the multiplicative group.
It can be considered as a presheaf of groups on $\sm/S$,
and one can check that it is a Nisnevich sheaf. Moreover,
for any smooth $S$-scheme $X$, one has
\begin{equation}\label{computeGmcoh}
H^{i}_{\nis}(X,\GG_{m})=\begin{cases}
\bundle O^*(X)&\text{if $i=0$,}\\
\pic(X)&\text{if $i=1$,}\\
0&\text{otherwise.}
\end{cases}
\end{equation}
As $S$ is assumed to be regular,
$\GG_{m}$ is $\AA^1$-local as a complex concentrated in degree $0$.
We deduce that we have the formula
\begin{equation}\label{GminDmeffk}
H^{i}_{\nis}(X,\GG_{m})=\Hom_{\DMte(k,\ZZ)}(\ZZ(X),\GG_{m}[i]) \ .
\end{equation}
In particular, it follows from \ref{changeringcoeff} that for any smooth
$S$-scheme $X$, one has a canonical morphism of abelian groups
\begin{equation}\label{PicDmeffkR}
\pic(X)\To
\Hom_{\DMte(S,R)}(R(X),\GG_{m}\otimes^\derL_{\ZZ}R[1]) \ .
\end{equation}
If moreover $R$ is flat over $\ZZ$, we get the formula
\begin{equation}\label{PicDmeffkRflat}
\pic(X)\otimes^{}_{\ZZ}R\simeq
\Hom_{\DMte(S,R)}(R(X),\GG_{m}\otimes^\derL_{\ZZ}R[1]) \ .
\end{equation}
\end{ex}

\begin{prop}\label{A1filtcolimNiscoh}
Let $I$ be a small filtering category, and
$K$ a functor from $I$ to the category of complexes of
Nisnevich sheaves of $R$-modules. Then for any smooth
$S$-scheme $X$, the canonical map
$$\limind_{i\in I}\Hom_{\DMte(S,R)}(R(X),K_{i})\To
\Hom_{\DMte(k,R)}(R(X),\limind_{i\in I}K)$$
is an isomorphism.
\end{prop}

\begin{proof}
This follows from Corollary \ref{filtcolimNiscoh}
once we see that $\AA^1$-homotopy invariant
complexes are stable by filtering colimits.
\end{proof}

\begin{cor}\label{repcompactdmtildeeff}
For any smooth $S$-scheme $X$, $R(X)$
is a compact object of $\DMte(S,R)$.
\end{cor}

\subsection{Derived tensor product and derived $\Hom$}

\begin{paragr}
We consider a full subcategory $\V$ of $\sm/S$ as in \ref{axiomesVV}
and a commutative ring $R$.
The category of sheaves of $R$-modules on $\V$ has a tensor product $\otimes^{}_R$
defined in the usual way: if $F$ and $G$ are two sheaves of $R$-modules,
$F\otimes^{}_R G$ is the Nisnevich sheaf associated to the presheaf
$$X\longmapsto F(X)\otimes^{}_R G(X) \ . $$
the unit of this tensor product is $R=R(S)$.
This makes the category $\mathit{Sh}(\V,R)$ a closed symmetric
monoidal category. For two objects $X$ and $Y$ of $\V$, we have a canonical
isomorphism
$$R(X\times_S Y)\simeq R(X)\otimes^{}_R R(Y) \ . $$
Finally, an important property of this tensor product is that
for any $X$ in $\V$, the sheaf $R(X)$ is flat, by which we mean
that the functor
$$F\longmapsto R(X)\otimes^{}_R F$$
is exact. This implies that the family of the sheaves
$R(X)$ for $X$ in $\V$ is \emph{flat} in the sense of \cite[2.1]{HCD}.
Hence we can apply
Corollary 2.6 of \emph{loc. cit.} to get that
the $\V_\nis$-local model structure of \ref{existmodcatlocalnis}
is compatible with the tensor product
in a very (rather technical but also)
gentle way: define the tensor product of two complexes
of sheaves of $R$-modules $K$ and $L$ on $\V$ by the formula
$$(K\otimes^{}_R L)^n=\bigoplus_{p+q=n}K^p\otimes^{}_R L^q$$
with differential $d(x\otimes y)=dx\otimes y+(-1)^{\mathrm{deg}(x)} x\otimes dy$.
This defines a structure of symmetric monoidal category
on $\Comp(\mathit{Sh}(\V,R))$
(the unit is just $R$ seen as complex
concentrated in degree $0$, and the symmetry rule is given by
the usual formula
$x\otimes y\longmapsto (-1)^{\mathrm{deg}(x)\mathrm{deg}(y)}y\otimes x$).
A consequence of \emph{loc. cit.} Corollary 2.6 is that the functor
$(K,L)\longmapsto K\otimes^{}_R L$
is a left Quillen bifunctor, which implies in particular that it has
a well behaved total left derived functor
$$\Der(\V,R)\times\Der(\V,R)\To\Der(\V,R)\quad ,
\qquad (K,L)\longmapsto K\otimes^\derL_R L \ .$$
Moreover, for a given complex $L$, the functor
$$\Der(\V,R)\To\Der(\V,R)\quad ,
\qquad K\longmapsto K\otimes^\derL_R L$$
is the total left derived functor of the functor
$K\longmapsto K\otimes^{}_R L$ (see Remark 2.9 of \emph{loc. cit.}).
This means that for
any $\V$-cofibrant complex $K$ (that is, a complex such that
the map $0\To K$ is a $\V$-cofibration), the canonical map
$$K\otimes^\derL_R L\To K\otimes^{}_R L$$
is an isomorphism in $\Der(\V,R)$ for any complex $L$. In particular,
if $F$ is a direct factor of some $R(X)$ (with $X$ in $\V$),
then for any complex of sheaves $L$, the map
$$F\otimes^\derL_R L\To F\otimes^{}_R L$$
is an isomorphism in $\Der(\V,R)$.
This derived tensor product makes $\Der(\V,R)$ a
closed symmetric monoidal triangulated category.
This means that for two objects $L$ and $M$ of
$\Der(\V,R)$, there is an object $\derR\sHom(L,M)$ of $\Der(\V,R)$
that is defined by the universal property
$$\forall K\in\Der(\V,R)\ , \quad
\Hom_{\Der(\V,R)}(K\otimes^\derL_{R}L,M)\simeq
\Hom_{\Der(\V,R)}(K,\derR\sHom(L,M)).$$
The functor $\derR\sHom$ can also be characterized as the
total right derived functor of the internal Hom of the category
of complexes of Nisnevich sheaves of $R$-moddules on $\V$.
If $L$ is $\V$-cofibrant and if $M$ if $\V_{\nis}$-local, then
$\derR\sHom(L,M)$ can be represented by the complex
of sheaves
$$X\longmapsto\mathrm{Tot}\big [\Hom_{\mathit{Sh}(\V,R)}(R(X)\otimes^{}_{R}L,M)\big ].$$

The derived tensor product on $\Der(\V,R)$
induces a derived tensor product on $\DMte(S,R)$
as follows.
\end{paragr}

\begin{prop}\label{A1derivedtensor}
The tensor product of complexes $\otimes^{}_R$ has a total
left derived functor
$$\DMte(S,R)\times\DMte(S,R)
\To\DMte(S,R)\quad ,
\qquad (K,L)\longmapsto K\otimes^\derL_R L $$
that makes $\DMte(S,R)$ a closed symmetric
tensor triangulated category.
Moreover, the localization functor
$\Der(\V,R)\To\DMte(S,R)$
is a triangulated symmetric monoidal functor.
\end{prop}

\begin{proof}
This follows easily from \cite[Corollary 3.14]{HCD} applied to the classes $\G$
and $\H$ defined in the proof of \ref{existmodcatlocalnis}
and to the class $\T$ of complexes of shape
$$\cdots\To 0\To R(X\times_S\AA^1_S)\To R(X)\To 0\To\cdots$$
with $X$ in $\V$.
\end{proof}

\begin{paragr}\label{A1compactRHomeffprep}
It follows from Proposition \ref{A1derivedtensor} that
the category $\DMte(S,R)$ has an internal Hom that
we still denote by $\derR\sHom$. Hence for three objects $K$, $L$ and $M$
in $\DMte(S,R)$, we have a canonical isomorphism
\begin{equation}\label{A1RHomeff}
\Hom_{\DMte(S,R)}(K\otimes^\derL_{R}L,M)\simeq
\Hom_{\DMte(S,R)}(K,\derR\sHom(L,M)).
\end{equation}
If $L$ is $\V$-cofibrant and $M$ is $\V_{\nis}$-local and $\AA^1$-local,
then
\begin{equation}\label{A1RHomeff2}
\derR\sHom(L,M)=\mathrm{Tot}\big [ \Hom_{\mathit{Sh}(\V,R)}(R(-)\otimes^{}_{R}L,M)\big ].
\end{equation}
\end{paragr}

\begin{prop}\label{A1compactRHomeff}
If $L$ is a compact object of $\DMte(S,R)$, then
for any small family $(K_{\lambda})_{\lambda\in \Lambda}$
of objects of $\DMte(S,R)$, the canonical map
$$\bigoplus_{\lambda\in\Lambda}\derR\sHom\big(L,K_{\lambda}\big)\To
\derR\sHom\big(L,\bigoplus_{\lambda\in\Lambda} K_{\lambda}\big)$$
is an isomorphism in $\DMte(S,R)$.
\end{prop}

\begin{proof}
Once the family $(K_{\lambda})_{\lambda\in \Lambda}$ is fixed, this map defines
a morphism of triangulated functors from the triangulated category of compact objects of
$\DMte(S,R)$ to $\DMte(S,R)$.
Therefore, it is sufficient to check this property when $L=R(Y)$
with $Y$ in $\V$. This is equivalent to say that for any $X$ in $\V$,
the map
$$\Hom
\big(R(X),\bigoplus_{\lambda\in\Lambda}\derR\sHom\big(R(Y),K_{\lambda}\big)\big)\To
\Hom
\big(R(X),\derR\sHom\big(R(Y),\bigoplus_{\lambda\in\Lambda} K_{\lambda}\big)\big)$$
is bijective.
As $R(X)$ is compact (\ref{repcompactdmtildeeff}), we have
$$\Hom
\big(R(X),\bigoplus_{\lambda\in\Lambda}\derR\sHom\big(R(Y),K_{\lambda}\big)\big)
\simeq\bigoplus_{\lambda\in\Lambda}
\Hom
\big(R(X),\derR\sHom\big(R(Y),K_{\lambda}\big)\big),$$
and as $R(X\times_{S}Y)\simeq R(X)\otimes^\derL_{R}R(Y)$ is compact as well,
we have
$$\begin{aligned}
\Hom
\big(R(X),\derR\sHom\big(R(Y),\bigoplus_{\lambda\in\Lambda} K_{\lambda}\big)\big)&\simeq
\Hom
\big(R(X)\otimes^\derL_{R}R(Y),\bigoplus_{\lambda\in\Lambda} K_{\lambda}\big)\\
&\simeq\bigoplus_{\lambda\in\Lambda}
\Hom\big(R(X)\otimes^\derL_{R}R(Y),K_{\lambda}\big)\\
&\simeq\bigoplus_{\lambda\in\Lambda}
\Hom\big(R(X),\derR\sHom(R(Y),K_{\lambda})\big).
\end{aligned}$$
This implies our claim immediately.
\end{proof}

\begin{paragr}\label{quillenconstantsheaf}
Let $R\Mod$ be the category of $R$-modules.
If $M$ is an $R$-module, we still denote by $M$
the constant Nisnevich sheaf of $R$-modules on $\V$
generated by $M$.
This defines a symmetric monoidal functor from
the category of (unbounded) complexes of $R$-modules $\Comp(R)$
to the category $\Comp(\V,R)$
\begin{equation}\label{quillenconstantsheaf1}
\Comp(R)\To\Comp(\V,R)\quad , \qquad
M\longmapsto M.
\end{equation}
This functor is a left adjoint to the global sections functor
\begin{equation}\label{quillenconstantsheaf2}
\Gamma : \Comp(\V,R)\To\Comp(R)\quad , \qquad
M\longmapsto \Gamma(M)=\Gamma(S , M).
\end{equation}

The category $\Comp(R)$ is a Quillen model category
with the quasi-isomorphims as weak equivalences and the
degreewise surjective maps as fibrations (see \emph{e.g.} \cite[Theorem 2.3.11]{Hov}).
We call this model structure the \emph{projective model structure}.
This implies that the constant sheaf functor \eqref{quillenconstantsheaf1}
is a left Quillen functor for the model structures of
Propositions \ref{existmodcatlocalnis} and \ref{A1cmfeff} on
$\Comp(\V,R)$. Therefore, the global sections functor \eqref{quillenconstantsheaf2}
is a right Quillen functor and has total right derived functor
\begin{equation}\label{quillenconstantsheaf3}
\derR\Gamma : \DMte(S,R)\To\Der(R)
\end{equation}
where $\Der(R)$ denotes the derived category of $R$.
For two objects $M$ and $N$ of $\DMte(S,R)$,
we define
\begin{equation}\label{quillenconstantsheaf4}
\derR\Hom(M,N)=\derR\Gamma\big(\derR\sHom(M,N)\big).
\end{equation}
We invite the reader to check that $\derR\Hom$ is the derived Hom
of $\DMte(S,R)$. In particular, for any integer $n$, we have
a canonical isomorphism
\begin{equation}\label{quillenconstantsheaf5}
H^n\big(\derR\Hom(M,N)\big)\simeq\Hom_{\DMte(S,R)}(M,N[n]).
\end{equation}
\end{paragr}

\subsection{Tate object and purity}

\begin{paragr}
Let $\GG_m=\AA^1_S-\{0\}$ be the multiplicative group scheme over $S$.
The unit of $\GG_m$ defines a morphism $R=R(S)\To R(\GG_m)$,
and we define the \emph{Tate object} $R(1)$ as the cokernel
$$R(1)=\mathrm{coker}\big(R\To R(\GG_m)\big)[-1]$$
(this definition makes sense in the category of complexes of
$\mathit{Sh}(\V,R)$ as well as in $\Der(\V,R)$
or in $\DMte(S,R)$ as we take the cokernel
of a split monomorphism). By definition, $R(1)[1]$
is a direct factor of $R(\GG_m)$, so that $R(1)$
is $\V$-cofibrant. Hence for any integer $n\geq 0$,
$R(n)=R(1)^{\otimes n}$ is also $\V$-cofibrant.
For a complex $K$ of $\mathit{Sh}(\V,R)$,
we define $K(n)=K\otimes^{}_R R(n)$.
As $R(n)$ is $\V$-cofibrant, the map
$$K\otimes^\derL_RR(n)\To K\otimes^{}_RR(n)=K(n)$$
is an isomorphism in $\DMte(k,R)$. Another
description of $R(1)$ is the following.
\end{paragr}

\begin{prop}\label{motivicsphere}
The inclusion of $\GG_{m}$ in $\AA^{1}$ induces
a canonical split distinguished triangle in $\DMte(S,R)$
$$R(\GG_{m})\To R(\AA^1_{S})\overset{0}{\To} R(1)[2]\To R(\GG_{m})[1]$$
that gives the canonical decomposition $R(\GG_{m})=R\oplus R(1)[1]$.
\end{prop}

\begin{proof}
This follows formally from the definition
of $R(1)$ and from the fact that $R(\AA^1)=R$ in $\DMte(S,R)$.
\end{proof}

\begin{paragr}\label{link_homotopy_cat}
Let $\op\Delta\mathit{Sh(\sm/S)}$ be the category of
Nisnevich sheaves of simplicial sets on $Sm/S$.
Morel and Voevodsky defined in \cite{MV} the $\AA^1$-homotopy
theory in $\op\Delta\mathit{Sh(\sm/S)}$. In particular, we have a notion
of $\AA^1$-weak equivalences of simplicial sheaves
that defines a proper model category structure (with the
monomorphisms as cofibrations).
Furthermore, we have a canonical functor
$$\op\Delta\mathit{Sh}(\sm/S)\To\Comp(\sm/S,R)\quad ,
\qquad \X\longmapsto R(\X)$$
which has the following properties; see e.g. \cite{A11,A12}.
\begin{itemize}
\item[(1)] The functor $R$ above preserves colimits.
\item[(2)] The functor $R$ preserves monomorphisms.
\item[(3)] The functor $R$ sends $\AA^1$-weak equivalences
to $\AA^1$-equivalences.
\end{itemize}
We deduce from these properties that the functor $R$
sends homotopy pushout squares of $\op\Delta\mathit{Sh(\sm/S)}$
to homotopy pushout squares of $\Comp(\sm/S,R)$ and induces
a functor
$$R:\mathscr H(S)\To\DMte(S,R)$$
where $\mathscr H(S)$ denotes the localization of $\op\Delta\mathit{Sh(\sm/S)}$
by the $\AA^1$-weak equivalences.

This implies that all the results of \cite{MV} that are formulated
in terms of $\AA^1$-weak equivalences (or isomorphisms
in $\H(S)$) and in terms of homotopy pushout have their counterpart
in $\DMte(S,R)$. We give below the results
we will need that come from this principle.
\end{paragr}

\begin{paragr}\label{defthomspaces}
Let $X$ be a smooth $S$-scheme
 and $\bundle V$ a vector bundle over $X$.
Consider the open immersion $j:\bundle V^\times \rightarrow \bundle V$
of the complement of the zero section of $\bundle V/X$.
We define the \emph{Thom space of $\bundle V$} as the quotient
\begin{equation}\label{defthomspacescoker}
R(\thom\bundle V)
=\mathrm{coker}\big(R(\bundle V^\times)
   \xrightarrow{j_*} R(\bundle V)\big) \ .
\end{equation}
We thus have a short exact sequence of sheaves of $R$-modules
\begin{equation}\label{defthomspacesshort}
0\To R(\bundle V^\times)\To R(\bundle V)\To R(\thom\bundle V)\To 0 \ .
\end{equation}
\end{paragr}

\begin{prop}\label{trivialthom99}
Let $\bundle O^n$ be the trivial vector bundle of dimension $n$
on a smooth $S$-scheme $X$. Then we have a canonical
isomorphism in $\DMte(S,R)$:
$$R(\thom \bundle O^n)\simeq R(X)(n)[2n] \ .$$
\end{prop}

\begin{proof}
This follows from Proposition \ref{motivicsphere}
and from the second statement of
\cite[Proposition 2.17, page 112]{MV}.
\end{proof}

For a given vector bundle $\bundle V$, over a $S$-scheme $X$,
we will denote by $\PP(\bundle V)\To X$ the corresponding projective
bundle.

\begin{prop}\label{projthom}
Let $\bundle V$ be a vector bundle on a smooth $S$-scheme $X$.
Then we have a canonical distinguished triangle in $\DMte(S,R)$
$$R(\PP(\bundle V))\To R(\PP(\bundle V\oplus \bundle O))\To R(\thom\bundle V)
\To R(\PP(\bundle V))[1] \ .$$
\end{prop}

\begin{proof}
This follows from Proposition \ref{motivicsphere}
and from the third statement of
\cite[Proposition 2.17, page 112]{MV}.
\end{proof}

\begin{cor}\label{trivialprojthom}
We have a canonical distinguished triangle in $\DMte(S,R)$
$$R(\PP^n_{S})\To R(\PP^{n+1}_{S})\To R(n+1)[2n+2]\To R(\PP^n_{S})[1] \ .$$
Moreover, this triangle splits canonically for $n=0$ and gives the decomposition
$$R(\PP^1_{S})=R\oplus R(1)[2] \ .$$
\end{cor}

\begin{proof}
This is a direct consequence of Propositions \ref{trivialthom99}
and \ref{projthom}. The splitting of the case $n=0$ comes obviously from
the canonical map from $\PP^1_{S}$ to $S$.
\end{proof}

\begin{paragr}\label{defProjinfty}
The inclusions $\PP^n_{S}\subset\PP^{n+1}_{S}$ allow us to define the
Nisnevich sheaf of sets
\begin{equation}\label{eqdefprojinfty}
\PP^\infty_{S}=\textstyle\limind_{\, n\geq 0}\PP^n_{S}
\end{equation}
We get a Nisnevich sheaf of $R$-modules
\begin{equation}\label{eqdefRprojinfty}
R(\PP^\infty_{S})=\textstyle\limind_{\, n\geq 0}R(\PP^n_{S})\ .
\end{equation}
For a complex $K$ of sheaves of $R$-modules, we define
the hypercohomology of $\PP^\infty_{S}$ with coefficients in $K$
to be
\begin{equation}\label{defhypercohprojinfty}
\hypercoh^i_\nis(\PP^\infty_{S},K)=\Hom_{\Der(\sm/S,R)}(R(\PP^\infty_{S}),K[i]) \ .
\end{equation}
\end{paragr}

\begin{prop}\label{MilnorProj}
There is a short exact sequence
$$\textstyle 0\To\limproj^1_{\, n\geq 0}\hypercoh^{i-1}_\nis(\PP^n_{S},K)
\To\hypercoh^i_\nis(\PP^\infty_{S},K)\To
\limproj_{\, n\geq 0}\hypercoh^i_\nis(\PP^n_{S},K)\To 0 \ .$$
\end{prop}

\begin{proof}
As the filtering colimits are exact in $\mathit{Sh}(\sm/S,R)$
we have an isomorphism
$\smash{\hocolim R(\PP^n_{S})\simeq R(\PP^\infty_{S})}$
in $\Der(\sm/S,R)$. This result is thus a direct application
of the Milnor short exact sequence applied to this homotopy
colimit (see e.g. \cite[Proposition 7.3.2]{Hov}).
\end{proof}

\begin{prop}[Purity Theorem]\label{purity}
Let $i:Z\To X$ a closed immersion of smooth $S$-schemes,
and $U=X-i(Z)$. Denote by $\bundle N_{X,Z}$ the normal vector bundle of $i$.
Then there is a canonical distinguished triangle in $\DMte(S,R)$
$$\xymatrix{R(U)\To R(X)\To R(\thom\bundle N_{X,Z})\To R(U)[1] \ .}$$
\end{prop}

\begin{proof}
This follows from \cite[Theorem 2.23, page 115]{MV}.
\end{proof}

\begin{cor}\label{anmoinszerosphere}
There is a canonical decomposition $R(\AA^n_{S}-\{0\})=R\oplus R(n)[2n-1]$
in $\DMte(S,R)$.
\end{cor}

\begin{proof}
The Purity Theorem and Proposition \ref{trivialthom99}
give a distinguished triangle
$$\xymatrix{R(\AA^n_{S}-\{0\})\To R(\AA^n_{S})\To R(n)[2n]\To R(\AA^n_{S}-\{0\})[1] \ .}$$
But this triangle is isomorphic to the distinguished triangle
$$\xymatrix{R(\AA^n_{S}-\{0\})\To R\To Q[1]\To R(\AA^n_{S}-\{0\})[1]}$$
where $Q$ is the kernel of the obvious
map $R(\AA^n_{S}-\{0\})\To R$, which shows that these triangle
split.
\end{proof}
%

\subsection{Tate spectra}

\begin{paragr}
We want the derived tensor product by $R(1)$
to be an equivalence of categories. As this is not the case
in $\DMte(S,R)$, we will modify
the category $\DMte(S,R)$ and construct
the triangulated category of real motives
${\DMt}(S,R)$ in which this will occur by definition.
For this purpose, we will define the model category of
symmetric Tate spectra. We will give only
the minimal definitions we will need to
work with. We invite the interested reader
to have look at \cite[Section 6]{HCD}
for a more complete account. The main
properties of ${\DMt}(S,R)$ are listed in \ref{proprietesDMt}.

We consider given a category of smooth $S$-schemes $\V$
as in \ref{axiomesVV}.
\end{paragr}

\begin{paragr}
A \emph{symmetric Tate spectrum} (in $\mathit{Sh}(\V,R)$)
is a collection
$E=(E_{n},\s_{n})_{n\geq 0}$, where for each integer $n\geq 0$,
$E_{n}$ is a complex of Nisnevich sheaves on $\V$
endowed with an action of the symmetric group $\mathfrak S_{n}$,
and $\s_n : R(1)\otimes^{}_{R}E_{n}\To E_{n+1}$ is a morphism
of complexes, such that the induced maps obtained by composition
$$R(1)^{\otimes m}\otimes^{}_{R} E_{n}\To R(1)^{\otimes m-1}\otimes^{}_{R} E_{n+1}\To\cdots\To
R(1)\otimes^{}_{R} E_{m+n-1}\To E_{m+n}$$
are $\mathfrak S_{m}\times\mathfrak S_{n}$-equivariant.
We have to define the actions to be precise: $\mathfrak S_{m}$
acts on $R(m)=R(1)^{\otimes m}$ by permutation,
and the action on $E_{m+n}$ is induced by the diagonal
inclusion $\mathfrak S_{m}\times\mathfrak S_{n}\subset \mathfrak S_{m+n}$.
A morphism of symmetric Tate spectra
$u:(E_{n},\s_{n})\To(F_{n},\tau_{n})$ is a collection of $\mathfrak S_{n}$-equivariant maps
$u_{n}:E_{n}\To F_{n}$ such that the squares
$$\xymatrix{
R(1)\otimes^{}_{R} E_{n}\ar[r]^{\s_{n}}\ar[d]_{R(1)\otimes u_{n}}&E_{n+1}\ar[d]^{u_{n+1}}\\
R(1)\otimes^{}_{R} F_{n}\ar[r]_{\tau_{n}}&F_{n+1}
}$$
commute. We denote by $\Spt(\V,R)$ the category of
symmetric Tate spectra.
If $A$ is a complex of sheaves of $R$-modules
on $\V$, we define its \emph{infinite suspension} $\Sigma^\infty(A)$
as the symmetric Tate spectrum that consists of the collection
$(A(n),1_{A(n+1)})_{n\geq 0}$ where $\mathfrak S_{n}$ acts on
$A(n)=R(1)^{\otimes n}\otimes^{}_{R}A$
by permutation on $R(n)=R(1)^{\otimes n}$. This defines
the infinite suspension functor
\begin{equation}\label{definftysusptate}
\Sigma^\infty : \Comp(\V,R)\To\Spt(\V,R)
\end{equation}
This functor has a right adjoint
\begin{equation}\label{defomegainftytate}
\Omega^\infty : \Spt(\V,R)\To\Comp(\V,R)
\end{equation}
defined by $\Omega^\infty(E_{n},\s_{n})_{n\geq 0}=E_{0}$.
According to \cite[6.14 and 6.20]{HCD}
we can define a ($R$-linear) tensor product of symmetric spectra
$E\otimes^{}_{R}F$ satisfying the following properties (and these
properties determine this tensor product up to a canonical isomorphism).
\begin{itemize}
\item[(1)] This tensor product makes the category of symmetric
Tate spectra a closed symmetric monoidal category
with $\Sigma^\infty(R(S))$ as unit.
\item[(2)] The infinite suspension functor \eqref{definftysusptate}
is a symmetric monoidal functor.
\end{itemize}
Say that a map of symmetric Tate spectra $u:(E_{n},\s_{n})\To(F_{n},\tau_{n})$
is a \emph{quasi-isomorphism} if the map $u_{n}:E_{n}\To F_{n}$
is a quasi-isomorphism of complexes of Nisnevich sheaves of $R$-modules
for any $n\geq 0$. We define the \emph{Tate derived category
of $\mathit{Sh}(\V,R)$} as the localization of $\Spt(\V,R)$ by the class
of quasi-isomorphisms. We will write $\Der_{\tate}(\V,R)$ for this
``derived category''. One can check that $\Der_{\tate}(\V,R)$
is a triangulated category (according to \cite[Remark 6.19]{HCD},
this is the homotopy category of a stable model category)
and that the functor induced by $\Sigma^\infty$
is a triangulated functor (because this is a left Quillen
functor between stable model categories).

A symmetric Tate spectrum $E=(E_{n},\s_{n})_{n\geq 0}$ is
a \emph{weak $\Omega^\infty$-spectrum} if for any integer $n\geq 0$,
the map $\s_{n}$ induces an isomorphism
$E_{n}\simeq\derR\sHom(R(1),E_{n+1})$ in $\DMte(S,R)$.
A symmetric Tate spectrum $E=(E_{n},\s_{n})_{n\geq 0}$ is
a \emph{$\Omega^\infty$-spectrum} if it is a weak $\Omega^\infty$-spectrum
and if, for any integer $n\geq 0$, the complex $E_{n}$ is $\V_{\nis}$-local
and $\AA^1$-homotopy invariant.

A morphism of symmetric Tate spectra $u:A\To B$ is a
\emph{stable $\AA^1$-equivalence} if for any weak
$\Omega^\infty$-spectrum $E$, the map
$$u^* : \Hom_{\Der_{\tate}(\V,R)}(B,E)\To\Hom_{\Der_{\tate}(\V,R)}(A,E)$$
is an isomorphism of $R$-modules.

A morphism of Tate spectra is a \emph{stable $\AA^1$-fibration} if it is termwise $\V_{\nis}$-surjective and if its kernel is a
$\Omega^\infty$-spectrum.

A morphism of Tate spectra is a \emph{stable $\V$-cofibration} if it has the left lifting property with respect to the stable $\AA^1$-fibrations
which are also stable $\AA^1$-equivalences. A symmetric Tate spectrum $E$ is \emph{stably $\V$-cofibrant} if the map $0\To E$ is a stable
$\V$-cofibration.
\end{paragr}

\begin{prop}\label{cmfA1stable}
The category of symmetric Tate spectra is a stable proper symmetric
monoidal model category with the stable $\AA^1$-equivalences
as the weak equivalences, the stable $\AA^1$-fibrations as fibrations
and the stable $\V$-cofibrations as cofibrations. The
infinite suspension functor is a symmetric left Quillen functor
that sends the $\AA^1$-equivalences to the stable
$\AA^1$-equivalences. Moreover, the tensor product by
any stably $\V$-cofibrant symmetric Tate spectrum
preserves the stable $\AA^1$-equivalences.
\end{prop}

\begin{proof}
The first assertion is an application of \cite[Proposition 6.15]{HCD}.
The fact that the functor $\Sigma^\infty$ preserves weak
equivalences comes from \emph{loc. cit.}, Proposition  6.18.
The last assertion follows from \emph{loc. cit.}, Proposition 6.35.
\end{proof}

\begin{paragr}\label{proprietesDMt}
The proposition above means the following.

Define the \emph{triangulated category of real
mixed motives} ${\DMt}(S,R)$ as the localization of
the category $\Spt(\V,R)$ by the class of stable
$\AA^1$-equivalences. Then $\DMt(S,R)$ is a triangulated category
with infinite direct sums and products. To be more precise,
any short exact sequence in $\Spt(\V,R)$
gives rise canonically to an exact triangle in $\DMt(S,R)$,
and any distinguished triangle is isomorphic to
an exact triangle that comes from a short exact sequence.
Furthermore, this triangulated category does not depend on the
category $\V$: the category $\V$ is only a technical tool
to define a model category structure that is well
behaved with the tensor product and Nisnevich descent in $\V$.

The infinite suspension functor sends $\AA^1$-equivalences to
stable $\AA^1$-equivalences and thus induces a functor
\begin{equation}\label{DMtsigmainfty}
\Sigma^\infty : \DMte(S,R)\To\DMt(S,R) \ .
\end{equation}
The right adjoint of the infinite suspension functor has a total right
derived functor
\begin{equation}\label{DMtomegainfty}
\derR\Omega^\infty : \DMt(S,R)\To\DMte(S,R) \ .
\end{equation}
For a (weak) $\Omega^\infty$-spectrum $E$, one has
\begin{equation}\label{DMtomegaspectra}
\derR\Omega^\infty(E)=E_{0} \ .
\end{equation}
The tensor product on $\Spt(\V,R)$ has a total left derived functor
\begin{equation}\label{DMttensor}
\DMt(S,R)\times\DMt(S,R)\To\DMt(S,R)
\quad , \qquad (E,F)\longmapsto E\otimes^\derL_{R}F \ .
\end{equation}
If $E$ is stably $\V$-cofibrant, then the canonical map
$E\otimes^\derL_{R}F\To E\otimes^{}_{R}F$ is an isomorphism in $\DMt(S,R)$.
Moreover, the functor \eqref{DMtsigmainfty} is symmetric monoidal.
In particular, for two complexes of Nisnevich sheaves of $R$-modules $A$ and $B$
we have a canonical isomorphism
\begin{equation}\label{DMttensor2}
\Sigma^\infty(A\otimes^\derL_{R}B)\simeq
\Sigma^\infty(A)\otimes^\derL_{R}\Sigma^\infty(B)\ .
\end{equation}
The category $\DMt(S,R)$ has also an internal Hom that we denote by
$\derR\sHom(E,F)$.

We will write $\tR=\Sigma^\infty(R(S))$,
and for a smooth $S$-scheme $X$, we define
$\tR(X)$ to be $\Sigma^\infty(R(X))$.
We also define $\tR(n)=\Sigma^\infty(R(n))$
for $n\geq 0$. Note that $\tR$ is the unit of the (derived) tensor product.
We define the symmetric Tate spectrum $\tR(-1)$ by
the formula $\tR(-1)_{n}=R(n+1)$ with the action of $\mathfrak S_{n}$
defined as the action by permutations on the first $n$ factors of
$R(n+1)=R(1)^{\otimes n}\otimes^{}_{R} R(1)$. The maps $\tR(-1)_{n}\otimes R(1)\To \tR(-1)_{n+1}$
are just the identities. One can check that $\tR(-1)$ is $\V$-cofibrant.
\end{paragr}

\begin{prop}\label{TateinverseDMt}
The object $\tR(1)$ is invertible in $\DMt(S,R)$ and
we have an isomorphism $\tR(-1)\simeq\tR(1)^{-1}$.
In other words, there are isomorphisms
$$\tR(1)\otimes^{\derL}_{R}\tR(-1)\simeq\tR\quad
\text{and}\quad\tR(-1)\otimes^\derL_{R}\tR(1)\simeq\tR \ . $$
\end{prop}

\begin{proof}
This follows from \cite[Proposition 6.24]{HCD}.
\end{proof}

\begin{paragr}
For an integer $n\geq 0$, we define $\tR(-n)=\tR(-1)^{\otimes n}$.
For an integer $n$, and a symmetric Tate spectrum $E$, we define
$$E(n)=E\otimes\tR(n) \ .$$
As $\tR(n)$ is $\V$-cofibrant, the canonical maps
$E\otimes^{\derL}_{R}\tR(n)\To E\otimes^{}_{R}\tR(n)=E(n)$
are isomorphisms in $\DMt(S,R)$.

We will say that a symmetric Tate spectrum
$E=(E_{n},\s_{n})_{n\geq 0}$ is a \emph{weak
$\Omega^\infty$-spectrum} if for any integer $n\geq 0$,
the map $\s_{n}$ induces by adjunction an isomorphism
$$E_{n}\simeq\derR\sHom(R(1),E_{n+1})$$
in $\DMte(S,R)$.
\end{paragr}

\begin{prop}\label{suspDMt}
Let $E$ be a weak $\Omega^\infty$-spectrum.
Then for any integer $n\geq 0$ and any complex of Nisnevich
sheaves of $R$-modules $A$, there is a canonical isomorphism of
$R$-modules
$$\Hom_{\DMt(S,R)}(\Sigma^\infty(A),E(n))\simeq
\Hom_{\DMte(S,R)}(A,E_{n}) \ .$$
In particular, for any smooth $S$-scheme $X$, one has isomorphisms
$$\hypercoh^{i}_{\nis}(X,E_{n})\simeq
\Hom_{\DMt(S,R)}(\tR(X),E(n)[i]) \ .$$
\end{prop}

\begin{proof}
This  is an application of \cite[Proposition 6.28]{HCD}.
\end{proof}

\begin{cor}
A morphism of weak $\Omega^\infty$-spectra $E\To F$
is a stable $\AA^1$-equivalence if and only if
the map $E_{n}\To F_{n}$ is a $\AA^1$-equivalence
for all $n\geq 0$.
\end{cor}

\begin{prop}\label{compactdmtilde}
For any smooth $S$-scheme $X$ and any integer $n$,
$\tR(X)(n)$ is a compact object of $\DMt(S,R)$.
\end{prop}

\begin{proof}
Let $(E_{\lambda})_{\lambda\in\Lambda}$ be a small family of
Tate spectra. We want to show that the map
$$\bigoplus_{\lambda\in\Lambda}\Hom_{\DMt(S,R)}\big(\tR(X)(n),E_{\lambda}\big)
\To\Hom_{\DMt(S,R)}\big(\tR(X)(n),\bigoplus_{\lambda\in\Lambda}E_{\lambda}\big)$$
is bijective. Replacing the spectra $E_{\lambda}$ by the spectra $E_{\lambda}(-n)$,
we can suppose that $n=0$. Furthermore, we can assume that the spectra
$E_{\lambda}$ are weak $\Omega^\infty$-spectra.
As $R(1)$ is a compact object of $\DMte(S,R)$
(this is by definition a direct factor of the object $R(\GG_{m})$
which is compact by \ref{repcompactdmtildeeff}),
it follows from Proposition \ref{A1compactRHomeff}
that $\bigoplus_{\lambda\in\Lambda}E_{\lambda}$
is a weak $\Omega^\infty$-spectrum as well. Therefore, by Proposition \ref{suspDMt},
we have
$$\begin{aligned}
\bigoplus_{\lambda\in\Lambda}\Hom_{\DMt(S,R)}\big(\tR(X),E_{\lambda}\big)&\simeq
\bigoplus_{\lambda\in\Lambda}\Hom_{\DMte(S,R)}\big(R(X),E_{0,\lambda}\big)\\
&\simeq
\Hom_{\DMte(S,R)}\big(R(X),\bigoplus_{\lambda\in\Lambda}E_{0,\lambda}\big)\\
&\simeq
\Hom_{\DMte
(S,R)}\big(R(X),\big(\bigoplus_{\lambda\in\Lambda}E_{\lambda}\big)_{0}\big)\\
&\simeq
\Hom_{\DMt(S,R)}\big(\tR(X),\bigoplus_{\lambda\in\Lambda}E_{\lambda}\big).
\end{aligned}$$
This proves the result.
\end{proof}

\begin{paragr}\label{constantsheafnoneff}
The functor
\begin{equation}\label{constantsheafnoneff1}
\Comp(R)\To\Spt(\V,R)\quad , \qquad
M\longmapsto \Sigma^\infty(M).
\end{equation}
is a left Quillen functor from the projective model structure on $\Comp(R)$
(see \ref{quillenconstantsheaf}) to the model structure of Proposition \ref{cmfA1stable}.
The right adjoint of \eqref{constantsheafnoneff1}
\begin{equation}\label{constantsheafnoneff2}
\Spt(\V,R)\To\Comp(R)\quad , \qquad
M\longmapsto \Gamma(\Omega^\infty(M))
\end{equation}
is a right Quillen functor. The corresponding total right derived functor
is canonically isomorphic to the composed
functor $\derR\Gamma\circ\derR\Omega^\infty$.

For two objects $E$ and $F$ of $\DMt(S,R)$, we define
\begin{equation}\label{constantsheafnoneff3}
\derR\Hom(E,F)=\derR\Gamma\Big(\derR\Omega^\infty\big(\derR\sHom(E,F)\big)\Big).
\end{equation}
For any integer $n$, we have
a canonical isomorphism
\begin{equation}\label{qconstantsheafnoneff4}
H^n\big(\derR\Hom(E,F)\big)\simeq\Hom_{\DMt(S,R)}(E,F[n]).
\end{equation}
\end{paragr}

\subsection{Ring spectra}

\begin{paragr}\label{defringspectrum}
A \emph{ring spectrum} is a monoid object in the category
of symmetric Tate spectra. A ring spectrum is
\emph{commutative} if it is commutative as a monoid
object of $\Spt(\V,R)$.

Given a ring spectrum $\ste$,
one can form the category of left $\ste$-modules.
These are the symmetric Tate spectra $M$ endowed with
a left action of $\ste$
$$\ste\otimes^{}_{R} M\To M$$
satisfying the usual associativity and unit properties.

We denote by $\Spt(\V,\ste)$
the category of left $\ste$-modules. There is a \emph{base change functor}
\begin{equation}\label{basechangestetate}
\Spt(\V,R)\To\Spt(\V,\ste)\quad , \qquad F\longmapsto \ste\otimes^{}_{R}F
\end{equation}
which is a left adjoint of the forgetful functor
\begin{equation}\label{forgetstetate}
\Spt(\V,\ste)\To\Spt(\V,R)\quad , \qquad M\longmapsto M \ .
\end{equation}
If $\ste$ is commutative, the category $\Spt(\V,\ste)$
is canonically endowed with a closed symmetric
monoidal category structure such that the functor \eqref{basechangestetate}
is a symmetric monoidal functor. We denote by $\otimes^{}_{\ste}$
the corresponding tensor product. The unit of this monoidal structure
is $\ste$ seen as an $\ste$-module.

A morphism of $\ste$-modules is
a \emph{stable $\AA^1$-equivalence}
(resp. a \emph{stable $\AA^1$-fibration})
if it is so as a morphism of symmetric Tate spectra.
A morphism of $\ste$-modules is a stable $\V$-cofibration
if it is has the left lifting property with respect to the
stable $\AA^1$-fibrations which are also
stable $\AA^1$-equivalences.
\end{paragr}

\begin{prop}\label{cmftatestemod}
For a given ring spectrum $\ste$, the category
$\Spt(\V,\ste)$ is endowed with a stable proper
model category structure with the
stable $\AA^1$-equivalences as weak equivalences,
the stable $\AA^1$-fibrations as fibrations, and the
stable $\V$-co\-fibra\-tions as cofibrations.
The base change functor \eqref{basechangestetate}
is a left Quillen functor. Moreover, if $\ste$ is commutative,
then this model structure is symmetric monoidal.
\end{prop}

\begin{proof}
See \cite[Corollary 6.39]{HCD}.
\end{proof}

\begin{paragr}\label{defringspectrum2}
Let $\ste$ be a commutative ring spectrum.
We define $\DMt(S,\ste)$ to be the localization
of the category $\Spt(\V,\ste)$ by the class of
stable $\AA^1$-equivalences. It follows from the
proposition above that this category is canonically endowed with
a triangulated category structure. The base change functor
has a total left derived functor
\begin{equation}\label{derbasechangestetate}
\DMt(S,R)\To\DMt(S,\ste)\quad , \qquad F\longmapsto \ste\odmt_{R}F
\end{equation}
which is a left adjoint of the forgetful functor
\begin{equation}\label{derforgettate}
\DMt(S,\ste)\To\DMt(S,R)\quad , \qquad M\longmapsto M \ .
\end{equation}
The forgetful functor \eqref{derforgettate} is conservative
(which means that an object of $\DMt(S,\ste)$
is null if and only if it is null in $\DMt(S,R)$.
There is a derived tensor product
\begin{equation}\label{dertensorstemod}
\DMt(S,\ste)\times\DMt(S,\ste)\To\DMt(S,\ste)
\quad , \qquad (M,N)\longmapsto M\odmt_{\ste}N
\end{equation}
that turns $\DMt(S,\ste)$ into a symmetric monoidal
triangulated category (by applying \cite[Theorem 4.3.2]{Hov}
to the model structure of Proposition \ref{cmftatestemod}).
The derived base change
functor \eqref{derbasechangestetate} is of course
a symmetric monoidal functor. The category
$\DMt(S,\ste)$ also has an internal Hom that we denote by
$\derR\sHom_{\ste}(M,N)$. We thus have the formula
\begin{equation}\label{internalhomste}
\Hom_{\DMt(S,\ste)}(L\otimes^\derL_{\ste}M,N)\simeq
\Hom_{\DMt(S,\ste)}(L,\derR\sHom_{\ste}(M,N)).
\end{equation}

It follows from \ref{constantsheafnoneff} that the functor
\begin{equation}\label{constantsheafste1}
\Comp(R)\To\Spt(\V,\ste)\quad , \qquad
M\longmapsto\ste\otimes^{}_{R}\Sigma^\infty(M).
\end{equation}
is a left Quillen functor. The right derived functor of
its right adjoint is the composition of the forgetful functor \eqref{derforgettate}
with the functor $\derR\Gamma\circ\derR\Omega^\infty$.
For two objects $M$ and $N$ of $\DMt(S,\ste)$, we define
\begin{equation}\label{constantsheafste2}
\derR\Hom_{\ste}(M,N)=
\derR\Gamma\Big(\derR\Omega^\infty\big(\derR\sHom_{\ste}(M,N)\big)\Big).
\end{equation}
For any integer $n$, we have
a canonical isomorphism
\begin{equation}\label{constantsheafste3}
H^n\big(\derR\Hom_{\ste}(M,N)\big)\simeq\Hom_{\DMt(S,\ste)}(M,N[n]).
\end{equation}

For a smooth $S$-scheme $X$, we define
the free $\ste$-module generated by $X$ as
$$\ste(X)=\ste\otimes^{}_{R}\tR(X)=\ste\otimes^{}_{R}\Sigma^\infty(R(X)) \ .$$
As $\tR(X)$ is stably $\V$-cofibrant, the canonical map
$$\ste\otimes^{\derL}_{R}\tR(X)\To\ste\otimes^{}_{R}\tR(X)=\ste(X)$$
is an isomorphism in $\DMt(S,R)$ (hence in $\DMt(S,\ste)$
as well). This implies that for any $\ste$-module $M$,
we have canonical isomorphisms
\begin{equation}\label{calculcohstemod}
\Hom_{\DMte(S,R)}(R(X),\derR\Omega^\infty(M))\simeq
\Hom_{\DMt(S,\ste)}(\ste(X),M) \ .
\end{equation}
Note that as the forgetful functor \eqref{derforgettate} preserves direct sums,
Proposition \ref{compactdmtilde} implies that $\ste(X)(n)$
is a compact object of $\DMt(S,\ste)$ for all smooth $S$-scheme $X$
and integer $n$.
\end{paragr}


\section{Modules over a Weil spectrum}

\emph{From now on, we assume the given scheme $S$ is regular.}

Let $\V$ be a full subcategory of smooth $S$-schemes satisfying the
hypothesis of \ref{axiomesVV}.
We also fix a field of characteristic zero $\KK$ called the \emph{field of
coefficients}.

\subsection{Mixed Weil theory}

\begin{paragr} \label{introducing_substratum}
Let $E$ be a complex of presheaves of $\KK$-module on $\V$ which has
the Brown-Gersten property and is $\AA^1$-homotopy invariant.
From Proposition \ref{bgtilde}, the first property means that
$H^i(E(X))=\hypercoh^i_\nis(X,E_\nis)$ for any scheme $X$ in $\V$,
and any integer $i$.
The second property implies the complex of sheaves $E_\nis$ is
quasi-isomorphic to a $\AA^1$-local complex.
In the sequel we will write $E$ for the corresponding object of $\DMte(S,\KK)$.
Whence we obtain, for any smooth $S$-scheme $X$ in $\V$,
a canonical isomorphism
$$
H^i(E(X))=\Hom_{\DMte(S,\KK)}(\KK(X),E[i]) \, .
$$
Suppose moreover that $E$ has a structure of a presheaf of
commutative differential graded $K$-algebras.
This structure corresponds to morphisms of presheaves
$$
E \otimes_K E \xrightarrow \mu E, \qquad
 \KK \xrightarrow \eta E
$$
satisfying the usual identities (corresponding to the associativity
and commutativity properties of the multiplication $\mu$ and to the the
fact $\eta$ is a unit).
Applying the associated Nisnevich sheaf functor, we obtain
 in $\Comp(\Sh(S,\KK))$ the following morphisms
$$
E_\nis \otimes_K E_\nis \xrightarrow \mu E_\nis, \qquad
 \KK \xrightarrow \eta E_\nis.
$$
As the sheafifying functor and the tensor product over $\KK$ are exact,
these morphisms indeed induce a commutative monoid structure
on $E$, as an object of $\DMte(S,\KK)$.
\end{paragr}

\begin{paragr} \label{df_cohomology}
Consider now a merely commutative monoid object $E$ of $\DMte(S,\KK)$.

Let us denote by $\mu:E \odmte_\KK E \To E$
and $\eta:\KK \To E$
respectively the multiplication and the unit maps.

If $M$ is an object of $\DMte(S,\KK)$,
we set $H^i(M,E)=\Hom_{\DMte(S,\KK)}(M,E[i])$.
For two objects $M$ and $N$ of $\DMte(S,\KK)$,
we define the external cup product
$$H^p(M,E) \otimes_\KK H^q(N,E)
 \To H^{p+q}(M \odmte_\KK N,E)$$
as follows. Considering two morphisms
$\alpha : M\To E[p]$ and $\beta : N\To E[q]$ in $\DMte(S,\KK)$,
we define a map $\alpha \otimes_\mu \beta$ as
the composite
$$
M \odmte_{\KK} N
 \xrightarrow{\alpha \odmte_{\KK} \beta} E[p]\odmte_{\KK}E[q]
 \xrightarrow{\mu[p+q]}E[p+q]
$$
that is the expected product of $\alpha$ and $\beta$.

For a smooth $S$-scheme $X$, we simply write $H^i(X,E)=H^i(\KK(X),E)$.
We can consider the diagonal embedding $X \To X \times_S X$ which induces a comultiplication $\delta_*:\KK(X) \To \KK(X) \odmte_\KK \KK(X)$. This
allows to define as usual a `cup product' on $H^*(X,E)$ by the formula
$$\alpha \, . \,  \beta=(\alpha \otimes_\mu \beta) \circ \delta_* \ .$$ We
will always consider $H^*(X,E)$ as a graded $\KK$-algebra with this cup product.

We introduce 
the following axioms~:
\begin{enumerate}
\item[W1] \textit{Dimension}.--- $H^{i}(S,E)\simeq
\begin{cases}
\KK&\text{if $i=0$,}\\
0&\text{otherwise.}
\end{cases}$
\item[W2] \textit{Stability}.---
$\mathrm{dim}_{\KK} H^i(\GG_{m},E)=
\begin{cases}
1&\text{if $i=0$ or $i=1$,}\\
0&\text{otherwise.}
\end{cases}$
\item[W3] \textit{K\"unneth formula}.---
For any smooth $S$-schemes $X$ and $Y$, the exterior cup product
induces an isomorphism
$$
\bigoplus_{p+q=n}H^p(X,E) \otimes_\KK H^q(Y,E)
 \overset{\sim}{\To} H^{n}(X\times_{S}Y,E) \ .
$$
\item[W3$^\prime$] \textit{Weak K\"unneth formula}.---
For any smooth $S$-scheme $X$, the exterior cup product
induces an isomorphism
$$
\bigoplus_{p+q=n}H^p(X,E) \otimes_\KK H^q(\GG_{m},E)
 \overset{\sim}{\To} H^{n}(X\times_{S}\GG_{m},E) \ .
$$
\end{enumerate}
\end{paragr}

\begin{paragr}\label{defstabstructure}
Under assumptions W1 and W2, we will call any non zero element
$c \in H^1(\GG_m,E)$ a \textit{stability class}. 
Note that such a class corresponds to a non trivial map
$$c:\KK(1)\To E$$
in $\DMte(S,\KK)$.
In particular, if $E$ is a presheaf
of commutative diffential graded $\KK$-algebras
which has the $\text{B.-G.}$-property
and is $\AA^1$-homotopy invariant, then such a stability class
can be lifted to an actual map of complexes of presheaves.
Such a lift will be called a \emph{stability structure} on $E$.

Remark that, in the formulation of axiom W3 (resp. W3$^\prime$)
we might require the K\"unneth formula to hold only
for $X$ and $Y$ (resp. $X$) in $\V$: as any smooth $S$-scheme is locally in $\V$
for the Nisnevich topology, this apparently weaker condition
implies the general one by a Mayer-Vietoris argument.
\end{paragr}

\begin{df}
A \emph{mixed Weil theory}\footnote{
In what follows, we will prove this terminology is not usurpated: a consequence
of the main results of this paper is that, when $S$ is the spectrum
of a field $k$, the restriction of the functor $H^*(\, . \, ,E)$
to smooth and projective $k$-schemes is a Weil cohomology in the sense
defined in \cite{andremotifs}.}
is a presheaf $E$ of commutative differential graded $\KK$-algebras
on $\V$ which has the 
Brown-Gersten property
(or equivalently the excision property, see \ref{bgtilde}), 
is $\AA^1$-homotopy invariant,
and satisfies the properties W1, W2 and W3 stated above.

A \emph{stable theory} is a presheaf $E$ of commutative differential graded $\KK$-algebras
on $\V$ which has the
Brown-Gersten property, is $\AA^1$-homotopy invariant,
and satisfies the properties W1, W2 and W3$^\prime$.
\end{df}

\begin{paragr}\label{defsteWeil}
Any stable theory $E$ gives rise canonically
to a commutative ring spectrum $\ste$, as explained below.
The idea to define the spectrum $\ste$
consists essentially to consider a weighted version of $E$
(this should be clearer considering the comments given in \ref{defweightedWeil}).

Let $\Hom^*(\KK(1),E)$ be the complex of maps
of complexes of sheaves from $\KK(1)$ to $E$
(the category $\Comp(\V,\KK)$ is naturally enriched in
complexes of $\KK$-vector spaces). As $\KK(1)$ is $\V$-cofibrant
and as $E$ is fibrant with respect to the model category structure
of Proposition \ref{A1cmfeff}, we have for any integer $i$
\begin{equation}\label{cohHomK1E}
H^{i}(\Hom^*(\KK(1),E))=
\begin{cases}
H^1(\GG_{m},E)&\text{if $i=0$,}\\
0&\text{otherwise.}
\end{cases}
\end{equation}
Consider the constant sheaf of complexes on $\V$
\begin{equation}\label{deflefshetzste}
L=\Hom^*(\KK(1),E)_{S}
\end{equation}
associated to the complex $\Hom^*(\KK(1),E)$. We can now define a
symmetric Tate spectrum $\ste=(\ste_{n},\sigma_{n})_{n\geq 0}$
as follows. Put first $\ste_{n}=\sHom(L^{\otimes n},E)$ (here $\sHom$
stands for the internal Hom in the category $\Comp(\V,\KK)$).
We have a canonical map
$$L=\Hom^*(\KK(1),E)_{S}\To\sHom(\KK(1),E)$$
which gives a map
\begin{equation}\label{defsigmansteintermediaire0}
\KK(1)\otimes_{\KK}L\To E\, ,
\end{equation}
and tensoring with $\sHom(L^{\otimes n},E)$
gives a map
\begin{equation}\label{defsigmansteintermediaire1}
\KK(1)\otimes_{\KK}L\otimes_{\KK}\sHom(L^{\otimes n},E)
\To E\otimes_{\KK}\sHom(L^{\otimes n},E)\, .
\end{equation}
The product on $E$ induces a canonical action of $E$ on $\sHom(L^{\otimes n},E)$~:
\begin{equation}\label{defsigmansteintermediaire2}
E\otimes_{\KK}\sHom(L^{\otimes n},E)\To\sHom(L^{\otimes n},E)\, .
\end{equation}
The composition of \eqref{defsigmansteintermediaire1} and \eqref{defsigmansteintermediaire2}
finally leads to a morphism
\begin{equation}\label{defsigmansteintermediaire3}
\KK(1)\otimes_{\KK}L\otimes_{\KK}\sHom(L^{\otimes n},E)
\To \sHom(L^{\otimes n},E)\, .
\end{equation}
The map $\sigma_{n}:\ste_{n}(1)\To\ste_{n+1}$
is defined at last from \eqref{defsigmansteintermediaire3} by transposition,
using the isomorphism $\sHom(L,\sHom(L^{\otimes n},E)\simeq\sHom(L^{\otimes (n+1)},E)$.
The action of $\mathfrak S_{n}$ on $\ste_{n}$ is by permutation of factors in
$L^{\otimes n}$. Note that the fact $\ste$ is well defined relies heavily on the fact
$E$ is commutative as a differential graded algebra.
We define in the same spirit a commutative ring spectrum structure on $\ste$.
The unit map $\tKK\To\ste$ is determined by a sequence of maps
$\eta_{n}:\KK(n)\To\ste_{n}$. The map $\eta_{0}$ is of course the
unit of $E$, 
and the rest of the sequence is then obtained easily by induction:
if $\eta_{n-1}$ is defined, then $\eta_{n}$ is obtained as the composition
$$\KK(n)\xrightarrow{\eta_{n-1}(1)}\ste_{n-1}(1)\xrightarrow{\sigma_{n-1}}\ste_{n}\, .$$
The multiplication of $\ste$ is determined by
maps $\mu_{m,n}:\ste_{m}\otimes_{\KK}\ste_{n}\To\ste_{m+n}$ which are
defined by composition of the obvious maps below.
$$\sHom(L^{\otimes m},E)\otimes_{\KK}\sHom(L^{\otimes m},E)
\to\sHom(L^{\otimes(m+n)},E\otimes_{\KK}E)\to\sHom(L^{\otimes(m+n)},E)$$
\end{paragr}

\begin{prop}\label{orientomegaspectre}
Let $E$ be a stable theory.
The associated commutative ring spectrum $\ste$
is a weak $\Omega^\infty$-spectrum, and
there is a canonical
isomorphism $E\simeq\derR\Omega^\infty(\ste)$ in $\DMte(S,\KK)$.
In other words, for any sheaf of complexes $M$,
we have canonical isomorphisms
$$\Hom_{\DMte(S,\KK)}(M,E)
\simeq\Hom_{\DMte(S,\KK)}(M,\derR\Omega^\infty(\ste))
\simeq\Hom_{\DMt(S,\KK)}(\derL\Sigma^\infty(M),\ste) .$$
Furthermore, any stability structure on $E$ defines
an isomorphism $\ste(1)\simeq \ste$ in $\DMt(S,\KK)$.
\end{prop}

\begin{proof}
It follows from \eqref{cohHomK1E} that
the complex $\Hom^*(\KK(1),E)$ is quasi-isomorphic to
the constant sheaf associated to the vector space $H^1(\GG_{m},E)$.
As a consequence, the constant sheaf $L$ is a $\V$-cofibrant\footnote{
As we work with a field of coefficients $\KK$, any constant
sheaf of complexes of vector spaces is $\V$-cofibrant.}
complex which is (non canonically) isomorphic
to $\KK$ in $\DMte(S,\KK)$.
Taking into account that $E$ is
quasi-isomorphic, as a presheaf, to its fibrant
replacement in the model structure of Proposition \ref{A1cmfeff},
we also have a canonical isomorphism in $\DMte(S,\KK)$
$$\ste_{n}=\sHom(L^{\otimes n},E)\simeq\derR\sHom(L^{\otimes n},E)\, .$$
Hence we can get a non canonical isomorphism $\ste_{n}\simeq E$
which corresponds to the choice of a generator $c$ of $H^1(\GG_{m},E)$.
Under such an identification, the structural maps
$$\ste_{n}\To\derR\sHom(\KK(1),\ste_{n+1})$$
all correspond in $\DMte(S,\KK)$ to the map
$$\tau_{c} : E\To\derR\sHom(\KK(1),E)$$
induced by transposition of the map $E(1)\To E$, obtained as the cup product of
the identity of $E$ and of the map $\KK(1)\To E$
coming from the chosen generator $c$.
The weak K\"unneth formula and the stability axiom thus imply that the
map $\tau_{c}$ above is an isomorphism in $\DMte(S,\KK)$.
This proves that $\ste$ is indeed a weak $\Omega^\infty$-spectrum.
The reformulation of this assertion comes directly from Proposition \ref{suspDMt}.

Consider now a stability structure $c:\KK(1)\To E$ on $E$.
We have to define a morphism of symmetric Tate spectra $u:\ste(1)\To \ste$,
which corresponds to $\mathfrak S_{n}$-equivariant maps
commuting with the $\sigma_{n}$'s
$$u_{n}:\sHom(L^{\otimes n},E)(1)\To\sHom(L^{\otimes n},E)\, .$$
Such a map $u_{n}$ is determined by a map
$$v_{n}:L^{\otimes n}\otimes_{\KK}\sHom(L^{\otimes n},E)(1)\To E\, .$$
We already have an evaluation map twisted by $\KK(1)$
$$L^{\otimes n}\otimes_{\KK}\sHom(L^{\otimes n},E)(1)\To E(1)\, ,$$
so that to define $v_{n}$, we are reduced to define a map
$$E(1)\To E\, ;$$
this is obtained as the cup product of the identity of $E$ with the given map $c$.
The fact that $u$ is an isomorphism in $\DMte(S,\KK)$
comes again from the stability axiom and from the weak K\"unneth formula.
\end{proof}

\begin{paragr}\label{defweightedWeil}
Given a stable theory and its associated commutative ring spectrum $\ste$,
for a smooth $S$-scheme $X$ and two integers $p$ and $q$, we define the
\emph{$q^{\text{th}}$ group of cohomology of $X$ of twist $p$ with coefficients in $\ste$}
to be
\begin{equation}\label{defweightedWeil1}
H^q(X,\ste(p))=\Hom_{\DMt(S,\KK)}(\tKK(X),\ste(p)[q])\, .
\end{equation}
We obviously have
\begin{equation}\label{defweightedWeil2}
H^q(X,E)=H^q(X,\ste)\, ,
\end{equation}
and more generally, if $p\geq 0$,
$H^*(X,\ste(p))$ is just the Nisnevich hypercohomology
of $X$ with coefficients in the sheaf of complexes $\sHom(L^{\otimes p},E)$.
Hence for any integer $p$, any choice of a generator of $H^1(\GG_{m},E)$
determines a non canonical (but still functorial) isomorphism
$H^q(X,E)\simeq H^q(X,\ste(p))$.

We also define complexes
\begin{equation}\label{defweightedWeil3}
\derR\Gamma(X,\ste(p))=\derR\Hom_{\KK}(\tKK(X),\ste(p))\simeq\derR\Hom_{\ste}(\ste(X),\ste(p))
\end{equation}
and we get by definition
\begin{equation}\label{defweightedWeil4}
\HH^q(\derR\Gamma(X,\ste(p)))=H^q(X,\ste(p))\, .
\end{equation}
\end{paragr}

\begin{paragr}
Let $X$ be a smooth $S$-scheme and
$\alpha:\ste(X) \To \ste(p)[i]$
and
$\beta:\ste(X) \To \ste(q)[j]$
be morphisms of $\ste$-modules, corresponding
to cohomological classes.
The cup product of $\alpha$ and $\beta$ over $X$
then corresponds to a map of $\ste$-modules
$$\alpha \, . \, \beta:\ste(X)\To\ste(p+q)[i+j]$$
defined as the composition
$$
\ste(X) \xrightarrow{\delta_*}\ste(X\times_{S}X)\simeq
 \ste(X) \otimes^\derL_\ste \ste(X) \xrightarrow{\alpha \otimes^\derL_\ste \beta}
  \ste(p)[i] \otimes^\derL_\ste \ste(q)[j]\simeq\ste(p+q)[i+j].
$$
\end{paragr}

\subsection{First Chern classes}

We assume a stable theory $E$ is given. We will consider
its associated commutative ring spectrum $\ste$ (\ref{defsteWeil}), and the
corresponding cohomology groups (\ref{defweightedWeil}).

\begin{paragr}\label{blablastabclass}
Recall we have a canonical decomposition $\KK(\GG_m)=\KK \oplus \KK(1)[1]$
in the category $\DMte(S,\KK)$.
The unit map $\KK\To\ste$ determines by twisting and shifting a map
\begin{equation}\label{blablastabclass1}
c:\KK(1)[1]\To\ste(1)[1]\, .
\end{equation}
The morphism \eqref{blablastabclass1}, seen in $\DMte(S,\KK)$
corresponds to a non trivial cohomology class
in $\Hom_{\DMt(S,\KK)}(\KK(1)[1],\ste(1)[1])=H^1(\GG_{m},\ste(1))$.

We also have a decomposition $\KK(\PP^1_S)=\KK \oplus \KK(1)[2]$,
so that \eqref{blablastabclass1} also corresponds to a cohomology class $c$ in
$H^2(\PP^1_S,\ste(1))$
that will be called the \emph{canonical orientation of $\ste$.}

Note also that the decomposition $\KK(\PP^1_S)=\KK \oplus \KK(1)[2]$
and the weak K\"unneth formula
implies the K\"unneth formula holds with respect to
products of type $\PP^1_S\times_S X$
for Nisnevich cohomology with coefficients in $E$.
We will still refer to this as the `weak K\"unneth formula'.
\end{paragr}

\begin{lm}\label{calculcohtateobj}
For any integer $n \geq 0$, the graded vector space $H^*(\KK(n),E)$
is (non canonically) isomorphic to $\KK$ concentrated in degree zero.
\end{lm}
\begin{proof}
The case $n=0$ is precisely W1. Assume $n\geq 1$.
We can begin by a choice of a stability structure on $E$,
which defines, using W2 and the weak K\"unneth formula,
an isomorphism in $\DMte(S,\KK)$:
$$E\simeq\derR\sHom(\KK(1),E)\, .$$
This gives
$$\begin{aligned}
\derR\Hom_{\DMte(S,\KK)}(\KK(n),E)
&\simeq\derR\Hom_{\DMte(S,\KK)}(\KK(n-1),\derR\sHom(\KK(1),E))\\
&\simeq\derR\Hom_{\DMte(S,\KK)}(\KK(n-1),E)\, .
\end{aligned}$$
We conclude by induction on $n$.
\end{proof}

For any integer $1 \leq n \leq m$,
we let $\iota_{n,m}:\PP^n_S \To \PP^m_S$ be the embedding
given by $(x_0:\ldots:x_n)\longmapsto(x_0:\ldots:x_n:0:\ldots:0)$.
\begin{lm}
\label{lm:pre_proj_bundle_th}
For any integer $n>0$, the cohomology group $H^*(\PP^n_S,E)$
is concentrated in degrees $i$ such that $i$ is even and $i \in [0,2n]$.

For any integer $0 \leq n \leq m$,
$\iota_{n,m}^*:H^*(\PP^m_S,E) \To H^*(\PP^n_S,E)$
is an isomorphism in degrees $i \in [0,2n]$.
\end{lm}
\begin{proof}
The case where $m=1$ is already known (\ref{blablastabclass}).
The remaining assertions follow then by induction
from the canonical distinguished triangle
$$\KK(\PP^{n-1}_S) \xrightarrow{\iota_{n-1,n}} \KK(\PP^n_S)
 \To \KK(n)[2n] \To \KK(\PP^{n-1}_S)[1]$$
in $\DMte(S,\KK)$ (see Corollary \ref{trivialprojthom}), and
using Lemma \ref{calculcohtateobj}.
\end{proof}

\begin{lm}
\label{lm:premutation_induce_iso_proj}
Let $n \geq 2$ and
$\sigma$ be a permutation of the set $\{0,\dots,n\}$.

Consider the morphism $\sigma:\PP^n_S \To \PP^n_S,
(x_0:\ldots :x_n) \longmapsto (x_{\sigma(0)}:\ldots:x_{\sigma(n)})$.

Then $\sigma^*:H^2(\PP^n_S,E) \To H^2(\PP^n_S,E)$
is the identity.
\end{lm}
\begin{proof}
We consider first the case $n\geq 3$.
We can assume $\sigma$ is the transposition $(n-1,n)$.
Then $\sigma \iota_{1,n}=\iota_{1,n}$ and the claim
follows from the preceeding lemma.

It remains to prove the case $n=2$.
Let $\sigma$ a transposition of $\{0,1,2\}$.
There is then a transposition $\tau$ of $\{0,1,2,3\}$
such that $\iota_{2,3}\sigma=\tau\iota_{2,3}$.
As we already know that $\tau$ induces the identity
in degree $2$ cohomology. By applying Lemma \ref{lm:pre_proj_bundle_th},
we see that the map $\iota_{2,3}$ induces an isomorphism in
degree $2$ cohomology as well. We thus get, by functoriality,
$\sigma^*\iota^*_{2,3}=\iota^*_{2,3}\tau^*=\iota^*_{2,3}$,
with $\iota^*_{2,3}$ invertible, which ends the proof.
\end{proof}

\begin{paragr}\label{defSegre}
For any integer $n,m \geq 0$, we will consider the Segre embedding
\begin{equation}\label{defSegre1}
\sigma_{n,m}:\PP^n_S \times \PP^m_S \To \PP^{nm+n+m}_S
\end{equation}
and the $n$-fold Segre embedding
\begin{equation}\label{defSegre2}
\sigma^{(n)}:(\PP^1_S)^n \To \PP^{2^n-1}_S \, .
\end{equation}
\end{paragr}

\begin{prop} \label{prop:fundamental_chern_classes}
There exists a unique sequence $(c_{1,n})_{n>0}$ of cohomology classes,
with $c_{1,n} \in H^2(\PP^n_S,\ste(1))$, such that~:
\begin{enumerate}
\item[(i)] $c_{1,1}=c$ is the canonical orientation of $\ste$;
\item[(ii)] for any integers $1 \leq n \leq m$,
 $\iota_{n,m}^*(c_{1,m})=c_{1,n}$.
\end{enumerate}
Moreover, the following formulas hold~:
\begin{enumerate}
\item[(iii)] for any integer $n>0$, $c_{1,n}^n \neq 0$ and
$c_{1,n}^{n+1}=0$;
\item[(iv)] for any integers $n,m>0$,
$\sigma_{n,m}^*(c_{1,nm+n+m})=\pi^*_n(c_{1,n})+\pi^*_m(c_{1,m})$,
where $\pi_n$ and $\pi_m$ denote the projections
from $\PP^n_S\times_S\PP^m_S$ to $\PP^n_S$ and $\PP^m_S$
respectively.
\end{enumerate}
\end{prop}
\begin{proof}
The unicity statement is clear from \ref{lm:pre_proj_bundle_th}.

Let $n \geq 2$ be an integer and consider the embedding
$$
p:\PP^1_S \To (\PP^1_S)^n \ , \quad (x:y) \longmapsto ((x:y),(0:1),\dots,(0:1)) \, .
$$
The morphism $\iota_{1,2^n-1}$ factors\footnote{This factorization might
hold eventually only up to a permutation of coordinates in $\PP^{2^n-1}_S$
(depending on the choices made to define the Segre embeddings), but this will
be harmless by Lemma \ref{lm:premutation_induce_iso_proj}.} as
$$
\PP^1_S \xrightarrow p (\PP^1_S)^n
 \xrightarrow{\sigma^{(n)}} \PP^{2^n-1}_S
$$
which induces in cohomology
$$
H^2(\PP^{2^n-1}_S,\ste(1)) \xrightarrow{\sigma^{(n)*}} H^2((\PP^1_S)^n,\ste(1))
 \xrightarrow{p^*} H^2(\PP^1_S,\ste(1)).
$$
Let $t$ be the unique class in $H^2(\PP^{2^n-1}_S,\ste(1))$ which is sent to
$c$ by the isomorphism $\iota^*_{1,2^n-1}$ (in degree $2$
cohomology).
Using the weak K\"unneth formula and Lemma \ref{calculcohtateobj},
we obtain a decomposition
$$H^2((\PP^1_S)^n,\ste(1))=\KK .u_1 \oplus \cdots \oplus \KK .u_n$$
(where $u_i$ stands for the the generator corresponding to $c$
in the $i$th copy of $H^2(\PP^1_S,\ste(1))$ in $H^2((\PP^1_S)^n,\ste(1))$).
The map $p^*$ is then the $\KK$-linear map sending $u_1$ to $c$ and $u_i$ to $0$ for $i>1$.
This implies $\sigma^{(n)*}(t)=u_1+\lambda_2.u_2+\cdots+\lambda_n.u_n$
for $\lambda_i \in \KK$.

But from Lemma \ref{lm:premutation_induce_iso_proj}, the class $t$
is invariant under permutations of the coordinates of $\PP^{2^n-1}_S$.
This implies $\sigma^{(n)*}(t)$ is invariant under permutations
of the factors of $(\PP^1_S)^n$ which gives $\lambda_2=\dots=\lambda_n=1$.

Then $\sigma^{(n)*}(t^n)=n!\, u_1\cdots u_n$ which is non zero by the weak K\"unneth
formula and Lemma \ref{calculcohtateobj} (remember $\KK(\PP^1_{S})=\KK\oplus\KK(1)[2]$).
This implies $t^n \neq 0$.

We put $c_{1,n}=\iota_{n,2^n-1}^*(t)$ so that we have
$\iota^*_{1,n}(c_{1,n})=c_{1,1}$.

As $\iota_{n,2^n-1}$ induces an isomorphism
in degree less than $2n$, we see that $c_{1,n}^n \neq 0$
and $c_{1,n}^{n+1}=0$ (indeed, $u_i^2=0$ for any $i=1,\ldots,n$).

The existence of the sequence and property (iii) are then proved.
Moreover, by the unicity statement,
we see that the class $t$ introduced in the preceeding reasoning
is $t=c_{1,2^n-1}$.

Property (iv) follows from the relation $\sigma^{(n)*}(c_{1,2^n-1})=u_1+\ldots+u_n$
and from property (ii).
\end{proof}

\begin{paragr}\label{paragr:representing_pic}
Remember from \ref{defProjinfty} the ind-scheme $\PP^\infty_S$
defined by the tower of inclusions
$$
\PP^1_S \To\cdots
 \To  \PP^n_S \xrightarrow{{\ \iota_n \ }}
 \PP^{n+1}_S\To\cdots
$$
We set
$H^q(\PP^\infty_S,\ste(p))
=\Hom_{\DMt(S,\KK)}(\derL\Sigma^\infty(\KK(\PP^\infty_k)),\ste(p)[q])$.
\end{paragr}

\begin{cor}\label{cohstePinfty}
The sequence $(c_{1,n})_{n >0}$ of the previous proposition gives an element $\mathfrak{c}$ of $H^2(\PP^\infty_S,\ste(1))$ which induces an
isomorphism of graded $\KK$-algebras
$$
\KK[[\mathfrak{c}]]=\prod_{n\geq 0}H^{2n}(\PP^\infty_S,\ste(n)).
$$
\end{cor}
\begin{proof}
Using Proposition \ref{suspDMt},
we have the Milnor short exact sequences (\ref{MilnorProj})
$$\textstyle 0\To\limproj^1_{\, n\geq 0}{H}^{i-1}(\PP^n_{S},\ste(p))
\To{H}^i(\PP^\infty_{S},\ste(p)) \xrightarrow{(*)} \limproj_{\, n\geq 0}{H}^i(\PP^n_{S},\ste(p))\To 0 \ .$$
The $\limproj^1$ of a constant functor is
null, and thus Lemma \ref{lm:pre_proj_bundle_th}
implies $(*)$ is an isomorphism. The corollary then follows directly
from the previous proposition.
\end{proof}

\begin{paragr}\label{deffirstchern}
The sequence $(c_{1,n})_{n>0}$ induces a morphism in $\DMt(S,\KK)$
\begin{equation}\label{deffirstchern1}
\mathfrak c:\Sigma^\infty\KK(\PP^\infty_S) \To \ste(1)[2]\, .
\end{equation}
As a consequence, using the functor
$\KK:\mathscr H (S) \To \DMte(S,\KK)$
introduced in \ref{link_homotopy_cat},
for any smooth $S$-scheme $X$, we obtain
a canonical map
$$
[X,\PP^\infty_S]
\To \Hom_{\DMt(S,\KK)}(\Sigma^\infty\KK(X),\Sigma^\infty\KK(\PP^\infty_S))\, .
$$
The map \eqref{deffirstchern1} induces a map
$$\Hom_{\DMt(S,\KK)}(\Sigma^\infty\KK(X),\Sigma^\infty\KK(\PP^\infty_S))
\To \Hom_{\DMt(S,\KK)}(\Sigma^\infty\KK(X),\ste(1)[2])\, .$$
As the base scheme $S$ is regular, it follows
from \cite[Proposition 3.8, p. 138]{MV} that
we have a natural bijection
\begin{equation}\label{deffirstchern3}
[X,\PP^\infty_S]\simeq\pic(X) \, .
\end{equation}
We have thus associated to \eqref{deffirstchern1} a canonical map
\begin{equation}\label{deffirstchern2}
c_1:\pic(X) \To H^2(X,\ste(1))
\end{equation}
called the \emph{first Chern class}.
Note that this map is just defined as a map of pointed sets (it
obviously preserves zero), but we have
more structures, as stated below.
\end{paragr}

\begin{prop} \label{cor:first_chern_class}
The map \eqref{deffirstchern2} introduced above is a morphism of abelian groups
which is functorial in $X$ with respect to pullbacks.
\end{prop}
\begin{proof}
Functoriality is obvious.

The family of Segre embeddings
$\sigma_{n,m}:\PP^n_S \times \PP^m_S \To \PP^{nm+n+m}_S$
defines a morphism of ind-schemes
$$
\PP^\infty_S \times_S \PP^\infty_S \xrightarrow{\sigma}
 \PP^\infty_S.
$$
which in turn defines an H-group structure on $\PP^\infty_S$
as an object of $\mathscr H(S)$, and put a group structure
on $[X,\PP^\infty_S]$.

Let $\lambda$ (resp. $\lambda'$, $\lambda''$) be the canonical dual
line bundle on $\PP^\infty_S$ (resp. the two canonical dual line
bundles on $\PP^\infty_S \times_S \PP^\infty_S$).
An easy computation gives
$$
\sigma^*:\pic(\PP^\infty_S)
 \To \pic(\PP^\infty_S \times_S \PP^\infty_S)\, , \
  \lambda \longmapsto \lambda'\otimes\lambda''
$$
which implies the preceding group structure
coincides with the usual group structure on the Picard group
via \eqref{deffirstchern3}.

But similarly, from property (iv) in \ref{prop:fundamental_chern_classes},
we obtain the map
$$
\sigma^*:H^*(\PP^\infty_S;E)=\KK[[\mathfrak c]]
 \To H^*(\PP^\infty_S \times_S \PP^\infty_S,E)=\KK[[\mathfrak c',\mathfrak c'']]\, , \
 \mathfrak c \longmapsto \mathfrak c'+ \mathfrak c''
$$
which gives precisely the result we need.
\end{proof}

\subsection{Projective bundle theorem and cycle class maps}

\begin{paragr}
We consider given a stable theory $E$ and its associated commutative ring spectrum $\ste$ (\ref{defsteWeil}).

Recall from \ref{defringspectrum} the symmetric monoidal category
of $\ste$-modules is endowed with a
notion of stable $\AA^1$-weak equivalence. The associated localized
category is denoted by $\DMt(S,\ste)$; see \ref{defringspectrum2}.
We have an adjoint pair of functors
$$
\ste\otimes^\derL_{\KK}(\, -\, ):\DMt(S,\KK) \rightleftarrows \DMt(S,\ste):\mathcal U
$$
where $\ste\otimes^\derL_{\KK}(\, -\, )$ is the total left derived fuinctor
of the free $\ste$-module functor and $\mathcal U$
forgets the $\ste$-module structure. The functor $\mathcal U$
is conservative.

The study of the cohomology theory associated to $\ste$
follows obviously from the study of the triangulated category
$\DMt(S,\ste)$. It follows from the existence of the first Chern class \eqref{deffirstchern2}
that the results and constructions of \cite{Deg6} apply to $\DMt(S,\ste)$.
This leads to the following classical results.
\end{paragr}

\begin{paragr}\label{deftrivmapbundleste}
Consider now a vector bundle $\bundle V$ of rank $n$ over a smooth $S$-scheme $X$,
$p:\bundle \PP(\bundle V) \To X$ be the canonical projection.
Consider the first Chern class \eqref{deffirstchern2}
$$
\pic(\PP(\bundle V)) \xrightarrow{c_1}
H^2(\PP(\bundle V),\ste(1))=\Hom_{\DMt(S,\ste)}(\ste(\PP(\bundle V)),\ste(1)[2]).
$$
Thus the canonical dual invertible sheaf $\lambda=\bundle O(-1)$ on $\PP(\bundle V)$
induces a morphism of $ \ste$-modules
\begin{equation}\label{deftrivmapbundleste1}
c_1(\lambda):\ste(\PP(\bundle V)) \To \ste(1)[2].
\end{equation}
This defines a map
\begin{equation}\label{deftrivmapbundleste2}
a_{\PP(\bundle V)}: \ste(\PP(\bundle V))\To  \bigoplus_{i=0}^{n-1} \ste(X)(i)[2i]
\end{equation}
 by the formula
\begin{equation}\label{deftrivmapbundleste3}
a_{\PP(\bundle V)}=\sum_{i=0}^{n-1} c_1(\lambda)^i  . \,  p_* \, .
\end{equation}
\end{paragr}

\begin{prop}[Projective Bundle Formula]\label{projbundlethmstemod}
The map \eqref{deftrivmapbundleste2} is an isomorphism in $\DMt(S,\ste)$.
\end{prop}

\begin{proof}
See \cite[Theorem 3.2]{Deg6}.
\end{proof}

\begin{paragr}\label{defhigherChernste}
We can now come to the classical definition of Chern classes. Let
$X$ be a smooth scheme and $\bundle V/X$ a vector bundle of rank $n$.
Let $\lambda$ be the canonical dual line bundle on $\PP(\bundle V)$.

By virtue of Proposition \ref{projbundlethmstemod},
the canonical map
$$
\bigoplus_{i=0}^{n-1} H^{2j-2i}(X,\ste(j-i)).\, c_1(\lambda)^i \To H^{2j}(\PP(\bundle V),\ste(j))
$$
is an isomorphism for all $j$.

Define the Chern classes $c_i(\bundle V)$ of $\bundle V$ in $H^{2i}(X,\ste(i))$ by the relations
\begin{itemize}
\item[(a)] $\sum_{i=0}^n p^*c_i(\bundle V) \, . \, c_1(\lambda)^{n-i}=0$;
\item[(b)] $c_0(\bundle V)=1$;
\item[(c)] $c_i(\bundle V)=0$ for $i>n$.
\end{itemize}
These Chern classes are functorial with respect to pullback
and extends the first Chern classes given by \eqref{deffirstchern2}.
Following a classical argument, we obtain the additivity for these
Chern classes (see \cite[Lemma 3.13]{Deg6}):
\end{paragr}

\begin{lm}
\label{additivity_chern_classes}
Let $X$ be a smooth $S$-scheme, and consider an exact sequence
$$
0 \To \bundle V' \To \bundle V \To \bundle V'' \rightarrow 0
$$
of vector bundles over $X$. Then $c_r(\bundle V)=\sum_{i+j=r} c_i(\bundle V') \, . \, c_j(\bundle V'')$.
\end{lm}

\begin{paragr}\label{defLefshetzembed}
Let $\bundle V$ be a vector bundle of rank $n+1$ over a smooth $S$-scheme $X$.

For any integer $r \in [0,n]$, we define the \emph{Lefschetz embedding}
\begin{equation}\label{defLefshetzembed1}
\mathfrak l_r(\PP(\bundle V)):\ste(X)(r)[2r]\To\ste(\PP(\bundle V))
\end{equation}
as the composition
\begin{equation}\label{defLefshetzembed2}
\ste(X)(r)[2r] \xrightarrow{(*)} \oplus_{i=0}^n \ste(X)(i)[2i]
 \xrightarrow{a_{\PP(\bundle V)}^{-1}} \ste(\PP(\bundle V))
\end{equation}
where $(*)$ is the obvious embedding.

Recall the morphism
\begin{equation}\label{defLefshetzembed3}
\pi:\ste(\PP(\bundle V \oplus \bundle O)) \To \ste (\thom\bundle V)
\end{equation}
appearing from the distinguished triangle of Proposition \ref{projthom}.
\end{paragr}

\begin{lm}\label{trivialthomstemod}
Let $\bundle V$ be a vector bundle of rank $n$ over a smooth $S$-scheme $X$,
and $P=\PP(\bundle V \oplus \bundle O)$ be its projective completion.
The composite morphism
$$
\ste(X)(n)[2n]
 \xrightarrow{\mathfrak l_n(P)} \ste(\PP(\bundle V \oplus \bundle O))
  \xrightarrow{\pi} \ste (\thom\bundle V)
$$
is an isomorphism.
\end{lm}

\begin{proof}
Simply use the distinguished triangle of Proposition \ref{projthom},
the definition of $a_P$ and the compatibility of the first Chern
class with pullback.
\end{proof}

\begin{paragr}\label{definvtrivthomste}
In the situation of the preceding lemma,
we will denote by
$$\mathfrak p_{\bundle V}:\ste(\thom\bundle V) \To \ste(X)(n)[2n]$$
the inverse isomorphism of $\pi\, \mathfrak l_n(P)$.
\end{paragr}

\begin{prop}[Purity Theorem]\label{purityste}
Let $i:Z\To X$ be a closed immersion of smooth $S$-schemes
of pure codimension $n$. We denote by $j:U=X-Z\To X$
the complementary open immersion.
There is a canonical distinguished triangle
$$\ste(U)\xrightarrow{j_{*}}\ste(X)\xrightarrow{i^*}\ste(Z)(n)[2n]
\xrightarrow{\partial_{X,Z}}\ste(U)[1]$$ in $\DMt(S,\ste)$.
\end{prop}

\begin{proof}
By applying the triangulated functor
$\ste\otimes^\derL_{\KK}\derL\Sigma^\infty(-)$ to the distinguished triangle of
Proposition \ref{purity}, we obtain a distinguished triangle
$$\ste(U)\To\ste(X)\To\ste(\thom\bundle N_{X,Z})\To\ste(U)[1]$$
with $\bundle N_{X,Z}$ the normal vector bundle of the immersion $i$.
We conclude using the isomorphism $\mathfrak p_{\bundle N_{X,Z}}$ introduced above.
\end{proof}

\begin{paragr}\label{defGysintriangle}
The distinguished triangle of the proposition is called
the \emph{Gysin triangle associated to the pair $(X,Z)$},
and the map $i^*$ is called the \emph{Gysin morphism associated to $i$}.
The precise definition of the Gysin triangle and its main
functorialities are described and proved in \cite{Deg6}
in a more general context. We recall the main properties we will need below.
\end{paragr}



\begin{prop}\label{functGysinste}
Given a cartesian square of smooth $S$-schemes
$$
\xymatrix{
T\ar^j[r]\ar_g[d] & Y\ar^f[d] \\
Z\ar^i[r] & X
}
$$
where $i$ and $j$ are closed immersions of pure codimension $n$,
we have the following commutative diagram.
$$\xymatrix{
{\ste}(Y)\ar^/-3pt/{j^*}[r]\ar^{f_*}[d]
 & {\ste}(T)(n)[2n]\ar^/-5pt/{\partial_{Y,T}}[r]\ar^{g_*(m)[2m]}[d]
 & {\ste}(Y-T)[1]\ar^{h_*}[d] \\
{\ste}(X)\ar^/-3pt/{i^*}[r]
 & {\ste}(Z)(n)[2n]\ar^/-5pt/{\partial_{X,Z}}[r]
 & {\ste}(X-Z)[1]
}$$
\end{prop}

\begin{proof}
See \cite[Proposition 4.10]{Deg6}.
\end{proof}

\begin{rem}
In cohomology, the Gysin morphism introduced above induces
a morphism
$i_*:H^*(Z,\ste) \rightarrow H^{*+2n}(X,\ste(n))$.

The commutativity of the left hand square above gives the usual
projection formula in the transversal case~:
\begin{equation} \label{transversal_case}
f^*i_*=j_*g^*.
\end{equation}

Note that, as explained in \cite[Corollary 4.11]{Deg6}, the above projection formula
implies easily the following projection formula for $\ste$-modules
$$
(1_Z \, . \,  i_*) \circ i^*=i^*  \, . \,  1_X (n)[2n]\, ,
$$
which implies the usual projection formula in cohomology:
\begin{equation} \label{product_projection_formula}
\forall \, a \in H^i(Z,\ste), \forall \, b \in H^i(X,\ste), \quad i_*(a \, . \, i^*b)=i_*(a) \, . \, b \in H^{i+2n}(X,\ste(n))\ .
\end{equation}
\end{rem}

\begin{df} \label{df:fund_class}
Let $i:Z \To X$ be a codimension $n$ closed immersion
between smooth $S$-schemes.
We put $\eta_X(Z)=i_*(1)$ as an element in $H^{2n}(X,\ste(n))$
and call $\eta_X(Z)$ the \emph{fundamental class of $Z$ in $X$}.
\end{df}
Note that this fundamental class corresponds uniquely to a morphism of
$\ste$-modules
$$
\eta_X(Z):\ste(X) \To \ste(n)[2n]
$$
which we refer also as the fundamental class.

\begin{rem} \label{fund_class&section}
Suppose that $i$ admits a retraction $p:X \To Z$.
Let $p_Z:Z \To S$ be the structural morphism.
Then the projection formula gives
$$
i^*=(p_{Z*}  \, . \,  p_*i_*) \circ i^*=(p_{Z*}i^*) \, . \,  p_*
 =\eta_X(Z)  \, . \,  p_*.
$$
The Gysin morphism in this case is completely determined by the
fundamental class $\eta_Z(X)$.
\end{rem}

\begin{ex} \label{ex:thom_class}
Let $\bundle V$ be a vector bundle of rank $n$ over $X$,
and $P=\PP(\bundle V\oplus \bundle O)$ be its projective completion.
Let $p:P \To X$ be the canonical projection
and $i:X \To P$ the zero section.
If $\lambda$ denotes the canonical dual line bundle on $P$, the
\emph{Thom class of $\bundle V$} is the cohomological class in $H^{2n}(P,\ste(n))$
$$
t(\bundle V)=\sum_{i=0}^n p^*(c_i(\bundle V))  \, . \, c_1(\lambda)^{n-i}.
$$
By a purity argument (see \cite[4.14]{Deg6}), one gets
 $\eta_P(X)=t(\bundle V)$.

Suppose now given a section $s$ of $\bundle V/X$ transversal
to the zero section $s_0$ of $\bundle V/X$. Put $Y=s^{-1}(s_0(X))$ and
consider the pullback square
$$
\xymatrix{
Y\ar^i[r]\ar_j[d] & X\ar^s[d] \\ X\ar^{s_0}[r] & P.
}
$$
From the projection formula and the identities
$s^*p^*c_i(\bundle V)=c_i(\bundle V)$ and $s^*c_1(\lambda)=0$,
we obtain
$$
\eta_X(Y)=i_*j^*(1)=s^*s_{0*}(1)=s^*(t(\bundle V))=c_n(\bundle V).
$$
\end{ex}

Following \cite[Proposition 4.16]{Deg6}, we also obtain the excess
intersection formula~:

\begin{prop}
Consider a cartesian square of smooth schemes
$$
\xymatrix{
T\ar^j[r]\ar_g[d] & Y\ar^f[d] \\
Z\ar^i[r] & X
}
$$
where $i$ and $j$ are closed immersions of respective codimension $n$ and $m$. Let $e=n-m$ and put $\xi=g^*{\bundle N_{X,Z}}/\bundle N_{Y,T}$ as a
$T$-vector bundle. Let $c_e(\xi)$ be the $e^{\text{th}}$ Chern class of $\xi$.

Then, the following square is commutative.
$$
\xymatrix{
{\ste}(Y)\ar^/-3pt/{j^*}[r]\ar_{f_*}[d]
 & {\ste}(T)(n)[2n]\ar^/-5pt/{\partial_{Y,T}}[r]
 \ar[d]^{c_e(\xi)  . \, g_*(n)[2n]}   
 & {\ste}(Y-T)[1]\ar^{h_*}[d] \\
{\ste}(X)\ar^/-3pt/{i^*}[r]
 & {\ste}(Z)(n)[2n]\ar^/-5pt/{\partial_{X,Z}}[r]
 & {\ste}(X-Z)[1]
}
$$
\end{prop}

\begin{rem}
In particular, we obtain the excess intersection formula in cohomology.
$$
\forall \, a \in H^*(Z,\ste), \quad f^*i_*(a)=j_*(c_e(\xi)\, . \, g^*(a))
$$
This also implies the self intersection formula
for a closed immersion $i:Z \rightarrow X$ of codimension $n$
between smooth schemes.
$$
\forall \, a \in H^*(Z,\ste), \quad i^*i_*(a)=c_n(\bundle N_{X,Z}) \, . \, a
$$
\end{rem}

\begin{paragr}
Let $\bundle V$ be a vector bundle of rank $n$ over $X$,
$p:\bundle V \To X$ the canonical projection,
and
$i:X \To \bundle V$ the zero section. Note that
$p_*:\ste(\bundle V) \To \ste(X)$ is an isomorphism
(by a Mayer-Vietoris argument, one can suppose that
$\bundle V$ is trivial, so that this follows from $\AA^1$-homotopy invariance).
Hence $i_*:\ste(X) \To \ste(\bundle V)$ is the reciprocal
isomorphism.
The self intersection formula implies $\eta_{\bundle V}(X)=p^*c_n(\bundle V)$
in $H^{2n}(\bundle V,\ste(n))$.
Thus, we obtain the computation of the Gysin morphism associated
with the zero section~:
$p_*i^*=c_n(\bundle V) \, . \, 1_X:\ste(X) \rightarrow \ste(X)(n)[2n]$.

We deduce from that the \emph{Euler long exact sequence in cohomology}
$$
\to H^{r-2n}(X,\ste(-n)) \xrightarrow{c_n(\bundle V).} H^r(X,\ste)
 \xrightarrow{q^*} H^r(\bundle V-X,\ste)
  \to H^{r-2n+1}(X,\ste(-n)) \to
$$
\end{paragr}

\begin{thm}
\label{thm:assoc}
Consider a cartesian square of smooth $S$-schemes
$$
\xymatrix{
Z\ar^k[r]\ar_l[d] & Y'\ar^j[d] \\
Y\ar^i[r] & X
}
$$
such that $i$,$j$,$k$,$l$ are closed immersions of respective pure
codimension $n$, $m$, $s$, $t$. We consider the open immersions
$i':Y-Z \rightarrow X-Y'$, $j':Y'-Z \rightarrow X-Y$
and we put put $d=n+s=m+t$.
Then the following diagram is commutative.
$$
\xymatrix{
\ste(X)\ar^{j^*}[r]\ar_{i^*}[d]
 & \ste(Y')(m)[2m]\ar^/-2pt/{\partial_j}[r]\ar^{k^*(m)[2m]}[d]
 & \ste(X-Y')[1]\ar^{{i'}^*[1]}[d] \\
\ste(Y)(n)[2n]\ar_{l^*[2n]}[r]
 & \ste(Z)(d)[2d]\ar^/-7pt/{\partial_l(n)[2n]}[r]
     \ar^{\partial_k(m)[2m]}[d]
 & \ste(Y-Z)(n)[2n][1]\ar^{\partial_{i'}[1]}[d] \\
& \ste(Y-Z)(m)[2m][1]\ar_{-\partial_{j'}[1]}[r]
 & \ste(X-Y \cup Y')[2]
}
$$
\end{thm}

\begin{proof}
This is an application of \cite[Theorem 4.32]{Deg6}.
\end{proof}

\begin{rem}
Indeed, the commutativity of the first two squares
asserts the functoriality of the Gysin triangle with respect
to the Gysin morphism of a closed immersion. The next square
is an associativity result for residues. This theorem also
ensures the compatibility of Gysin morphisms with
tensor product of $\ste$-modules (this will ensure
the compatibility of Gysin morphisms with cup product in
cohomology).
\end{rem}

\begin{paragr}\label{defcherncharactste}
As a conclusion, we have proved in particular the axioms of
Grothendieck for the existence of a cycle class map
(cf. paragraph 2 of \cite{Gro2}): A1
is proved in \ref{projbundlethmstemod},
A2 in example \ref{ex:thom_class}, A3 in \ref{thm:assoc},
and A4 in \ref{product_projection_formula}.
Moreover, the projection fomula \ref{transversal_case} implies
the axiom A5 of 5 in \emph{loc. cit.}. Hence, following the method of \cite{Gro2}
and the theory of $\lambda$-operations, 
we obtain for any smooth $S$-scheme $X$,
a unique homomorphism of rings
\begin{equation}\label{defcherncharactste1}
\mathit{ch}:K_0(X)_\QQ \To H^{2*}(X,\ste(*))
\end{equation}
which is natural in $X$ and such that for any line bundle $\lambda$ on $X$,
the identity below holds.
\begin{equation}\label{defcherncharactste2}
\mathit{ch}([\lambda])=\sum_i \frac 1 {i!} c_1(\lambda)^i
\end{equation}
\end{paragr}

\begin{paragr}\label{recollection}
Let $\mathit{SH}(S)$ be the $\PP^1$-stable homotopy
category of schemes over $S$; we refer to \cite{Jard2,DRO2,pizero}
for different (but equivalent) definitions of $\mathit{SH}(S)$.
According to \cite[5.2]{pizero}\footnote{In \cite{pizero}, $\DMt(S,R)$
is defined using non symmetric spectra instead of symmetric spectra.
But it follows from Voevodsky's Lemma \cite[Lemma 3.13]{Jard2}
and from \cite[Theorem 10.1]{Hov3} that the two definitions lead to
equivalent categories.},
there is a canonical symmetric monoidal triangulated functor
\begin{equation}\label{compfctSHandDMt}
\tilde{R}:\mathit{SH}(S)\To\DMt(S,R)
\end{equation}
that preserves direct sums. It is essentially defined
by sending $\Sigma^{\infty}_{\PP^1}(X_{+})$ to the Tate
spectrum $\Sigma^\infty(R(X))=\tR(X)$
for any smooth $S$-scheme $X$.

For $R=\QQ$, the functor
\eqref{compfctSHandDMt}
induces an equivalence of categories\footnote{The functor \eqref{compfctSHandDMt}
and the equivalence of triangulated  categories \eqref{compfctSHandDMteqrat}
hold without the regularity assumption on $S$.}
\begin{equation}\label{compfctSHandDMteqrat}
\mathit{SH}(S)_{\QQ}\simeq\DMt(S,\QQ),
\end{equation}
where $\mathit{SH}(S)_{\QQ}$ denotes the localization of
$\mathit{SH}(S)$ by the rational equivalences (that are the maps
inducing an isomorphism of stable motivic homotopy groups up to torsion);
see e.g. \cite[Remark 4.3.3 and 5.2]{pizero}. Hence there is no
essential difference to work with $\mathit{SH}(S)_{\QQ}$ or
with $\DMt(S,\QQ)$, which allows to apply here results proved in $\mathit{SH}(S)_{\QQ}$.
In particular, by virtue of  \cite[Definition \textsc{iv}.54]{thriou},
there exists an object $\BGLQ$ in $\DMt(S,\QQ)$ which
represents algebraic $K$-theory\footnote{This is reasonable
because we assumed $S$ to be regular: $K$-theory is homotopy
invariant only for regular schemes.}. Hence, for any smooth $S$-scheme $X$
and any integer $n$, we have
\begin{equation}\label{defkthrat}
\Hom_{\DMt(S,\QQ)}(\tQQ(X)[n],\BGLQ)=K_{n}(X)_{\QQ}\, .
\end{equation}
Moreover, Riou defines for any integer $k$ a morphism
\begin{equation}\label{defadams}
\adams k : \BGLQ\To\BGLQ
\end{equation}
which induces the usual Adams operation on $K$-theory
(see Definition \textsc{iv}.59 of \emph{loc. cit.}).
For an $S$-scheme $X$, define
\begin{equation}\label{motcohbeilinson1}
K_{q}(X)^{(p)}_{\QQ}=\{ x\in K_{q}(X)_{\QQ}\, | \, \adams k (x)=k^px \
\text{for all $k\in\ZZ$}\}
\end{equation}
Recall that \emph{Beilinson motivic cohomology} is defined by
the following formula (see \cite{beilinson}).
\begin{equation}\label{motcohbeilinson2}
\HH^q_{\mathcyr{B}}(X,\QQ(p))=K_{2p-q}(X)^{(p)}_{\QQ}
\end{equation}
Remember from \cite[Expos\'e X, Theorem 5.3.2]{SGA6} that we have
\begin{equation}\label{motcohbeilinson3}
\HH^2_{\mathcyr{B}}(X,\QQ(1))=\pic(X)_{\QQ}\, .
\end{equation}
By virtue of \cite[Theorem \textsc{iv}.72]{thriou}, there exists for each integer $p$,
a projector $\pi_{p}:\BGLQ\To\BGLQ$ such that if we denote
by $\motcoh^{(p)}$ the image of $\pi_{p}$, the canonical map
\begin{equation}\label{motcohbeilinson4}
\bigoplus_{p\in\ZZ}\motcoh^{(p)}\To\BGLQ
\end{equation}
is an isomorphism\footnote{It follows from \cite[Theorem \textsc{v}.31]{thriou}
that this is the unique decomposition of $\BGLQ$ which lifts
the Adams decomposition of $K$-groups \eqref{motcohbeilinson1}.}.
Furthermore, we have
\begin{equation}\label{motcohbeilinson5}
\Hom_{\DMt(S,\QQ)}(\tQQ(X),\motcoh^{(p)}[q])=\HH^{2p+q}_{\mathcyr{B}}(X,\QQ(p))
\end{equation}
for any smooth $S$-scheme $X$. The periodicity theorem
for algebraic $K$-theory \cite[Proposition \textsc{iv}.2]{thriou}
now translates into the existence of canonical isomorphisms
\begin{equation}\label{motcohbeilinson6}
\motcoh^{(0)}(p)[2p]\simeq\motcoh^{(p)}\, .
\end{equation}
In the sequel of this paper,
we will write simply $\motcoh$ for $\motcoh^{(0)}$.
The object $\motcoh$ is called the \emph{Beilinson motivic cohomology
spectrum}.
\end{paragr}

\begin{thm}\label{higherChernste}
There exists a canonical isomorphism in $\Der(\KK)$:
$$\derR\Hom(\motcoh,\ste)\simeq\derR\Gamma(S,\ste)\, .$$
In particular, we have
$$\Hom_{\DMt(S,\QQ)}(\motcoh,\ste[i])\simeq\HH^{i}(S,\ste)=
\begin{cases}
\KK&\text{if $i=0$,}\\
0&\text{otherwise.}
\end{cases}$$
Moreover, there is a unique morphism
$\mathit{cl}_{\mathcyr{B}}:\motcoh\To\ste$ 
inducing the Chern character \eqref{defcherncharactste1}.
This is the unique morphism from $\motcoh$ to $\ste$
which preserves the unit.
\end{thm}

\begin{proof}
This is a rather straightforward application of the nice results
and methods of Riou in \cite{thriou}. The main remark is that
the theory of Chern classes allows to compute the
cohomology of Grassmanians (e.g. following the method
of \cite{Gro2bis}), which in turn shows that we can apply
\cite[Theorem \textsc{iv}.48]{thriou}. Hence using corollary \ref{cohstePinfty},
we see that the arguments to prove \cite[Theorem \textsc{v}.31]{thriou}
can be followed \emph{mutatis mutandis} to give
the expected computation\footnote{The proof of \cite[Theorem \textsc{v}.31]{thriou}
works over any regular base scheme, and holds if we replace motivic cohomology
by any oriented $\QQ$-linear cohomology theory}.

Note that, given a map $\motcoh\To\ste$, we get
in particular morphisms
$$K_{0}(X)^{(n)}_{\QQ}=\Hom_{\DMt(S,\QQ)}(\tQQ(X),\motcoh(n)[2n])
\To\HH^{2n}(X,\ste(n))\, .$$
Hence there is at most one map
$\motcoh\To\ste$ inducing the Chern character \eqref{defcherncharactste1}.
The fact that \eqref{defcherncharactste1} determines
a map $\motcoh\To\ste$ comes
from \cite[Lemma \textsc{iii}.26 and Theorem \textsc{iv}.11]{thriou}
applied to $\ste$.
\end{proof}

\begin{paragr}\label{higherChernste2}
The preceding theorem allows to produce cycle class maps
in the case where the base $S$ is the spectrum of a field $k$.

Let $\HQ$ be Voevodsky's motivic cohomology spectrum (see e.g. \cite{MZ,coniveau}).
According to \cite[Proposition \textsc{v}.36]{thriou}, the
Chern character (which, according to \cite[Section 2.6]{thriou},
is the unique map which preserves the unit)
$\mathit{ch}:\BGLQ\To\HQ$
factorizes through $\motcoh$. Furthermore, it can be shown
that the map $\motcoh\To\HQ$ is an isomorphism in $\DMt(S,\QQ)$:
this follows from the coniveau spectral sequence of the $K$-theory
spectrum, which degenerates rationaly; see \cite{FrSus,coniveau}).
In particular, we get isomorphisms
\begin{equation}\label{motcohbeilinson9}
H^q_{\mathcyr{B}}(X,\QQ(p))\simeq\HH^q(X,\QQ(p))\, .
\end{equation}
We obtain a solid commutative diagram
\begin{equation}\begin{split}
\xymatrix{
\HQ\ar@{-->}[rr]^{\mathit{cl}}&&\ste\\
&\motcoh\ar[ul]_\simeq\ar[ur]^{\mathit{cl}_\mathcyr{B}}&\\
&\BGLQ\ar[uul]^{\mathit{ch}}\ar[u]_{}\ar[uur]_{\mathit{ch}}&
}\end{split}\end{equation}
which defines the cycle class map
\begin{equation}\label{VoecycleclassmapHQste}
\mathit{cl}:\HQ\To\ste\, .
\end{equation}
It follows from Theorem \ref{higherChernste} and from
\cite[Theorem \textsc{v}.31]{thriou} that
$\mathit{cl}$ is the unique map which preserves the unit.
It induces functorial maps (the genuine cycle class maps)
\begin{equation}\label{cycleclassmapsste}
H^q(X,\QQ(p))\To\HH^q(X,\ste(p))\, .
\end{equation}
These cycle class maps are completely determined by
the fact they are functorial and compatible with
cup products and with first Chern classes (this is proved
by applying \cite[Lemma \textsc{iii}.26 and Theorem \textsc{iv}.11]{thriou}).

%
\end{paragr}

\subsection{Gysin morphisms}

\begin{paragr}\label{gysin}
We still consider given a stable theory $E$ and its associated
commutative ring spectrum $\ste$ (\ref{defsteWeil}).

We will now introduce the Gysin morphism of a projective
morphism between smooth $S$-schemes in the setting
of $\ste$-modules (which corresponds to push forward
in cohomology), and recall some of its main properties.

Let $f:Y \To X$ be a projective morphism of codimension $d\in\ZZ$
between smooth $S$-schemes.

Let us choose a factorisation of $f$ into $Y \xrightarrow i \PP^n_{X}
\xrightarrow p X$, where $i$ is a closed immersion of pure codimension
$n+d$, the map $p$ being the canonical projection.

We define the \emph{Gysin morphism
associated to $f$} in $\DMt(S,\ste)$
\begin{equation}\label{gysin1}
f^*:\ste(X)\To\ste(Y)(d)[2d]
\end{equation}
as the following compositum
\begin{equation}\label{gysin2}
f^*=\left\lbrack\ste(X)(n)[2n] \xrightarrow{\mathfrak l_n(\PP^n_{X})} \ste(\PP^n_{X})
 \xrightarrow{i^*} \ste(Y)(n+d)[2(n+d)] \right\rbrack (-n)[-2n].
\end{equation}
\end{paragr}
One can show $f^*$ is independent of the chosen factorization;
see \cite[Lemma 5.11]{Deg6}.

\begin{prop}\label{gysinfctste1}
Consider $Z \xrightarrow g Y \xrightarrow f X$ be
projective morphisms, of codimension $n$ and $m$ respectively,
between smooth $S$-schemes.
Then the following triangle commutes.
$$\xymatrix@C=10pt@R=15pt{
\ste(X)\ar[rr]^(.4){(fg)^*}\ar[dr]_{f^*}&&\ste(Z)(n+m)[2(n+m)]\\
&\ste(Y)(m)[2m]\ar[ur]_{{g^*(n)[2n]}}&
}$$
\end{prop}

\begin{proof}
See \cite[Proposition 5.14]{Deg6}.
\end{proof}

\begin{prop}\label{gysinfctste2}
Consider a cartesian square of smooth $S$-schemes
$$\xymatrix{T\ar_q[d]\ar^g[r] & Z\ar^p[d] \\ Y\ar^f[r] & X}$$
such that $f$ (resp. $g$) is a projective morphism of codimension
$n$ (resp. $m$). Let $\xi$ be the excess bundle over $T$ associated to
that square, and let $e=n-m$ be its rank (cf. \cite[proof of Proposition 6.6]{Ful}).
Then,  $f^*p_*=\big(c_e(\xi) \, . \,  q_*(m)[2m]\big) g^*$ as maps from
$\ste(Z)$ to $\ste(Y)(n)[2n]$.
\end{prop}

\begin{proof}
See \cite[Proposition 5.17]{Deg6}.
\end{proof}

\begin{prop}\label{gysinfctste3}
Consider a cartesian square of smooth $S$-schemes
$$
\xymatrix{
T\ar^j[r]\ar_g[d] & Y\ar^f[d] \\
Z\ar^i[r] & X
}
$$
such that $f$ and $g$ are projective morphism
of respective relative codimension $p$ and $q$, and
such that $i$ and $j$ are
closed immersion of respective codimension $n$ and $m$.
Denote by $h:Y-T \To X-Z$ the morphism
induced by $f$. Then the following square is commutative
(in which the two arrows ${\partial_{X,Z}}$ and $\partial_{Y,T}$
are the one appearing in the obvious Gysin triangles).
$$\xymatrix@C=15pt@R=20pt{
\ste(T)(m+q)[2m+2q]\ar^/2pt/{\partial_{Y,T}(p)[2p]}[rr]&&\ste(Y-T)(p)[2p+1] \\
\ste(Z)(n)[2n]\ar^{\partial_{X,Z}}[rr]\ar_{g^*(n)[2n]}[u]&&\ste(X-Z)[1]\ar_{h^*[1]}[u]
}$$
\end{prop}

\begin{proof}
See \cite[Proposition 5.15]{Deg6}.
\end{proof}

\subsection{Poincar\'e duality}

\begin{paragr}\label{defdualizableobjects}
We first recall the abstract definition of duality in monoidal
categories.
Let $\C$ be a symmetric monoidal category.
We let $\unit$ and $\otimes$ denote respectively the unit object and the
tensor product of $\C$.
An object $X$ of $\C$ is said to have a \emph{strong dual}
if there exists an object $\dual{X}$ of $\C$ and two maps
$$\eta:\unit\To\dual{X}\otimes{X}\quad\text{and}\quad\e:{X}\otimes\dual{X}\To\unit$$
such that the following diagrams commute.
\begin{equation}\begin{split}
\xymatrix{
X\ar[r]^(.3){X\otimes\eta}\ar[dr]_{1_X}&X\otimes\dual{X}\otimes X\ar[d]^{\e\otimes X}&
&\dual{X}\ar[r]^(.35){\eta\otimes\dual{X}}\ar[dr]_{1_{\dual{X}}}&\dual{X}\otimes X\otimes\dual{X}\ar[d]^{\dual{X}\otimes\e}\\
&X&&&\dual{X}}
\end{split}\end{equation}
For any objects $Y$ and $Z$ of $\C$, we then have a canonical bijection
\begin{equation}
\Hom_\C(Z\otimes X,Y)\simeq\Hom_\C(Z,\dual{X}\otimes Y).
\end{equation}
In other words, $\dual{X}\otimes Y$ is in this case the internal $\Hom$
of the pair $(X,Y)$ for any $Y$. In particular, such a strong dual,
together with the maps $\e$ and $\eta$, is unique
up to a unique isomorphism. It is clear that for any symmetric
monoidal category $\D$ and any monoidal functor $F:\C\To\D$,
if $X$ has a strong dual $\dual{X}$, then $F(X)$ has a strong dual canonically
isomorphic to $F(\dual{X})$. If $\dual{X}$ is a strong dual of $X$, then
$X$ is a strong dual of $\dual{X}$. \\
\indent Let $\T$ be a closed symmetric monoidal
triangulated category\footnote{We just mean that the category $\T$
is endowed with a symmetric monoidal structure and
with a triangulated category structure, such that for any object $X$ of $\T$,
the functor $Y\longmapsto X\otimes Y$ is triangulated.}.
Denote by $\sHom$ its internal Hom.
For any objects $X$ and $Y$ in $\T$
the evaluation map
$$X\otimes\sHom(X,\unit)\To\unit$$
tensored with the identity of $Y$ defines by adjunction a map
\begin{equation}\label{evtordue}
\sHom(X,\unit)\otimes Y\To \sHom(X,Y).
\end{equation}
The object $X$ has a strong dual if and only if this map is an isomorphism
for all objects $Y$ in $\T$, and in this case, we have
$\dual{X}=\sHom(X,\unit)$: this follows from the fact that, essentially
by definition, $X$ has a strong dual if and only if
there exists an object $X^\vee$ in $\T$, such that the functor
$X^\vee\otimes(-)$ is right adjoint to the functor $(-)\otimes X$,
so that $X^\vee\otimes(-)$ has to be canonically isomorphic
to the functor $\sHom(X,-)$ (the canonical isomorphism being precisely \eqref{evtordue}).
For $Y$ fixed, the map \eqref{evtordue}
is a morphism of triangulated functors. Hence the objects $X$
such that \eqref{evtordue} is an isomorphism form a full
triangulated subcategory of $\T$. In other words,
the full subcategory $\T_\du$ of $\T$ that consists of the objects
which have a strong dual is a thick triangulated subcategory of $\T$.\\
\indent If moreover $\T$ has small sums, then to say that
any object of $\T_\du$ is compact is equivalent to say that the unit $\unit$
is compact. This is proved as follows. Suppose that $\unit$ is compact,
and let $X$ be an object of $\T$ which has a strong dual $\dual{X}$.
Then for any small family $(Y_{\lambda})_{\lambda\in\Lambda}$
of objects of $\T$, we get the following identifications.
$$\begin{aligned}
\bigoplus_{\lambda\in\Lambda}\Hom_{\T}\big(X,Y_{\lambda}\big)&\simeq
\bigoplus_{\lambda\in\Lambda}\Hom_{\T}\big(\unit,\dual{X}\otimes Y_{\lambda}\big)\\
&\simeq \Hom_{\T}\big(\unit,\bigoplus_{\lambda\in\Lambda}(\dual{X}\otimes Y_{\lambda})\big)\\
&\simeq \Hom_{\T}\big(\unit,\dual{X}\otimes \bigoplus_{\lambda\in\Lambda}Y_{\lambda}\big)\\
&\simeq \Hom_{\T}\big(X,\bigoplus_{\lambda\in\Lambda}Y_{\lambda}\big).
\end{aligned}$$
The converse is obvious. Note that it can happen that a compact
object of $\T$ doesn't have any strong dual; a counter-example
can be found in \cite{riou}. We will produce another counter-example
below, as a consequence of a comparison theorem: for any complete discrete valuation ring $V$
(of characteristic zero with perfect residue field), there are compact objects
of $\DMt(\spec V,\QQ)$ which don't have any strong dual; see \ref{nonstrongduals}.
\end{paragr}

\begin{ex}\label{examplestrongdualsinD(K)}
Recall that $\Der(\KK)$ denotes the derived category of
$\KK$-vector spaces. This is a closed symmetric monoidal
triangulated category with tensor product $\otimes^{}_{\KK}$
and derived (internal) Hom $\derR\Hom_{\KK}$.
Note that for a complex of $\KK$-vector spaces $C$, the following conditions
are equivalent: 
\begin{itemize}
\item[(a)] $C$ is compact in $\Der(\KK)$;
\item[(b)] $C$ has a strong dual in $\Der(\KK)$;
\item[(c)] the $\KK$-vector space $\bigoplus_{n}H^n(C)$
is finite dimensional;
\item[(d)] $C$ is isomorphic in $\Der(\KK)$ to a bounded
complex of $\KK$-vector spaces which is degreewise finite dimensional.
\end{itemize}
\end{ex}

\begin{paragr}\label{defpoincaredualste1}
We consider again a stable theory $E$ and its associated commutative ring spectrum $\ste$ (\ref{defsteWeil}).

Let $X$ be a smooth and projective $S$-scheme of pure dimension $d$,
and denote by
$p:X \To S$ the canonical projection,
$\delta:X \To X \times_S X$ the diagonal embedding.

Then we can define pairings
\begin{align*}
\eta: \ & \ste \xrightarrow{p^*} \ste(X)(-d)[-2d]
 \xrightarrow{\delta_*} \ste(X \times_S X)(-d)[-2d]
  =\ste(X)(-d)[-2d] \otimes^{\derL}_\ste \ste(X) \\
\e: \ & \ste(X) \otimes^\derL_\ste \ste(X)(-d)[-2d]=\ste(X \times_S X)(-d)[-2d]
 \xrightarrow{\delta^*(-d)[-2d]} \ste(X) \xrightarrow{p_*} \ste.
\end{align*}
\end{paragr}

\begin{thm}[Poincar\'e duality]\label{dualitestemod}
The pair $(\e,\eta)$ defined above turns the object
$\ste(X)(-d)[-2d]$ into the strong dual of $\ste(X)$.
\end{thm}

\begin{proof}
This follows from the functoriality properties of the
Gysin morphism; see \cite[Theorem 5.23]{Deg6}.
\end{proof}

\begin{paragr}\label{strictnormalcrossingdual}
It can happen that $\ste(X)$ has a strong dual
for a non projective smooth $S$-scheme.
A classical example is the case where
$X=\overline{X}-D$, for a smooth and projective $S$-scheme $\overline{X}$
and a \emph{relative strict normal crossings divisor} $D$ in $\overline{X}$ (which means
here that $D$ is a divisor in $S$ with irreducible components $D_{i}$, $i\in I$,
such that $D$ is a reduced closed subscheme of $\overline{X}$, and such that
for any subset $J\subset I$, $D_{J}=\cap_{j\in J}D_{j}$
is smooth over $S$, and of codimension $\# J$ in $X$). The case where $D$ is irreducible comes
from Proposition \ref{purityste} and Theorem \ref{dualitestemod},
applied to $\overline{X}$ and $D$, and the general case follows
by an easy induction on the number of irreducible components of $D$.
As we already noticed, we cannot expect the object $\ste(X)$
to have a strong dual for an arbitrary $S$-scheme $X$ (\ref{nonstrongduals}).
However, when $S$ is the spectrum of a perfect field, $\ste(X)$
has a strong dual for any smooth $S$-scheme $X$; see \ref{compactdualdmQste}.
\end{paragr}

\subsection{Homological realization}





\begin{paragr}\label{defhomologicalrealization}
Let $E$ be a stable theory, and $\ste$ its associated
commutative ring spectrum. Recall $\Der(\KK)$ denotes
the (unbounded) derived category of the category of $\KK$-vector spaces.

We define the \emph{homological realization functor associated to $\ste$}
to be
\begin{equation}\label{defhomologicalrealization1}
\DMt(S,\ste)\To\Der(\KK)\quad , \qquad
M\longmapsto\derR\Hom_{\ste}(\ste,M)
\end{equation}
(where $\derR\Hom_{\ste}$ denotes the total right derived functor
of the $\Hom$ functor; see \eqref{constantsheafste2}).
This functor is right adjoint to the functor
\begin{equation}\label{defhomologicalrealization2}
\Der(\KK)\To\DMt(S,\ste)\quad , \qquad
C\longmapsto\ste\otimes^\derL_\KK\derL\Sigma^\infty(C_{S})
\end{equation}
(where $C_{S}$ denotes the constant sheaf associated to $C$).
As the functor \eqref{defhomologicalrealization2}
is obviously a symmetric monoidal functor, the
homological realization functor \eqref{defhomologicalrealization1}
is a lax symmetric monoidal functor. This means that for any
$\ste$-modules $M$ and $N$, there are coherent and natural maps
\begin{align}\label{defhomologicalrealization3}
&\derR\Hom_{\ste}(\ste,M)\otimes^{}_{\KK}\derR\Hom_{\ste}(\ste,N)
\To\derR\Hom_{\ste}(\ste,M\otimes^\derL_{\ste}N)\\
\intertext{and}
&\KK\To\derR\Hom_{\ste}(\ste,\ste)\, .
\end{align}
We define the category $\DMtp(S,\ste)$ to be the localizing subcategory~(cf. \ref{defA1equiv})
of the triangulated category $\DMt(S,\ste)$ generated by the objects which have a strong dual.

Note that any isomorphism $\ste(1)\simeq\ste$ (cf. \ref{orientomegaspectre})
induces an isomorphism $M(1)\simeq M$ in $\DMt(S,\ste)$ for any $\ste$-module $M$.
We deduce that the category $\DMtp(S,\ste)$ is stable by Tate twists.
Moreover, if $M$ and $N$ have strong duals, their tensor product $M\otimes^\derL_{\ste}N$
share the same property. In other words, $\DMtp(S,\ste)$ is
generated by a family of objects which is
stable by tensor product. This implies that the category $\DMtp(S,\ste)$ itself
is stable by tensor product in $\DMt(S,\ste)$. As a consequence, $\DMtp(S,\ste)$
is a symmetric monoidal category, and the inclusion functor from $\DMtp(S,\ste)$
into $\DMt(S,\ste)$ is symmetric monoidal. It is also obvious that any object
$M$ of $\DMtp(S,\ste)$ which has a strong dual $\dual M$ in $\DMt(S,\ste)$
has a strong dual in $\DMtp(S,\ste)$ which happens to be $\dual M$
itself.  There is a rather nice feature of the category $\DMtp(S,\ste)$:
an object of $\DMtp(S,\ste)$ is compact if and only if it has a strong dual.
The reason why this category $\DMtp(S,\ste)$ remains interesting is that,
by virtue of Poincar\'e duality (\ref{dualitestemod}),
for any smooth and projective $S$-scheme $X$,
the $\ste$-module $\ste(X)$ is in $\DMtp(S,\ste)$.

We finaly get a homological realization functor
\begin{equation}\label{defhomologicalrealization4}
\DMtp(S,\ste)\To\Der(\KK)\quad , \qquad
M\longmapsto\derR\Hom_{\ste}(\ste,M)
\end{equation}
by restriction of \ref{defhomologicalrealization1}.
\end{paragr}

\begin{thm}\label{homologicalrealization}
If $E$ is a mixed Weil theory, then
the homological realization functor
$$\DMtp(S,\ste)\To\Der(\KK)\quad , \qquad
M\longmapsto\derR\Hom_{\ste}(\ste,M)$$
is an equivalence of symmetric monoidal triangulated categories.
As a consequence, an object $M$ of $\DMtp(S,\ste)$ is compact
if and only if $\derR\Hom_{\ste}(\ste,M)$ is compact.
Moreover, for any $\ste$-module $M$ in $\DMtp(S,\ste)$,
there is a canonical isomorphism
$$\derR\Hom^{}_\ste(M,\ste)\simeq\derR\Hom^{}_\KK(\derR\Hom^{}_\ste(\ste,M),\KK).$$
In particular, if $M$ is compact, then we have canonical isomorphisms
$$\derR\Hom^{}_\ste(M,\ste)\simeq\derR\Hom^{}_\ste(\ste,\dual{M})
\simeq\dual{\derR\Hom^{}_\ste(\ste,M)}.$$
\end{thm}

\begin{proof}
The first step in the proof consists to see that the K\"unneth
formula implies that for any compact objects $M$ and $N$ of $\DMt(S,\ste)$,
the pairing
$$\derR\Hom_{\ste}(M,\ste)\otimes_{\KK}\derR\Hom_{\ste}(N,\ste)
\To\derR\Hom_{\ste}(M\otimes^\derL_{\ste}N,\ste)$$
is an isomorphism (it is sufficient to check this on a family
of compact generators, which is true by assumption for the family
that consists of the objects of type $\ste(X)$ for any smooth $S$-scheme $X$).


We will now prove that the homological realization functor \eqref{defhomologicalrealization4}
is a symmetric monoidal functor. The only thing to prove is in fact
that the map \eqref{defhomologicalrealization3} is
an isomorphism whenever $M$ and $N$ are in $\DMtp(S,\ste)$.
As $\DMtp(S,\ste)$ is generated by its compact objects, it is
sufficient to check this property when $M$ and $N$ are compact.
But in this case, $M$ and $N$ have strong duals, so that we get the
following isomorphisms.
$$\begin{aligned}
\derR\Hom_{\ste}(\ste,M)\otimes^{}_{\KK}\derR\Hom_{\ste}(\ste,N)
&\simeq\derR\Hom_{\ste}(\dual M,\ste)\otimes_{\KK}\derR\Hom_{\ste}(\dual N,\ste)\\
&\simeq\derR\Hom_{\ste}(\dual M\otimes^\derL_{\ste}\dual N,\ste)\\
&\simeq\derR\Hom_{\ste}(\dual{(M\otimes^\derL_{\ste}N)},\ste)\\
&\simeq\derR\Hom_{\ste}(\ste,M\otimes^\derL_{\ste}N)
\end{aligned}$$

We are now able to prove that the homological realization functor \eqref{defhomologicalrealization4}
is fully faithful.
Using the fact $\ste$ is compact in $\DMt(S,\ste)$ (see the end of \ref{defringspectrum2}),
we reduce to problem to showing fully faithfulness on compact objects.
Let $M$ and $N$ be compact objects of $\DMtp(S,\ste)$. We already noticed that
they both have strong duals $\dual M$ and $\dual N$ respectively. Note also
that symmetric monoidal functors preserve strong duals, so that
we get the following computations.
$$\begin{aligned}
\derR\Hom_{\ste}(M,N)&\simeq\derR\Hom_{\ste}(\ste,\dual M\otimes^\derL_{\ste}N)\\
&\simeq\dual{\derR\Hom_{\ste}(\ste,M)}\otimes_{\KK}\derR\Hom_{\ste}(\ste,N)\\
&\simeq\derR\Hom_{\KK}(\derR\Hom_{\ste}(\ste,M),\derR\Hom_{\ste}(\ste,N))
\end{aligned}$$
To prove the essential surjectivity, it is sufficient to check that a generating family
of $\Der(\KK)$ is in the essential image of the homological realization
functor. But this is obvious, as the object $\KK$ (seen as a complex concentrated in degree $0$)
generates $\Der(\KK)$.

The other assertions of the theorem are obvious consequences of this
equivalence of symmetric monoidal triangulated categories.
\end{proof}

\begin{paragr}\label{defhomstegm}
Assume that $E$ is a mixed Weil theory. Denote as usual by $\ste$
the commutative ring spectrum associated to $E$.

By virtue of Example \ref{examplestrongdualsinD(K)},
we know that $\derR\Gamma(X,\ste)$ is compact in $\Der(\KK)$
if and only if $\HH^*(X,\ste)$ is a finite dimensional
vector space. This is how Theorem \ref{homologicalrealization}
implies a finiteness result for $\HH^*(X,\ste)$ whenever $\ste(X)$
has a strong dual in $\DMt(S,\ste)$.

For any object $M$ of $\DMt(S,\ste)$ and any integers $p$ and $q$,
we get a pairing
\begin{equation}\label{genpairing1}
\Hom_{\DMt(S,\ste)}(\ste,M(-p)[-q])
\otimes^{}_{\KK}
\Hom_{\DMt(S,\ste)}(M,\ste(p)[q])
\To\KK
\end{equation}
inducing an isomorphism
\begin{equation}\label{genpairing2}
\Hom_{\DMt(S,\ste)}(M,\ste(p)[q])\simeq
\Hom_{\KK}(\Hom_{\DMt(S,\ste)}(\ste,M(-p)[-q]),\KK).
\end{equation}
whenever $M$ is in $\DMtp(S,\ste)$.

If $M$ has strong dual, the pairing \eqref{genpairing1} thus happens to be
a perfect pairing between finite dimensional $\KK$-vector spaces.

For a smooth $S$-scheme $X$, define the \emph{homology of $X$ with coefficients in
$\ste$} by the formula
\begin{equation}\label{defsmoothhomology}
\HH_{q}(X,\ste(p))=\Hom_{\DMt(S,\ste)}(\ste(p)[q],\ste(X))
\end{equation}
We get a canonical pairing
\begin{equation}\label{genpairing3}
\HH_{q}(X,\ste(p))\otimes^{}_{\KK}
\HH^q(X,\ste(p))\To\KK
\end{equation}
which happens to be perfect whenever $\ste(X)$
has a strong dual (e.g. when $X$ is projective).
For a smooth and projective $S$-scheme $X$ of pure dimension $d$,
Poincar\'e duality gives an isomorphism
\begin{equation}\label{genpairing4}
\HH^{2d-q}(X,\ste(d-p))\simeq
\HH_{q}(X,\ste(p))
\end{equation}
so that we get a perfect pairing
\begin{equation}\label{genpairing5}
\HH^{2d-q}(X,\ste(d-p))\otimes^{}_{\KK}
\HH^q(X,\ste(p))\To\KK
\ , \quad \alpha\otimes\beta\longmapsto\langle\alpha,\beta\rangle\, .
\end{equation}
Note that according to the definition of the duality pairing
\ref{defpoincaredualste1} and the K\"unneth formula, this pairing has
the familiar form~:
$$\langle\alpha,\beta\rangle=p_*(\alpha.\beta)$$
where $p_*:\HH^{2d}(X,\ste(d)) \rightarrow \HH^{0}(S,\ste(0))=\KK$
is the Gysin morphism associated to the canonical projection of $X/S$
 --- the so called \emph{trace morphism}.

For a smooth $S$-scheme $X$ of pure dimension $d$,
we can also define the \emph{cohomology with compact support
with coefficients in $\ste$} by the formula
\begin{equation}\label{defcohcompactsupport0}
\derR\Gamma_{c}(X,\ste(p))=\derR\Hom_{\ste}(\ste,\ste(X)(p-d)[-2d])\, .
\end{equation}
Setting $\HH^{q}_{c}(X,\ste(p))=\HH^q(\derR\Gamma_{c}(X,\ste(p)))$, we obtain
\begin{equation}\label{defcohcompactsupport}
\HH^{q}_{c}(X,\ste(p))=\Hom_{\DMt(S,\ste)}(\ste,\ste(X)(p-d)[q-2d])=\HH_{2d-q}(X,\ste(d-p))\, .
\end{equation}
It follows from Proposition \ref{gysinfctste1} that the cohomology with compact support
is functorial (in a contravariant way) with respect to projective morphisms,
and functorial (in a covariant way) with respect to equidimensional morphisms.
We also have a map
\begin{equation}\label{genpairing61}
\e: \ste(X) \otimes^\derL_\ste \ste(X)(-d)[-2d]
 \xrightarrow{\delta^*(-d)[-2d]} \ste(X) \xrightarrow{p_*} \ste
 \end{equation}
which defines by transposition a map
\begin{equation}\label{genpairing62}
\ste(X)(p-d)[-2d]\To\derR\sHom_{\ste}(\ste(X),\ste(p))\, .
\end{equation}
Note that $\derR\Hom_{\ste}(\ste,\derR\sHom_{\ste}(\ste(X),\ste(p)))=
\derR\Gamma(X,\ste(p))$. Hence we obtain a morphism
\begin{equation}\label{genpairing63}
\derR\Gamma_{c}(X,\ste(p))\To\derR\Gamma(X,\ste(p))
\end{equation}
which is functorial with respect to projective morphisms of $S$-schemes (thanks to the good
functorial properties of the Gysin morphisms),
and an isomorphism whenever $X$ is projective.
We also get a canonical pairing of complexes
\begin{equation}\label{genpairing7comp1}
\derR\Gamma_{c}(X,\ste(p))\otimes^{}_{\KK}\derR\Gamma(X,\ste(d-p)[2d])\To\KK
\end{equation}
defined by the canonical map
\begin{equation*}
\derR\Hom_{\ste}(\ste,\ste(X)(p-d)[-2d])\otimes_{\KK}\derR\Hom_{\ste}(\ste(X),\ste(d-p)[2d])
\To\derR\Hom_{\ste}(\ste,\ste) \, .
\end{equation*}
This gives rise to a pairing
\begin{equation}\label{genpairing7}
\HH^q_{c}(X,\ste(p))\otimes^{}_{\KK}
\HH^{2d-q}_{}(X,\ste(d-p))\To\KK
\end{equation}
which happens to be perfect if $\ste(X)$ has a strong dual in $\DMt(S,\ste)$.

Note that Poincar\'e duality gives rise to the following classical
computation; see e.g. \cite[3.3.3]{andremotifs}.
\end{paragr}

\begin{cor}[Lefschetz trace formula]\label{lefschetz}
Let $X$ and $Y$ be smooth and projective $S$-schemes
of pure dimension $d_{X}$ and $d_{Y}$ respectively,
Then, given two integers $p$ and $q$, for
$$\alpha\in\HH^{2d_{Y}+q}(X\times_{S}Y,\ste(d_{Y}+p))
\quad\text{and}\quad
\beta\in\HH^{2d_{X}-q}(Y\times_{S}X,\ste(d_{X}-p))\, ,$$
we have the equality
$$\langle\alpha,{}^{t}\beta\rangle=
\sum_{i}(-1)^{i}\mathrm{tr}(\beta\circ\alpha | \HH^{i}(X,\ste))\, ,$$
where ${}^{t}\beta\in\HH^{2d_{X}-q}(X\times_{S}Y,\ste(d_{X}-p))$
is the class corresponding to $\beta$ through the pullback by
the isomorphism $X\times_{S}Y\simeq Y\times_{S}X$,
$\langle \, . \, ,\, . \, \rangle$ is the Poincar\'e duality pairing,
and $\beta\circ\alpha$ denotes the composition of $\alpha$ and $\beta$
as cohomological correspondences.
\end{cor}

\begin{thm}\label{comparisonste}
Let $E$ and $E'$ be a mixed Weil theory and a stable theory
respectively. Denote by
$\ste$ and $\ste'$ the commutative ring spectra associated to $E$
and $E'$ respectively.
Let $u:E\To E'$ be a morphism of sheaves of differential graded $\KK$-algebras.
We assume that the induced map
$$\HH^1(\GG_{m},E)\To\HH^1(\GG_{m},E')$$
is not trivial. Then there exists a commutative ring spectrum $\ste''$
and two morphisms of ring spectra (which means morphisms of monoids in the
category of symmetric Tate spectra)
$$\ste\xrightarrow{ \ a \ }\ste''\xleftarrow{\ b \ }\ste'$$
with the following properties:
\begin{itemize}
\item[(a)] The map $\ste'\xrightarrow{ \ b \ }\ste''$ is an isomorphism in $\DMt(S,\KK)$.
\item[(b)] For any smooth $S$-scheme $X$, and any integer $n$,
the following diagram commutes (in which the vertical arrows are the canonical
isomorphisms).
$$\xymatrix{
\HH^n(X,E)\ar[rr]^{u}\ar[d]&&\HH^n(X,E')\ar[d]\\
\HH^n(X,\ste)\ar[r]^{a}&\HH^n(X,\ste'')&\HH^n(X,\ste')\ar[l]^{\simeq}_{b} }$$
\item[(c)] The maps $a$ and $b^{-1}$ define for any smooth
$S$-scheme $X$ maps
$$\HH^q(X,\ste(p))\To\HH^q(X,\ste'(p))\quad\text{and}\quad
\HH^q_{c}(X,\ste(p))\To\HH^q_{c}(X,\ste'(p))$$
which are compatible with cup products
and cycle class maps. If moreover $\ste(X)$
has a strong dual (e.g., from \ref{strictnormalcrossingdual},
if $X$ is the complement
of a relative strict normal crossings divisor in a smooth and
projective $S$-scheme), then these maps are bijective.
\end{itemize}
In particular, if moreover for any smooth $S$-scheme $X$, the $\ste$-module $\ste(X)$
has a strong dual in $\DMt(S,\ste)$, $E'$ is a mixed Weil theory
and the map $u$ is a quasi-isomorphism of complexes of Nisnevich sheaves.
\end{thm}

\begin{proof}
We have to come back to the very construction
of the ring spectrum associated to a stable theory
given in \ref{defsteWeil}. Define
$$L=\Hom^*(\KK(1),E)_{S}\quad\text{and}\quad L'=\Hom^*(\KK(1),E')_{S}.$$
We know that the symmetric Tate spectra $\ste$ and $\ste'$
are defined respectively by the sheaves of complexes
$$\ste_{n}=\sHom(L^{\otimes n},E)\quad\text{and}\quad
\ste'_{n}=\sHom(L'^{\otimes n},E').$$
Define a third ring spectrum $\ste''=(\ste''_{n},\tau_{n})$ as follows.
Put $\ste''_{n}=\sHom(L^{\otimes n},E')$. We have maps
$$L\To L'\To\sHom(\KK(1),E')$$
from which we construct maps of type
$$\tau'_{n}:\KK(1)\otimes_{\KK}L'\otimes_{\KK}\sHom(L'^{\otimes n},E)
\To \sHom(L'^{\otimes n},E)$$
following the same steps as for the construction of the map \eqref{defsigmansteintermediaire3}.
The structural maps
$$\tau_{n}:\ste''_{n}(1)\To\ste''_{n+1}$$
are defined by transposition of the maps $\tau'_{n}$.
One can then check that $\ste''$ is a commutative ring spectrum.
The map $a$ is induced by the maps
$$a_{n}:\sHom(L^{\otimes n},E)\To\sHom(L^{\otimes n},E')$$
which correspond to the composition with $u$,
and the map $b$ is induced by the maps
$$b_{n}:\sHom(L'^{\otimes n},E')\To\sHom(L^{\otimes n},E')$$
which corresponds to the composition with the map $L\To L'$
obtained from $u$ by functoriality.
These define the expected morphisms of ring spectra.

Property (a) comes obviously from the fact the map $L\To L'$
has to be a quasi-isomorphism according to the assumption on $u$.
Indeed, this implies the maps $b_{n}$ are all quasi-isomorphisms as well.
In particular, the total left derived functor of the base change
functor induced by $b$ is an equivalence of triangulated categories
$$\DMt(S,\ste')\simeq\DMt(S,\ste'')\, .$$

As a consequence,
the total left derived functor of the base change functor induced by $a$
$$\DMt(S,\ste)\To\DMt(S,\ste')\quad , \qquad M\longmapsto \ste''\otimes^\derL_{\ste}M$$
is a triangulated functor which preserves small direct sums, and it is also
symmetric monoidal. 
We claim that this induces by restriction a fully faithful symmetric monoidal
triangulated functor
$$\DMtp(S,\ste)\To\DMt(S,\ste')\, .$$
To see this, we first remark that $\ste$ is a compact
generator of $\DMtp(S,\ste)$: applying Theorem \ref{homologicalrealization}
to $E$ implies that $\DMtp(S,\ste)$ is equivalent to $\Der(\KK)$.
As the base change functor sends $\ste$ to $\ste''\simeq\ste'$, it is sufficient to prove
that the induced maps
$$\Hom_{\DMt(S,\ste)}(\ste,\ste[n])\To\Hom_{\DMt(S,\ste')}(\ste',\ste'[n])$$
are bijective. For $n\neq 0$, the two terms are null, and for
$n=0$, this map is a morphism of $\KK$-algebras from $\KK$ to itself,
so that it has to be an identity.

Properties (b) and (c) follow immediately
from this fully faithfulness (the compatibility with cycle class maps
follows from Theorem \ref{higherChernste}).
\end{proof}

\subsection{Cohomology of motives}

\begin{paragr}
In this section, the base scheme $S$ is the spectrum of a perfect
field $k$.

We consider given a stable cohomology theory $E$,
as well as its associated ring spectrum $\ste$.
Let $\mathit{TD}_{\AA^1}(k,\ste)$ the localizing subcategory of
$\DMt(k ,\ste)$ generated by objects of type $\ste(p)[q]$,
$p,q\in\ZZ$.
\end{paragr}

\begin{prop}\label{Tatehomologicalrealization}
The functor
$$\mathit{TD}_{\AA^1}(k,\ste)\To\mathit{D}(\KK)\ ,
\quad M\longmapsto\derR\Hom_{\ste}(\ste,M)$$
is an equivalence of symmetric monoidal triangulated categories.
\end{prop}

\begin{proof}
This functor is a right adjoint to the symmetric monoidal
functor
$$\mathit{D}(\KK)\To\mathit{TD}_{\AA^1}(k,\ste)$$
which sends a complex $C$ to $\ste\otimes^\derL_\KK \Sigma^\infty C$.
It is sufficient to prove that the latter is an equivalence
of categories. This follows essentially from the Homotopy axiom W1:
this implies that this functor is fully faithful on the
set of compact generators given by the unit object of $\mathit{D}(\KK)$,
which is sent to $\ste$. As $\ste(p)\simeq\ste$ for any integer $p$,
and as $\ste$ is compact in $\mathit{TD}_{\AA^1}(k,\ste)$,
we get the essential surjectivity by definition of
$\mathit{TD}_{\AA^1}(k,\ste)$.
\end{proof}

\begin{cor}\label{Tatehomologicalrealization2}
For any object $M$ of $\mathit{TD}_{\AA^1}(k,\ste)$, we have
a canonical isomorphism
$$\derR\Hom^{}_\ste(M,\ste)\simeq\derR\Hom^{}_\KK(\derR\Hom^{}_\ste(\ste,M),\KK).$$
\end{cor}

\begin{proof}
This follows from a straightforward translation from the
equivalence of categories given by Proposition \ref{Tatehomologicalrealization}.
\end{proof}

\begin{prop}\label{Tatehomologicalrealization3}
The $\ste$-module $\ste\otimes^\derL_\QQ\HQ$
is in $\mathit{TD}_{\AA^1}(k,\ste)$.
\end{prop}

\begin{proof}
We know that $\HQ\simeq\motcoh$ is a direct factor of
the $K$-theory spectrum $\BGLQ$. Hence it is sufficient to
prove that $\ste\otimes^\derL_\QQ\BGLQ$ is in
$\mathit{TD}_{\AA^1}(k,\ste)$, which follows immediately
from \cite[Theorem 6.2]{dugisakcell}.
\end{proof}

\begin{paragr}
Remember from \ref{higherChernste2}
we have an isomorphism $\motcoh\simeq\HQ$ in the category $\DMt(k ,\QQ)$.
Let $\mathit{cl}:\HQ\To\ste$ be the cycle class map \eqref{VoecycleclassmapHQste}.
It induces by adjunction a $\ste$-linear map
$$\ste\otimes^\derL_\QQ\HQ\To\ste\, .$$
\end{paragr}

\begin{prop}\label{isostemotcohste}
The map $\ste\otimes^\derL_\QQ\HQ\To\ste$
is an isomorphism in the category $\DMt(k ,\ste)$.
\end{prop}

\begin{proof}
We know from Theorem \ref{higherChernste}
that there is a canonical isomorphism in $\mathit{D}(\KK)$:
$$\derR\Hom_\ste(\ste\otimes^\derL_\QQ\HQ,\ste)=\derR\Hom_\QQ(\HQ,\ste)\simeq\KK$$
(where $\KK$ is seen as a complex concentrated in degree $0$).
By virtue of Proposition \ref{Tatehomologicalrealization3},
we can apply Corollary \ref{Tatehomologicalrealization2}
to $\ste\otimes^\derL_\QQ\HQ$ to obtain an isomorphism
$$\derR\Hom_\ste(\ste\otimes^\derL_\QQ\HQ,\ste)\simeq
\derR\Hom_\KK(\derR\Hom(\ste,\ste\otimes^\derL_\QQ\HQ),\KK)\, ,$$
This implies that we have an isomorphism in $\mathit{D}(\KK)$:
$$\KK\simeq\derR\Hom(\ste,\ste\otimes^\derL_\QQ\HQ)\, .$$
As $\derR\Hom_\ste(\ste,\ste)\simeq\KK$, and as, by Proposition
\ref{Tatehomologicalrealization}, the
homological realization functor $\derR\Hom_\ste(\ste,-)$
is an equivalence of categories from $\mathit{TD}_{\AA^1}(k,\ste)$
to $\mathit{D}(\KK)$, to prove that the map
$\ste\otimes^\derL_\QQ\HQ\To\ste$ is an isomorphism, we are reduced to
check that it is not trivial, which is obvious, by definition
of the cycle class map $\mathit{cl}$.
\end{proof}

\begin{paragr}\label{steHqmodule}
The canonical map from $\ste$
to $\ste\otimes^\derL_\QQ\HQ$ is an inverse of the isomorphism
of Proposition \ref{isostemotcohste}, whence it is
an isomorphism as well. We deduce from this the following result.
\end{paragr}

\begin{prop}\label{HQtostemonoidal}
The functor
$$\DMt(k,\HQ)\To\DMt(k,\ste)\ , \
M\longmapsto \ste\otimes^\derL_\QQ M$$
is a symmetric monoidal triangulated functor.
\end{prop}

\begin{proof}
As we are working with rational coefficients, using
\cite[proposition 4.3.21]{DAG3} (see also \cite{BerMoe,hinich}),
we can see that there is a commutative monoid structure\footnote{Our purpose
is to deal with symmetric monoidal structures on homotopy categories of
modules over a commutative monoid. A natural setting for this is the
notion of $E_\infty$-algebra. But, as we are working with
rational coefficients, it is possible to strictify any $E_\infty$-algebra
into a commutative monoid, so that we have chosen to remain coherent with the
rest of these notes, by considering genuine commutative monoids.
One could also avoid any complication by working directly with
symmetric monoidal $\infty$-categories~\cite{DAG3}.}
on the derived tensor product $\ste\otimes^\derL_\QQ\HQ$.
Proposition \ref{isostemotcohste} tells us that
the canonical map $\ste\To\ste\otimes^\derL_\QQ\HQ$
is an isomorphism in the homotopy category of
of commutative ring spectra (defined by stable $\AA^1$-equivalences).
Notice that, by virtue of \cite[Proposition 6.35]{HCD}, we can apply
\cite[Theorem 4.3]{SS} to see that $\ste\To\ste\otimes^\derL_\QQ\HQ$
induces an equivalence of symmetric monoidal triangulated categories
$$\DMt(k,\ste)\simeq\DMt(k,\ste\otimes^\derL_\QQ\HQ)\, .$$
The base change functor along $\HQ\To\ste\otimes^\derL_\QQ\HQ$ thus
gives a symmetric monoidal triangulated functor
$$\DMt(k,\HQ)\To
\DMt(k,\ste\otimes^\derL_\QQ\HQ)\simeq\DMt(k,\ste)\, .$$
The formula
$$\ste\otimes^\derL_\QQ M\simeq\ste\otimes^\derL_\QQ\HQ\otimes^\derL_\HQ M$$
shows that the functor we constructed above is (isomorphic to)
the functor considered in the proposition.
\end{proof}

\begin{paragr}
Let $\DM(k)$ be the triangulated category of mixed motives
over $k$; see \cite[Example 7.15]{HCD} for its construction.
This is a symmetric monoidal
triangulated category (as the homotopy category of a stable symmetric
monoidal model category), and it is generated, as a triangulated category,
by its compact objects. Moreover, the full subcategory of compact
objects in $\DM(k)$ is canonically equivalent to
Voevodsky's triangulated category of mixed motives $\DM_\mathit{gm}(k)$,
constructed in \cite{FSV5}.
Different (but equivalent) constructions of $\DM(k)$ are given
by \cite[Theorem 35]{MZ}, and the relation with $\DM_\mathit{gm}(k)$
is described in \cite[Section 2.3]{MZ}; a systematic study of the
triangulated categories $\DM(S)$ will appear in \cite{TCMM}.
We will denote by $\DM(k,\QQ)$ the rational version of $\DM(k)$,
and by $\DM_\mathit{gm}(k,\QQ)$ the rational version of
$\DM_\mathit{gm}(k)$. By virtue of
\cite[Theorem 68]{MZ}, there is a canonical equivalence of
symmetric monoidal triangulated categories\footnote{The equivalence
of categories \eqref{eqHQmodDM} is proved in \cite{MZ} using
resolution of singularities by de Jong alterations \cite{dejong};
however, it will be shown in \cite{TCMM} that such an equivalence
of triangulated categories holds over a geometrically unibranch base scheme,
by very different methods (without any kind of resolution of singularities).}
\begin{equation}\label{eqHQmodDM}
\mathcal U:\DM(k,\QQ)\overset{\sim}{\To}\DMt(k,\HQ)
\end{equation}
which sends the motive of $X$ twisted by $p$
to the object $\HQ\otimes^\derL_\QQ\Sigma^\infty(\QQ(X))(p)$
(for $X/k$ smooth, and $p\in\ZZ$);
it is induced by the forgetful functor from the category of Nisnevich sheaves
with transfers to the category of Nisnevich sheaves on $\sm/k$.
\end{paragr}

\begin{thm}\label{strongdualinDMQ}
The motives of shape $M_\mathit{gm}(X)(p)$,
for $X$ smooth and projective, and $p\in\ZZ$, form
a set of compact generators in $\DM(k,\QQ)$.
In particular, an object of $\DM(k,\QQ)$ is compact if and only if it
has a strong dual.
\end{thm}

\begin{proof}
This is proven using de Jong's resolution
of singularities by alterations \cite{dejong};
see the proof of \cite[Theorem 68]{MZ}.
\end{proof}

\begin{cor}\label{compactdualdmQste}
The following equality holds.
$$\DMtp(k,\ste)=\DMt(k,\ste) \, .$$
If moreover $E$ is a mixed Weil theory, then
the homological realization functor \eqref{defhomologicalrealization1}
defines an equivalence of symmetric monoidal triangulated
categories
$$\DMt(k,\ste)\simeq\Der(\KK)\, .$$
In particular, for any smooth $k$-scheme $X$, $\ste(X)$ has a strong dual, so
that \eqref{genpairing7} is a perfect pairing between finite dimensional
vector spaces.
\end{cor}

\begin{proof}
The first assertion follows immediately from Theorem \ref{strongdualinDMQ}.
Theorem \ref{homologicalrealization} then ends the proof.
\end{proof}

\begin{cor}\label{compkperfectste}
Assume that $E$ is a mixed Weil theory.
For any $\KK$-linear stable theory $E'$ defined on smooth $k$-schemes,
a morphism of sheaves of differential graded $\KK$-algebras $E\To E'$
is a quasi-isomorphism (in the category of complexes of Nisnevich sheaves)
if and only if the induced map $\HH^1(\GG_{m},E)\To\HH^1(\GG_{m},E')$
is not trivial.
\end{cor}

\begin{proof}
Apply Theorem \ref{comparisonste} and Corollary \ref{compactdualdmQste}.
\end{proof}

\begin{cor}\label{stablevsmixedWeil}
Assume that, for any smooth and projective $k$-schemes $X$
and $Y$, the K\"unneth map
$$
\bigoplus_{p+q=n}H^p(X,E) \otimes_\KK H^q(Y,E)
 \overset{\sim}{\To} H^{n}(X\times_{k}Y,E) \ .
$$
is an isomorphism.

Then $E$ is a mixed Weil theory.
\end{cor}

\begin{proof}
We claim that for any compact objects $M$ and $N$ of $\DMt(k,\ste)$,
the map
$$\derR\Hom_{\ste}(M,\ste)\otimes_{\KK}\derR\Hom_{\ste}(N,\ste)
\To\derR\Hom_{\ste}(M\otimes^\derL_{\ste}N,\ste)$$
is an isomorphism: it is sufficient to check this on a set of compact generators,
which is true by assumption, by virtue of Theorem \ref{strongdualinDMQ}.
\end{proof}

\begin{thm}\label{sterealDMk}
Let $E$ be a mixed Weil theory on smooth $k$-schemes,
and $\ste$ its associated commutative ring spectrum.
Then the motivic homological realization functor
$$\DM(k,\QQ)\To\Der(\KK)\quad , \qquad
M\longmapsto\derR\Hom_\QQ(\QQ,\ste\otimes^\derL_{\QQ}\mathcal U(M))$$
is a symmetric monoidal triangulated functor which preserves
compact objects. In particular, if $\Der^b(\KK)$ denotes the
bounded derived category of the category of finite dimensional
$\KK$-vector spaces, it induces by restriction a symmetric
monoidal triangulated functor
$$R_{\ste} : \DM_{\mathit{gm}}(k,\QQ)\To\Der^b(\KK)$$
such that, for any smooth $k$-scheme $X$, one has canonical
isomorphisms
$$R_{\ste}(\dual{M_{\mathit{gm}}(X)})\simeq
\dual{R_{\ste}(M_{\mathit{gm}}(X))}\simeq\derR\Gamma(X,\ste)\, .$$
\end{thm}

\begin{proof}
Under the equivalence of categories \eqref{eqHQmodDM} this
functor corresponds to the composition of the functor
of Proposition \ref{HQtostemonoidal} with the homological
realization functor \eqref{defhomologicalrealization1}.
Hence the first assertion follows from Corollary \ref{compactdualdmQste}.
In particular, this functor
preserves strong duals. Theorem \ref{strongdualinDMQ} now
implies it sends $ \DM_{\mathit{gm}}(k,\QQ)$ to $\Der^b(\KK)$.
If $X$ is a smooth $k$-scheme, we have a natural isomorphism
$$\HQ\otimes^\derL_\QQ\Sigma^\infty\QQ(X)\simeq\mathcal U(M_\mathit{gm}(X))\, .$$
We deduce from Proposition \ref{isostemotcohste} that
$$\ste(X)\simeq\ste\otimes^\derL_\QQ\HQ\otimes^\derL_\QQ\Sigma^\infty\QQ(X)\simeq
\ste\otimes^\derL_\QQ\mathcal U(M_\mathit{gm}(X))\, ,$$
which implies that
$$\begin{aligned}
\derR\Hom_\QQ(\QQ,\ste\otimes^\derL_{\QQ}\mathcal U(M_\mathit{gm}(X)))
&\simeq \derR\Hom_\QQ(\QQ,\ste(X))\\
&\simeq \derR\Hom_\ste(\ste,\ste(X))\, .
\end{aligned}$$
By Theorem \ref{homologicalrealization}, we get isomorphisms
$$R_{\ste}(\dual{M_{\mathit{gm}}(X)})\simeq
\dual{R_{\ste}(M_{\mathit{gm}}(X))}\simeq\derR\Gamma(X,\ste)\, ,$$
which ends the proof.
\end{proof}

\begin{rem}\label{conclusion}
The functor $R_\ste$ induces cycle class maps
$$\HH^q(X,\QQ(p))\To\HH^q(X,\ste(p))=\HH^q(X,E)(p)$$
which coincide with the cycle class maps introduced
in \ref{higherChernste2}. These cycle class maps are compatible with
first Chern classes, hence with Gysin maps
(by the categorical construction of these; see \cite{Deg6}).

The reader might have noticed that, in the definition
of a mixed Weil cohomology, we didn't ask
the differential graded algebra $E$ to be concentrated in
non negative degrees. It would be natural to ask the
cohomology groups $H^n(X,E)$ to vanish for any (affine) smooth scheme
$X$ and any negative integer $n$ (which is true in practice).
We conjecture this vanishing property to hold in general.

The existence of cycle class maps compatible with cup products
and with Gysin morphisms finally proves that the cohomology groups $H^n(X,E)$,
for $X$ smooth and projective over $k$, define a Weil
cohomology in the sense of \cite[Definition 3.3.1.1]{andremotifs},
modulo the vanishing property discussed above.
\end{rem}

\begin{paragr}
The proof of Theorem \ref{sterealDMk} relies on
Proposition \ref{HQtostemonoidal} and on the description
of $\DM(k,\QQ)$ as the homotopy category of
modules on the rational motivic cohomology spectrum.
Another strategy to prove Theorem \ref{sterealDMk}
is to identify $\DM(k,\QQ)$ with the ``orientable part''
of $\DMt(S,\QQ)$. This is achieved using an unpublished
result of F.~Morel~\cite{ratmot}, which computes the rational
motivic sphere spectrum in terms of motivic
cohomology spectrum; see Theorem \ref{ratmotmorel}.
More precisely, another proof of Theorem \ref{sterealDMk} is given by
Corollary \ref{realmotivesste}, equality \eqref{motivesperfectfield1},
and Theorem \ref{eqbeilinsonvoevodsky} below. Moreover,
Morel's result gives a very straightforward proof
of the existence and unicity of the cycle class map
for a stable theory; see Remark \ref{remratmotmorel}.
We will now outline this alternative point of view.
\end{paragr}

\begin{paragr}\label{morel}
Let $S$ be a scheme.
The permutation isomorphism
\begin{equation}\label{morel1}
\tau : \tQQ(1)[1]\otimes^\derL_{\QQ}\tQQ(1)[1]\To\tQQ(1)[1]\otimes^\derL_{\QQ}\tQQ(1)[1]
\end{equation}
satisfies the equation $\tau^2=1$ in $\DMt(S,\QQ)$.
Hence it defines an element $\epsilon$ in $\mathrm{End}_{\DMt(S,\QQ)}(\tQQ)$
which also satisfies the relation $\epsilon^2=1$.
We define two projectors
\begin{equation}\label{morel2}
e_{+}=\frac{1-\epsilon}{2}
\quad\text{and}\quad
e_{-}=\frac{1+\epsilon}{2}\, .
\end{equation}
As the triangulated category $\DMt(S,\QQ)$ is pseudo abelian,
we can define two objects by the formul{\ae}:
\begin{equation}\label{morel3}
\QQ_{+}=\mathrm{Im}\, e_{+}
\quad\text{and}\quad
\QQ_{-}=\mathrm{Im}\, e_{-}\, .
\end{equation}
Then for an object $M$ of $\DMt(S,\QQ)$, we set
\begin{equation}\label{morel4}
M_{+}=\QQ_{+}\otimes^\derL_{\QQ}M
\quad\text{and}\quad
M_{-}=\QQ_{-}\otimes^\derL_{\QQ}M\, .
\end{equation}
It is obvious that for any objects $M$ and $N$ of $\DMt(S,\QQ)$,
one has
\begin{equation}\label{morel5}
\Hom_{\DMt(S,\QQ)}(M_{i},N_{j})=0\quad\text{for $i,j\in\{+,-\}$ with $i\neq j$.}
\end{equation}
Denote by $\DMt(S,\QQ_{+})$ (resp. $\DMt(S,\QQ_{-})$) the full subcategory
of $\DMt(S,\QQ)$ made of objects which are isomorphic to some $M_{+}$
(resp. some $M_{-}$) for an object $M$ in $\DMt(S,\QQ)$. Then
\eqref{morel5} implies that the direct sum functor induces an
equivalence of triangulated categories
\begin{equation}\label{morel6}
\DMt(S,\QQ_{+})\times\DMt(S,\QQ_{-})\simeq\DMt(S,\QQ)\, .
\end{equation}
Assume now that $S$ is a regular scheme. Recall from \ref{recollection} the Beilinson motivic cohomology spectrum $\motcoh$. A deep result announced
by F.~Morel in \cite{ratmot} takes the following form (taking into account the equivalence of categories \eqref{defkthrat}).
\end{paragr}

\begin{thm}\label{ratmotmorel}
We have a canonical identification
$\QQ_{+}=\motcoh$. Moreover, if $-1$ is a sum of squares in $\bundle O(S)$,
then $\tQQ=\motcoh$.
\end{thm}

A proof will be given in \cite{TCMM}.

\begin{paragr}\label{defmorelbeilinsonmotives}
For a general scheme $S$, we define the \emph{triangulated
category of Morel-Beilinson motives} to be
\begin{equation}\label{morel7}
\DMB(S)=\DMt(S,\QQ_{+})\, .
\end{equation}
Note that according to \cite{ayoub}, the Grothendieck six operations are defined on the categories $\DMt(S,\QQ)$. As all these operations commute
with Tate twists, it is obvious that they preserve Morel-Beilinson motives. Hence the categories $\DMB(S)$ for various schemes $S$ are stable by the
six operations as subcategories of $\DMt(S,\QQ)$. In particular, $\DMB(S)$ is a symmetric monoidal triangulated category, and the canonical functor
from $\DMt(S,\QQ)$ to $\DMB(S)$ is a symmetric monoidal triangulated functor.
\end{paragr}

\begin{paragr}\label{steoriented0}
Suppose now that $S$ is a regular scheme. Consider
a stable theory $E$ defined on smooth $S$-schemes,
and let $\ste$ be its associated commutative ring spectrum.
\end{paragr}

\begin{prop}\label{steoriented}
We have $\ste=\ste_{+}$.
\end{prop}

\begin{proof}
This is a translation from Lemma \ref{lm:premutation_induce_iso_proj}.
\end{proof}

\begin{rem}\label{remratmotmorel}
Theorem \ref{ratmotmorel} and Proposition \ref{steoriented}
give another proof of Theorem \ref{higherChernste}:
the unit map $\tQQ=\QQ_+\oplus\,\QQ_-\To\ste$ factors uniquely
through $\tQQ_+=\motcoh$, which gives the cycle class map
$\motcoh\To\ste$ (it clearly preserves the unit, so that it has to be
the map obtained from the Chern character by Theorem \ref{ratmotmorel}).
This construction has the advantage
of giving directly the compatibilities of the
cycle class map with the algebra structures.
\end{rem}

\begin{paragr}\label{dualityDMbeil}
Define $\DMB^\vee(S)$ as the localizing subcategory~(\ref{defA1equiv})
of $\DMB(S)$ generated by the objects which have a strong dual
(e.g. $\tQQ(X)_{+}(p)$ for a smooth and projective $S$-scheme $X$ and
an integer $p$; see \cite{ayoub,roe}).
\end{paragr}

\begin{cor}\label{realmotivesste}
If $E$ is a mixed Weil theory, then
the \emph{motivic homological realization functor}
$$\DMB^\vee(S)\To\Der(\KK)\quad , \qquad
M\longmapsto\derR\Hom_{\QQ}(\motcoh,\ste\otimes^\derL_{\QQ}M)$$
is a symmetric monoidal triangulated functor.
\end{cor}

\begin{proof}
By virtue of Theorem \ref{ratmotmorel} and of the preceding proposition,
this functor is isomorphic to
the composition of the symmetric monoidal triangulated functor
$$\DMB^\vee(S)\To\DMtp(S,\ste)\quad , \qquad
M\longmapsto\ste\otimes^\derL_{\QQ}M$$
with the homological realization functor \eqref{defhomologicalrealization4}.
Theorem \ref{homologicalrealization} concludes.
\end{proof}

\begin{paragr}\label{motivesperfectfield}
Assume now $S$ is the spectrum of a perfect field $k$.

It follows then from \cite{riou} that we have
\begin{equation}\label{motivesperfectfield1}
\DMB^\vee(\spec k)=\DMB(\spec k)\, .
\end{equation}
\end{paragr}

\begin{thm}[F. Morel]\label{eqbeilinsonvoevodsky}
There exists a canonical equivalence of symmetric monoidal triangulated
categories
$$\DMB(\spec k)\simeq\DM(k,\QQ)\, .$$
\end{thm}

\begin{proof}
We know (e.g. from \cite{MZ}) that we have a canonical
symmetric monoidal triangulated functor
\begin{equation}\label{eqbeilinsonvoevodsky1}
\DMt(\spec k,\QQ)\To\DM(k,\QQ)\quad , \qquad M\longmapsto M_{\mathit{tr}}
\end{equation}
which preserves Tate twists, direct sums, and compact objects.
By virtue of \cite[Corollary 2.1.5]{FSV5},
the functor \eqref{eqbeilinsonvoevodsky1} vanishes on $\DMt(\spec k,\QQ_{-})$,
so that it induces a symmetric monoidal triangulated functor
\begin{equation}\label{eqbeilinsonvoevodsky2}
\DMB(\spec k)\To\DM(k,\QQ)\quad , \qquad M\longmapsto M_{\mathit{tr}}\, .
\end{equation}
It then follows from Theorem \ref{ratmotmorel}
and \cite[Theorem \textsc{v}.31]{thriou} that
for a given smooth $k$-scheme $X$ and two integers $p$ and $q$,
the induced map
$$\Hom_{\DMB(\spec k)}(\tQQ(X)_{+},\motcoh(p)[q])\To
\Hom_{\DM(k,\QQ)}(\tQQ(X)_{\mathit{tr}},\tQQ_{\mathit{tr}}(p)[q])$$
is in fact the isomorphism \eqref{motcohbeilinson9}.
By \ref{motivesperfectfield1},
this implies that  the functor \eqref{eqbeilinsonvoevodsky2}
is fully faithful on compact objects which
proves the full faithfulness.
The essential surjectivity follows from the fact that,
by the very construction of $\DM(k,\QQ)$, the objects of shape
$\tQQ(X)_{\mathit{tr}}(p)[q]$ generate $\DM(k,\QQ)$.
\end{proof}

\begin{rem}
It will be proved in \cite{TCMM} that Morel's Theorem \ref{ratmotmorel}
implies that Theorem \ref{eqbeilinsonvoevodsky} is true over any geometrically
unibranch base scheme.
\end{rem}



\section{Some classical mixed Weil cohomologies}

\subsection{Algebraic and analytic de Rham cohomologies}

\begin{paragr}
Suppose $k$ is a field of characteristic $0$.
Let $X$ be a smooth $k$-scheme. We denote by $\Omega^1_{X/k}$
the locally free sheaf of algebraic differential forms on $X$ over $k$.
Then the de Rham complex is the complex of $\bundle O_{X}$-modules
obtained from the exterior $\bundle O_{X}$-algebra generated by $\Omega^1_{X/k}$:
$$\Omega^*_{X/k}=\bigwedge\Omega^1_{X/k}.$$
Remember from \cite{cohDR,crys}
that the algebraic de Rham cohomology of $X$ is defined to be
$$
\HH^*_{\dR}(X)=\hypercoh^*_{\zar}(X,\Omega^*_{X/k}).
$$
We will show here that de Rham cohomology is canonically
represented by a mixed Weil theory.
\end{paragr}

\begin{paragr}
Let $X/k$ be an affine smooth scheme. We simply put
$\Omega_{\dR}(X)=\Gamma(X,\Omega^*_{X/k})$.
Then $\Omega_{\dR}(X)$ is a commutative graded differential
algebra and it defines a presheaf of commutative differential graded $k$-algebras
$$
\Omega_{\dR}:X\longmapsto\Omega_{\dR}(X).
$$
In this context, the K\"unneth formula is obvious: the canonical map
$$
\Omega_{\dR}(X) \otimes_{k} \Omega_{\dR}(Y)
 \To \Omega_{\dR}(X \times_k Y)
$$
is an isomorphism.

As $\Omega^*_{X/k}$ is a complex of coherent sheaves
on $X$ and $X$ is affine, the vanishing theorem of Serre \cite[1.3.1]{EGA3}
and the spectral sequence
$$E^1_{p,q}=\HH^p_{\zar}(X,\Omega^q_{X/k})\Rightarrow\HH^{p+q}_{\dR}(X)$$
implies
$$
\HH^*_{\dR}(X)=\HH^*(\Omega_{\dR}(X)).
$$
\end{paragr}

\begin{paragr}
The complex $\Omega_{\dR}$ satisfies \'etale descent on smooth
$k$-schemes, thus Nisnevich
descent. This means the following.

Let $X=\spec A$ and $Y=\spec B$ be smooth affine schemes
and $f:Y \To X$ an \'etale morphism.
Then, $\Omega_{\dR}(Y)=\Omega_{\dR}(X) \otimes_A B$.
Suppose $f$ is an \'etale cover.
The augmented \v Cech complex
$\check C^+_*(Y/X)$
is associated to the differential graded $A$-algebra
$$
T^+_A(B)=(A \To B \To B \otimes_A B
 \To B \otimes_A B \otimes_A B
 \To\dots)
$$

Thus,
$\Omega_{\dR}(\check C^+_*(Y/X))=\Omega_{\dR}(Y) \otimes_A T^+_A(B)$.

As $f$ is faithfully flat, it is a morphism of effective descent
with respect to the fibred category of quasi-coherent modules
(see \cite[Expos\'e VIII, Theorem 1.1]{SGA1}),
so that the complex $T^+_A(B)$ is acyclic.
For any integer $r \geq 0$, $\Omega_{\dR}^r(Y)$ is flat
over $A$, thus $\Omega_{\dR}^r(Y) \otimes_A T^+_A(B)$ is acyclic.
Hence the spectral sequence of a bounded bicomplex shows
the complex
$\tot[\Omega_{\dR}(\check C^+_*(Y/X))]$ is acyclic.
This implies the \'etale descent for algebraic de Rham cohomology; see \cite{artin}.
We deduce easily from the computations above that
for any distinguished square as \eqref{distsquare}
which consists of smooth affine $k$-schemes, we get
a short exact sequence
$$0\To\Omega_{\dR}(X)\To \Omega_{\dR}(U)\oplus \Omega_{\dR}(V)
\To\Omega_{\dR}(U\times_{X}V)\To 0.$$
Hence $\Omega_{\dR}$ has the B.-G.-property with respect
to the Nisnevich topology on the category of affine smooth $k$-schemes.
\end{paragr}

\begin{paragr}\label{dRW1andW2}
Finally, the following computations are easy:
$$
\begin{aligned}
H^n_{\dR}(\AA^1_{k})
&=
\begin{cases}
k & \text{if } n=0 \\ 0 & \text{otherwise}
\end{cases}\\
H^n_{\dR}(\GG_m) &=
\begin{cases}
k & \text{if } n=0 \\
k.\dlog& \text{if } n=1 \\
0 & \text{otherwise}
\end{cases}
\end{aligned}
$$
where $\dlog$ is the differential form defined by $\dlog(t)=dt/t$.
In conclusion, we have proved:
\end{paragr}

\begin{prop}\label{dRmixedWeil}
The presheaf $\Omega_{\dR}$ is a mixed Weil theory.
\end{prop}

\begin{paragr}
We denote by $\ste_{\dR}$ the
corresponding commutative ring spectrum.
Recall the canonical map
$$H^*_{\dR}(X)\To\hypercoh^*_{\nis}(X,\Omega_{\dR})\simeq H^*(X,\ste_{\dR})$$
is an isomorphism for any smooth $k$-scheme $X$.
\end{paragr}

\begin{paragr}
Suppose that $k$ is an algebraically closed field
of characteristic zero, complete with respect
to an archimedian (resp. non archimedian)
absolute value $| \, - \, |$. Then we can associate to
any smooth $k$-scheme $X$ an analytic space
(resp. a rigid analytic space) $X^\an$.
Let $\Omega^{*}_{X^\an}$ be the
analytic de Rham complex of $X^\an$
(seen as a sheaf of complexes).
This defines a presheaf $\Omega^\an_{\dR}$
of differential graded $k$-algebras
on $\Sm/k$ by the formula
$$\Omega^\an_{\dR}(X)=\Omega^{*}_{X^\an}(X^\an)\, .$$
The \emph{analytic de Rham cohomology} of a smooth scheme $X$
is defined as the hypercohomology of $X^\an$ with coefficients
in the sheaf $\Omega^{*}_{X^\an}$.
$$\HH^*_{\dR}(X^\an)=\hypercoh^*(X^\an,\Omega^{*}_{X^\an})$$
As $X^\an$ is Stein (resp. quasi-Stein) whenever $X$
is affine, Cartan's Theorem B (resp. Kiehl's analog of this theorem)
implies that for an affine smooth
$k$-scheme $X$, one has
$$\HH^*_{\dR}(X^\an)=\HH^*(\Omega^*_{X^\an}(X^\an))\, .$$
As analytic de Rham cohomology satisfies \'etale descent and is $\AA^1$-homotopy invariant, this implies that $\Omega^\an_{\dR}$ has the
B.-G.-property on affine smooth $k$-schemes, and is $\AA^1$-local. In fact, the complex $\Omega^\an_{\dR}$ is even a stable theory\footnote{We leave
this as an exercise for the reader; the arguments used below to prove that rigid cohomology is a stable theory (essentially the proof of Theorem
\ref{MWcohMWConSpecV}) might give a hint.} so that, by virtue of
Corollary \ref{compkperfectste}, the canonical map
$$\Omega^{}_{\dR}\To\Omega^{\an}_{\dR}$$
is a quasi-isomorphism locally for the Nisnevich topology.
In other words, we get Grothendieck's theorem~\cite{cohDR}
(resp. Kiehl's theorem~\cite{kieDR}): for any smooth $k$-scheme $X$,
the canonical map
$$H^*_{\dR}(X)\To H^*_{\dR}(X^\an)$$
is an isomorphism.
\end{paragr}

\subsection{Variations on Monsky-Washnitzer cohomology}

\begin{paragr}
We consider here a complete discrete valuation ring $V$ with fraction field of
characteristic zero $K$ and perfect residue field $k$. We set
$S=\spec {V}$, $\eta=\spec K$, and $s=\spec k$. We have an open immersion $j:\eta\To S$ and a closed immersion $i:s\To S$. For a (smooth) $S$-scheme
$X$, we write $X_{\eta}$ and $X_{s}$ for the generic fiber and the special fiber of $X$ respectively.
\end{paragr}

\begin{paragr}\label{defcohMW}
Consider a smooth affine $S$-scheme $X=\spec A$.

We denote by $A^\dagger$ the weak completion of $A$
with respect to the $\idealmax$-adic topology, where $\idealmax$ stands for
the maximal ideal of $V$; see \cite[Definition 1.1]{MW}.
Recall $A^\dagger$ is a formally smooth
$V$-algebra \cite[Theorem 2.6]{MW}. Denote by
$\Omega^*(A^\dagger/V)$ the complex of differential
forms of $A^\dagger$ relative to $V$.
It can be defined as the universal $\idealmax$-separated
differential graded $V$-algebra associated to $A$; see \cite[Theorem 4.2]{MW}.
More precisely, it is obtained from the algebraic de Rham complex of $A^\dagger$
over $V$ by the formula
$$\Omega^*(A^\dagger/V)=\Omega^*_{A^\dagger /V}
/\cap^\infty_{i=0}\idealmax^{i}\Omega^*_{A^\dagger /V} \, .$$
The \emph{Monsky-Washnitzer complex} of $X$ is defined as
$$E_{\mw}(X)=\Omega^*(A^\dagger/V)\otimes_{V}K
\simeq A^\dagger\otimes_A\Omega^*_{A/V}\otimes_V K\, ,$$
and the \emph{Monsky-Washnitzer cohomology} of $X$ is
$$H^n_\mw(X)=H^n(E_{\mw}(X))\, .$$
(see \cite{MW,put}).
\end{paragr}

\begin{thm}\label{MWcohMWConSpecV}
The Monsky-Washnitzer complex is a stable theory on
smooth affine $S$-schemes.
\end{thm}

\begin{proof}
The complex $E_\mw(X)$ can be compared with Berthelot's
rigid cohomology; see \cite[Proposition 1.10]{Ber1}. More precisely,
once a closed embedding $X\To\AA^n_S$ is chosen,
let $W$ denotes the schematic closure of $X$ in
$\PP^n_S$, and $\hat W$ denotes the formal
$\idealmax$-adic completion of $W$.
The proof of \cite[Proposition 1.10]{Ber1} consists then
to check that we have a canonical isomorphisms
of complexes of $K$-vector spaces
$$E_\mw(X)\simeq\varinjlim_V\Gamma(V,\Omega^*_{V})\simeq
\Gamma(]W[_{\hat W},j^\dagger\Omega^*_{\hat W})\, ,$$
where $V$ ranges over the strict neighbourhoods of
the tube of $X$ in $\hat W$, and that the canonical map
$$H^n(\Gamma(]W[_{\hat W},j^\dagger\Omega^*_{\hat W}))
\To H^n(\derR\Gamma(]W[_{\hat W},j^\dagger\Omega^*_{\hat W}))
=H^n_\rig(X_s/K)$$
is an isomorphism. In other words, $E_\mw(X)$ is (up to a canonical
quasi-isomorphism) the rigid complex associated to the
embbedings
$$X_s\To W \To \hat W\, .$$
Using \cite[Proposition 2.2]{Ber1}, one can extend this
comparison results to cohomology with support: for
any closed subscheme $Z$ of $X$, one has a canonical
isomorphism
$$H^n_{\mw,Z}(X)\simeq H^n_{\rig,Z_s}(X_s/K)$$
(where $H^{i+1}_{\mw,Z}(X)$ denotes the $i$th
cohomology group of the cone
of the map $E_\mw(X)\To E_\mw(X-Z)$).
Hence to prove \'etale excision for Monsky-Washnitzer
cohomology, we are reduce to prove \'etale excision
for rigid cohomology. This follows immediately
from the \'etale descent theorem for rigid cohomology,
proved by Chiarellotto and Tsuzuki~\cite{ChTs}.

We also have the following computations:
\begin{equation*} 
H^n_{\mw}(\AA^1_{S})=
\begin{cases}
K & \text{if } n=0 \\
0 & \text{otherwise}
\end{cases}\\
\qquad H^n_{\mw}(\GG_{m})=
\begin{cases}
K & \text{if } n=0 \\
K.\dlog & \text{if } n=1 \\
0 & \text{otherwise}
\end{cases}
\end{equation*}
(where $\dlog$ is the differential form on $V[t,t^{-1}]^\dagger$
defined by $\dlog(t)=dt/t$).

It remains to prove that the K\"unneth map
$$E_\mw(X)\otimes^{}_{K}E_\mw(Y)
\To E_\mw(X\times_{S}Y)$$
is a quasi-isomorphism
for any affine smooth $S$-scheme $X$ and for $Y=\AA^1_{S}$ or $Y=\GG_{m}$.
If $Y=\AA^1_S$, this follows from Monsky and Washnitzer
Homotopy Invariance Theorem \cite[Theorem 5.4]{MW}.
The case of $Y=\GG_{m}$ is solved by considering
 the Gysin long exact sequence associated
to the closed immersion
$$i:X=X\times\{0\}\To X\times\AA^1$$
which is constructed explicitely from \cite[Theorem 3.5]{MW2}:
$$\cdots\To
\HH^{n-2}_{\mw}(X)\To\HH^n_{\mw}(X\times\AA^1)
\To\HH^n_{\mw}(X\times\GG_{m})\To\HH^{n-1}_{\mw}(X)
\To\cdots
$$
The homotopy invariance of Monsky-Washnitzer cohomology
allows then to split canonically the long exact sequence above (using the
projection of $X\times\GG_{m}$ onto $X$), and we finally get
an isomorphism of graded $\HH^{*}_{\rig}(X/K)$-modules
$$\HH^*_{\mw}(X\times\GG_{m})
\simeq\HH^{*-1}_{\mw}(X)\, . \, \dlog \oplus\HH^*_{\mw}(X)\, .$$
This implies immediately the K\"unneth formula above for $Y=\GG_{m}$.
\end{proof}

\begin{paragr}
We define a presheaf of commutative
differential graded $K$-algebras $j_*E_{\dR}$ on $\Sm/S$
by the formula below.
\begin{equation*}
j_*E_{\dR}(X)=\Omega_{\dR}(X_{\eta})
\end{equation*}
It follows immediately from Proposition \ref{dRmixedWeil}
that $j_*E_\dR$ is a mixed Weil cohomology on affine smooth $S$-schemes.
\end{paragr}

\begin{paragr}\label{defspdRrig}
Consider a smooth affine $S$-scheme $X=\spec A$.
By definition of the Monsky-Washnitzer complex,
we have a natural morphism of differential graded algebras
\begin{equation}\label{defspdRrig2}
\specialisation_{X} : j_*E_{\dR}(X)=\Omega^*_{A/V}\otimes_VK
\To A^\dagger\otimes_A\Omega^*_{A/V}\otimes_VK=E_{\mw}(X)
\end{equation}
called the \emph{specialisation map}. This defines a morphism of
presheaves of differential graded algebras
\begin{equation}\label{defspdRrig3}
\specialisation : j_*E_{\dR}\To E_{\mw}\, .
\end{equation}
Denote by $j_*\ste_{\dR}$ (resp. by $\ste_{\mw}$) the commutative ring spectra associated to
$j_*E_{\dR}$ and $E_{\mw}$ respectively.
It is clear that we have, for any affine smooth $S$-scheme $X$, the following
identifications in the derived category of the category of $K$-vector spaces.
\begin{equation}\label{defspdRrig4}
\derR\Gamma(X,j_*\ste_{\dR})\simeq\Omega^*_\dR(X_\eta) \quad\text{and}
\quad \derR\Gamma(X,\ste_{\mw})\simeq E_\mw(X)\, .
\end{equation}
\end{paragr}

\begin{thm}\label{berthelotogus}
There is a specialisation map
$$\specialisation : j_*\ste_{\dR}\To \ste_{\mw}$$
in $\DMt(S,K)$ which is compatible with
cup product, and induces isomorphisms
$$\derR\Gamma(X_{\eta},\ste_{\dR})
\xrightarrow{\specialisation_{X}}\derR\Gamma(X,\ste_{\mw})
\quad\text{and}\quad
\derR\Gamma_{c}(X_{\eta},\ste_{\dR})
\xrightarrow{\specialisation_{X,c}}\derR\Gamma_{c}(X,\ste_{\mw})$$
in $\Der(K)$ for any smooth $S$-scheme $X$
such that $j_*\ste_{\dR}(X)$ has a strong dual in $\DMt(S,j_*\ste_{\dR})$
(e.g. $X$ might be projective
or the complement of a relative strict normal crossings divisor
in a smooth and projective $S$-scheme).
\end{thm}

\begin{proof}
Apply Theorem \ref{comparisonste} to \eqref{defspdRrig3} to get directly the map $\specialisation$ from $j_*\ste_{\dR}$ to $\ste_{\mw}$ and the
isomorphism $\specialisation_{X}$.

We also obtain isomorphisms

$$\derR\Gamma_{c}(X,j_*\ste_{\dR})
\simeq\derR\Gamma_{c}(X,\ste_{\mw})\, .$$

Using the fact $j_*\ste_\dR(X)$ has a strong dual in $\DMt(S,j_*\ste_\dR)$,
we have the following computations (we
assume $X/S$ is of dimension $d$).
$$\begin{aligned}
\derR\Gamma_{c}(X,j_*\ste_{\dR})
&\simeq\dual{\derR\Gamma(X,j_*\ste_{\dR}(-d)[-2d])}\\
&\simeq\dual{\derR\Gamma(X_{\eta},\ste_{\dR}(-d)[-2d])}\\
&\simeq\derR\Gamma_{c}(X_{\eta},\ste_{\dR})
\end{aligned}$$
These identifications give the expected isomorphism $\specialisation_{X,c}$.
\end{proof}

\begin{cor}\label{nonstrongduals}
For any non empty smooth $S$-scheme $X$ with empty special fiber, $\tQQ(X)$
does not have any strong dual in $\DMB(S)$.
\end{cor}

\begin{proof}
Given such an $S$-scheme $X$, it is clear that the specialisation map
$$\derR\Gamma(X_{\eta},\ste_{\dR})
\xrightarrow{\specialisation_{X}}\derR\Gamma(X,\ste_{\mw})=0$$ is not an isomorphism.
But if $\tQQ(X)$ had a strong dual in $\DMB(S)$,
then $j_*\ste_{\dR}(X)$ would have a strong dual in $\DMt(S,j_*\ste_{\dR})$ as well,
so that, by virtue of Theorem \ref{berthelotogus},
$\specialisation_{X}$ would be an isomorphism in $\Der(K)$.
\end{proof}

\begin{paragr}\label{defrigidonkwithV}
Recall from \cite{ayoub,roe,TCMM} that we have two
pairs of adjoint functors
\begin{equation}
j_{!}:\DMte(\eta,\QQ)\leftrightarrows\DMte(S,\QQ):j^*
\end{equation}
\begin{equation}\label{derivedistar}
\derL i^{*}:\DMte(S,\QQ)\leftrightarrows\DMte(s,\QQ):i_{*}
\end{equation}
such that $j_{!}$ and $i_{*}$ are fully faithful,
and such that for any object $M$ of $\DMte(S,K)$, there is a canonical
distinguished triangle:
\begin{equation}
j_{!}j^*(M)\To M\To i_{*}\derL i^*(M)\To j_{!}j^*(M)[1]\, .
\end{equation}
As we obviously have $j^*(E_{\mw})=0$ (this just means the Monsky-Washnitzer
cohomology of an affine smooth $V$-scheme with empty special fiber
is trivial), we deduce that
\begin{equation}
E_{\mw}\simeq i_{*} \derL i^*(E_{\mw})\, .
\end{equation}
Let $X$ be a smooth affine $k$-scheme. Using \cite[Theorem 1.3.1]{arabia},
there exists a smooth and affine $V$-scheme $Y=\spec A$
such that $X=Y_{s}$.
In other words, we get $\QQ(X)=\derL i^*\QQ(Y)$. This leads to the following
computations.
\begin{equation}\label{comprigmw1}\begin{split}
\begin{aligned}
\derR\Gamma(X,\derL i^* E_{\mw})
&\simeq \derR\Hom_{\QQ}(\QQ(X),\derL i^* E_{\mw})\\
&\simeq \derR\Hom_{\QQ}(\derL i^*\QQ(Y), \derL i^* E_{\mw})\\
&\simeq \derR\Hom_{\QQ}(\QQ(Y),i_{*}\derL i^*E_{\mw})\\
&\simeq \derR\Hom_{\QQ}(\QQ(Y),E_{\mw})\\
&\simeq \derR\Gamma(Y,E_{\mw})\\
&\simeq E_\mw(Y)
\end{aligned}\end{split}
\end{equation}
Note this isomorphism is functorial with respect to $Y$, $X$ being identified
with $Y_s$. 

The cohomology theory represented by $\derL i^* E_\mw$
in $\DMte(s,\QQ)$ can be described as a stable cohomology
theory as follows. The main difficulty for this is to
represent it by a sheaf of commutative differential graded $K$-algebras.
This is achieved by having a closer look at the definition of
the functor $\derL i^*$ of \eqref{derivedistar}:
this is the total left derived functor of the functor
\begin{equation}\label{defrifiminv2}
i^*:\Comp(\mathit{Sh}(\Sm/S,\QQ))\To\Comp(\mathit{Sh}(\Sm/s,\QQ))
\end{equation}
which preserves colimits and sends $\QQ(X)$ to $\QQ(X_{s})$.
The functor \eqref{defrifiminv2} is a left Quillen functor
with respect to the model structures defined by Proposition \ref{A1cmfeff}.
Hence $\derL i^*E_{\mw}$ is defined by applying \eqref{defrifiminv2}
to a $\V$-cofibrant resolution of $E_{\mw}$, where
$\V$ is the category of smooth and affine $V$-schemes.
We can consider a quasi-isomorphism
$p:E'_{\mw}\To E_{\mw}$,
with $E'_{\mw}$ a commutative monoid
which is $\V$-cofibrant as a complex of sheaves
(using the model structure of \cite[proposition 4.3.21]{DAG3},
whose assumptions are trivially checked in the $\QQ$-linear setting).
We then put $E_{\rig}=i^*E'_{\mw}$.
By definition, we have a canonical isomorphism
$$E_{\rig}\simeq \derL i^*E_{\mw}\, .$$

We will call $E_{\rig}$ the \emph{rigid cohomology complex}.
By construction, for any smooth affine $k$-scheme $X$, we have an isomorphism
\begin{equation}\label{defrig01}
\derR\Hom(\QQ(X),\derL i^* E_\mw)\simeq E_{\rig}(X)\, .
\end{equation}
\end{paragr}

\begin{prop}\label{rigeqrig}
For any smooth affine $k$-scheme $X$, there is a functorial isomorphism
$$H^n(E_{\rig}(X))\simeq H^n_\rig(X/K)\, ,$$
where $H^n_\rig(X/K)$ denotes Berthelot's rigid cohomology of $X$.
This comparison map is compatible with cup product.
\end{prop}

\begin{proof}
We know that, given a smooth affine $S$-scheme $Y$, there is
a functorial isomorphism
$$H^n_\mw(Y)\simeq H^n_\rig(Y_s/K)\, ,$$
which is compatible with cup product; see \cite[Proposition 1.10]{Ber1}.
As, by definition, $E_\rig=\derL i^*E_\mw$,
given a smooth affine $k$-scheme $X$, once a smooth affine
$S$-scheme $Y$ with special fiber isomorphic to $X$ is chosen,
we obtain from the isomorphims \eqref{comprigmw1} that
$$H^n(E_{\rig}(X))\simeq H^n_\rig(X/K)\, .$$
It remains to prove this isomorphism is independent of the lift $Y$, and
is functorial in $X$. Let $f:X\To X'$ be a morphism of smooth affine $k$-schemes.
Choose two smooth affine $S$-schemes $Y$ and $Y'$ endowed with
two isomorphisms $Y_s\simeq X$ and $Y'_s\simeq X'$ (which exist, thanks to
\cite[Theorem 1.3.1]{arabia}).
By virtue of \cite[Theorem 2.1.3]{arabia}, there exists
a commutative diagram of $S$-schemes
$$\xymatrix{
&X\ar[r]^f\ar[d]^{i_\varepsilon}\ar[dl]_{i}&X'\ar[d]^{i'}\\
Y&Y_\varepsilon\ar[l]^\varepsilon\ar[r]_g&Y'
}$$
with $i$ (resp. $i_\varepsilon$, resp. $i'$) being a closed
immersion which identifies $X$ (resp. $X$, resp. $X'$) with the
special fiber of $Y$ (resp. of $Y_\epsilon$, resp. of $Y'$), and
with $\varepsilon: Y_\varepsilon\To Y$ \'etale and inducing the
identity on the special fibers.
Then the naturality of the isomorphisms \eqref{comprigmw1}
gives the following commutative diagram
$$\xymatrix{
&E_\rig(X)&E_\rig(X')\ar[l]_{f^*}\\
E_\mw(Y)\ar[ur]^\simeq\ar[r]_{\varepsilon^*}&E_\mw(Y_\varepsilon)\ar[u]_\simeq&
E_\mw(Y')\ar[l]^{g^*}\ar[u]_\simeq
}$$
in which the non-horizontal maps are the canonical isomorphisms.
\end{proof}

\begin{thm}\label{rigidMWC}
The sheaf of commutative differential graded algebras
$E_{\rig}$ is a mixed Weil cohomology on
smooth $k$-schemes.
\end{thm}

\begin{proof}
As $E_{\rig}$ is fibrant by definition, it is $\AA^1$-homotopy
invariant and has the B.-G. property. Using Theorem \ref{MWcohMWConSpecV}
and the comparison isomorphisms \eqref{comprigmw1}, we see
that $E_{\rig}$ is a stable cohomology theory.
It thus remains to prove the K\"unneth Formula.
This comes immediately from the comparison
with Berthelot's rigid cohomology (Proposition \ref{rigeqrig}),
the latter being known to satisfy the K\"unneth formula; see \cite{Ber2}.
\end{proof}

\begin{sch}
Let us denote by $\ste_{\rig}$ the commutative ring spectrum
associated to the mixed Weil cohomology $E_{\rig}$.

Theorem \ref{berthelotogus} can be made a little more precise
in the following way. Recall from \cite{roe,ayoub,TCMM}
that we have a pair of adjoint triangulated functors
\begin{equation}
\derL i^* : \DMt(S,\QQ)\leftrightarrows\DMt(s,\QQ) : i_{*}
\end{equation}
satisfying the following properties.
\begin{itemize}
\item[(i)] The functor $\derL i^*$ is symmetric monoidal and preserves Tate twists.
\item[(ii)] For any smooth $S$-scheme $X$, we have $\derL i^*\tQQ(X)=\tQQ(X_{s})$.
\item[(iii)] The functor $i_{*}$ is fully faithful.
\item[(iv)] For any objects $M$ of $\DMt(S,\QQ)$ and any object $N$
of $\DMt(s,\QQ)$, we have a canonical isomorphism
\begin{equation}
M\otimes^\derL_{\QQ} i_{*}(N)\simeq i_{*}(\derL i^*(M)\otimes^\derL_{\QQ}N)\, .
\end{equation}
\end{itemize}
It follows from property (ii) and the definition of $\ste_{\rig}$ that we have
an isomorphism
\begin{equation}
\ste_{\mw}\simeq i_{*}\ste_{\rig}\, ,
\end{equation}
so that we have a specialization map
\begin{equation}\label{goodspecial}
\specialisation : j_*\ste_{\dR}\To i_*\ste_{\rig}
\end{equation}
in $\DMt(S,K)$.
Moreover, we obtain from properties (i) and (iv)
the following identifications for a smooth $S$-scheme $X$
of pure dimension $d$.
$$\begin{aligned}
\derR\Gamma_{c}(X,\ste_{\mw})
&\simeq\derR\Hom_{\QQ}(\tQQ,\ste_{\mw}\otimes^\derL_{\QQ}\tQQ(X)(-d)[-2d])\\
&\simeq\derR\Hom_{\QQ}(\tQQ,i_{*}(\ste_{\rig})\otimes^\derL_{\QQ}\tQQ(X)(-d)[-2d])\\
&\simeq\derR\Hom_{\QQ}(\tQQ,i_{*}(\ste_{\rig}\otimes^\derL_{\QQ}i^{*}(\tQQ(X)))(-d)[-2d])\\
&\simeq\derR\Hom_{\QQ}(i^*(\tQQ),\ste_{\rig}\otimes^\derL_{\QQ}\tQQ(X_{s})(-d)[-2d])\\
&\simeq\derR\Hom_{\QQ}(\tQQ,\ste_{\rig}\otimes^\derL_{\QQ}\tQQ(X_{s})(-d)[-2d])\\
&\simeq\derR\Gamma_{c}(X_{s},\ste_{\rig})
\end{aligned}$$
By virtue of Theorem \ref{berthelotogus},
the specialisation map \eqref{goodspecial} induces isomorphisms
$$\derR\Gamma(X_{\eta},\ste_{\dR})
\xrightarrow{\specialisation_{X}}\derR\Gamma(X_s,\ste_{\rig})
\quad\text{and}\quad
\derR\Gamma_{c}(X_{\eta},\ste_{\dR})
\xrightarrow{\specialisation_{X_s,c}}\derR\Gamma_{c}(X_s,\ste_{\rig})$$
in $\Der(K)$ for any smooth $S$-scheme $X$
such that $j_*\ste_{\dR}(X)$ has a strong dual in $\DMt(S,j_*\ste_{\dR})$.

It can be proven that $E_\rig$ is quasi-isomorphic to the restriction
of Besser's rigid complex (see \cite[Definition 4.13]{Bes}) to the
category of smooth $k$-schemes. In other words, the object
$\ste_\rig$ represents Berthelot's rigid cohomology in $\DMt(s,\QQ)$.
In the case where $X$ is smooth and projective over $S$, using the
comparison isomorphism relating rigid cohomology and crystalline cohomology
(see \cite[Proposition 1.9]{Ber1}),
Theorem \ref{berthelotogus} gives back the comparison isomorphism of
Berthelot and Ogus~\cite{BerOg}.
\end{sch}

\subsection{\'Etale cohomology}

\begin{paragr}
For sake of completeness, we will finish by explaining how
$\ell$-adic cohomology fits in the picture of
mixed Weil cohomologies as they are defined here.

Consider a countable perfect field $k$, and choose a separable closure $\overline k$ of $k$. For a smooth $k$-scheme $X$, write $\overline
X=X\otimes_{k}\overline k$. Let $\ell$ be a prime which is distinct from the characteristic of $k$.

Deligne~\cite{WeilII} defines for any smooth $k$-scheme $X$
a commutative differential graded $\QQ_{\ell}$-algebra
which computes the $\ell$-adic cohomology of \smash{$\overline X$}.
We will modify slightly some steps of his construction to ensure its functoriality.
\end{paragr}

\begin{paragr}
Consider a pro-simplicial set $\mathbf X=\sideset{``}{\text{''}}\limproj X_{\alpha}$.
We can then define its singular cohomology with coefficients in $\ZZ/\ell^n$
by the formula
\begin{equation}
H^i(\mathbf X,\ZZ/\ell^n)=\limind H^i(X_{\alpha},\ZZ/\ell^n)\, .
\end{equation}
We will say that $\mathbf X$ is \emph{essentially $\ell$-finite}
if the groups $H^i(\mathbf X,\ZZ/\ell^n)$ are finite.

For an essentially $\ell$-finite pro-simplicial set $\mathbf X$,
formula (5.2.1.7) of \cite{WeilII} defines a commutative
differential graded $\QQ_{\ell}$-algebra $A(\mathbf X)$ such that
\begin{equation}
H^{i}(A(\mathbf X))=\QQ_{\ell}\otimes_{\ZZ_{\ell}}\limproj_{n}H^{i}(\mathbf X,\ZZ/\ell^n)\, .
\end{equation}
This construction is (contravariantly) functorial in $\mathbf X$.
\end{paragr}

\begin{paragr}
For an \'etale surjective morphism $X'\To X$, define $\check C(X'/X)$
to be the \v Cech simplicial scheme defined by the formula
$$\check C(X'/X)_{n}=\underbrace{X'\times_{X}\dots\times_{X}X'}_{{\text{$n+1$ times}}} \, .$$
Note that the map $\check C(X'/X)\To X$ is an  \'etale hypercovering.

Given a smooth $\overline k$-scheme $X$, define
an \emph{\'etale fundamental system} $\mathcal X$
of $X$ to be a tower of morphisms of smooth $\overline k$-schemes
indexed by integers $\alpha\geq 0$
$$
\mathcal X=
\Big[\cdots\To X_{\alpha+1}\To X_{\alpha}\To\cdots\To X_{1}\To X_{0}=X\Big]
$$
such that $X_{\alpha}\To X$ is \'etale surjective for all $\alpha\geq 0$,
and such that any \'etale surjective map $U\To X$ factors through $X_{\alpha}$
for $\alpha$ big enough. Such an \'etale fundamental system of $X$
defines a pro-simplicial scheme $\sideset{``}{\text{''}}\limproj \check C(X_{\alpha}/X)$,
whence a pro-simplicial set
\begin{equation}
\pi(\mathcal X)=\sideset{``}{\text{''}}\limproj \pi_{0}(\check C(X_{\alpha}/X))
\end{equation}
which is essentially $\ell$-finite, and such that there is a canonical isomorphism
\begin{equation}\label{coincideetale1}
H^{i}(A(\pi(\mathcal X)))\xrightarrow{\sim}H^{i}_{\etale}(X,\QQ_{\ell})
\end{equation}
(see \cite[5.2.2]{WeilII}). Given a non-empty finite family of \'etale
fundamental systems $\mathcal X=\{\mathcal X^{1},\ldots,\mathcal X^{n}\}$, we define an
\'etale fundamental system $\mathcal X_{\mathit{tot}}$ whose $\alpha^{\text{th}}$
stage is defined as the fiber product of the $\alpha^{\text{th}}$
stages of the $\mathcal X^{i}$'s over $X$. Given a non-empty subset $\mathcal X'$
of $\mathcal X$, it can be described as
$\mathcal X'=\{\mathcal X^{i_1},\ldots,\mathcal X^{i_m}\}$, with
$i_k\neq i_l$ whenever $k\neq l$, with $m\leq n$.
We then have a canonical morphism of pro-schemes
$\mathcal X_{\mathit{tot}}\To\mathcal X'_{\mathit{tot}}$
induced by the projections $\mathcal X_{\mathit{tot}}\To \mathcal X^{i_k}$.
Taking the filtering projective limit
of all the pro-simplicial sets $\pi(\mathcal X_{\mathit{tot}})$,
where $\mathcal X$ ranges over the non-empty finite families
of \'etale fundamental systems of $X$,
defines a pro-simplicial set. We define
\begin{equation}\label{defcohadgetale}
A(X)=\limind_{\mathcal X} A(\pi(\mathcal X_{\mathit{tot}}))\, .
\end{equation}
As filtering colimits are exacts,
we deduce from \eqref{coincideetale1} that we have a canonical isomorphism
\begin{equation}\label{coincideetale2}
H^{i}(A(X))=\limind_{\mathcal X} H^{i}(A(\pi(\mathcal X_{\mathit{tot}})))
\xrightarrow{\sim}H^{i}_{\etale}(X,\QQ_{\ell})\, .
\end{equation}
We claim formula \eqref{defcohadgetale} defines a presheaf of commutative differential graded $\QQ_{\ell}$-algebras $A$ on the category of smooth
$\overline k$-schemes. Consider a morphism $f:X\To Y$ of smooth $\overline k$-schemes.
Any non-empty finite family of \'etale
fundamental systems $\mathcal Y=\{\mathcal Y^{1},\ldots,\mathcal Y^{n}\}$ of $Y$
defines by pullback a non-empty finite family of \'etale
fundamental systems $f^*(\mathcal Y)=\{X\times_Y\mathcal Y^{1},\ldots,X\times_Y\mathcal Y^{n}\}$
of $X$, with a canonical morphism of pro-schemes
$f^*(\mathcal Y)_{\mathit{tot}}\To\mathcal Y_{\mathit{tot}}$.
This induces a map
$$A(\pi(\mathcal Y_{\mathit{tot}}))\To A(\pi(f^*(\mathcal Y)_{\mathit{tot}}))\To A(X)\, .$$
By passing to the colimit of the $A(\pi(\mathcal Y_{\mathit{tot}}))$'s,
we get the expected map
$$f^* : A(Y)\To A(X)\, .$$
\end{paragr}

\begin{paragr}
Define a presheaf of commutative differential graded $\QQ_{\ell}$-algebras
$E_{\etale,\ell}$ on $\Sm/k$ by the formula
\begin{equation}
E_{\etale,\ell}(X)=A(\overline X)\, .
\end{equation}
Then one has
\begin{equation}
H^{n}(E_{\etale,\ell}(X))=\HH^n_{\etale}(\overline X,\QQ_{\ell})\, .
\end{equation}
In particular, $E_{\etale,\ell}$ satisfies \'etale descent, whence
it has the B.-G.-property. The well known properties of
\'etale cohomology proved by Artin and Grothendieck thus imply:
\end{paragr}

\begin{thm}
$E_{\etale,\ell}$ is a mixed Weil theory over smooth $k$-schemes.
\end{thm}


\bibliographystyle{amsalpha}
\bibliography{common}

\end{document}